%
%

%
%




\newcount\chpno
\newcount\secno
\newcount\thmno
\newcount\eqnno


\magnification=1200
\hsize=6.5truein
\vsize=9truein
\hoffset=0.05truein
\voffset=-0.25truein
\nopagenumbers
\raggedbottom


\def\makeheadline{\vbox to 0pt{\vskip -.5 truein
  \line{\vbox to 8.5 truept{}\the\headline}\vss}
  \nointerlineskip}

\def\uppernumberson{\headline={\hfil\rlap{\hbox to 1 truein
                         {\tenrm\hfil\folio\hskip 0.5 truein}}}}


\font\runheadfont=cmr10 at 10 true pt
\def\rhead{{\runheadfont GEOMETRY OF MODULI SPACE ENDS}}
\def\lhead{{\runheadfont PAUL M. N. FEEHAN}}

\def\rightheadline{\tenrm\hfil\ifnum\pageno<3{}\else\rhead\fi\hfil\folio}
\def\leftheadline{\tenrm\folio\hfil\lhead\hfil}

\def\numberson{\headline={\ifodd\pageno\rightheadline\else\leftheadline\fi}}
\voffset=2\baselineskip


\def\singlespace{
   \normallineskip=1pt
   \normalbaselineskip=12pt
   \normallineskiplimit=0pt
   \normalbaselines
   \abovedisplayskip=12pt plus 3pt minus 9pt
   \abovedisplayshortskip=0pt plus 3pt
   \belowdisplayskip=12pt plus 3pt minus 9pt
   \belowdisplayshortskip=7pt plus 3pt minus 4pt
   \parskip=0pt plus 1pt
   \smallskipamount=3pt plus 1pt minus 1pt
   \medskipamount=6pt plus 2pt minus 2pt
   \bigskipamount=12pt plus 4pt minus 4pt}







\outer\def\chapter#1{\global\advance\chpno by 1
          \centerline{\bf\the\chpno.  #1}
          \nobreak\bigskip\hskip 0pt}

\def\endchapter{\bigskip\bigskip}


\def\section#1{{\bigbreak\global\advance\secno by 1
      \noindent{\bf\the\chpno.\the\secno. #1.}\nobreak}}
 


\def\subsection#1{{\medbreak\noindent{\bf #1.}\enspace\nobreak}}
\def\subsubsection#1{{\medbreak\noindent{\it #1.}\enspace\nobreak}}



\def\assm{\medbreak\noindent\global\advance\thmno by 1 
 {\bf Assumption \the\chpno.\the\thmno.\enspace}}

\def\case#1{\smallbreak\noindent{\bf Case #1\enspace}}

\def\cond{\medbreak\noindent\global\advance\thmno by 1 
 {\bf Condition \the\chpno.\the\thmno.\enspace}}

\def\data{\medbreak\noindent\global\advance\thmno by 1 
 {\bf Data \the\chpno.\the\thmno.\enspace}}

\def\defn{\medbreak\noindent\global\advance\thmno by 1
 {\bf Definition \the\chpno.\the\thmno.\enspace}}

\def\expl{\medbreak\noindent\global\advance\thmno by 1
 {\bf Example \the\chpno.\the\thmno.\enspace}}

\def\step#1{\smallbreak\noindent{\bf Step #1\enspace}}

\def\rmk{\medbreak\noindent\global\advance\thmno by 1
{\bf Remark \the\chpno.\the\thmno.\enspace}}

\def\pf{\noindent{\it Proof}.\enspace}
\def\qed{\hfill{$\square$}\medskip}

\outer\def\proclaim#1#2
 {\global\advance\thmno by 1
 \medbreak\noindent
 {\bf#1 \the\chpno.\the\thmno.\enspace}{\sl#2}\par
 \ifdim\lastskip<\medskipamount\removelastskip\penalty55\medskip\fi}


\def\eqlabel{\global\advance\eqnno by 1
             \leqno{(\the\chpno.\the\eqnno)}}


\def\numeqn{\global\advance\eqnno by 1
            (\the\chpno.\the\eqnno)}


\def\chplabel#1{\global\edef#1{\the\chpno}}
\def\eqnlabel#1{\global\edef#1{\the\chpno.\the\eqnno}}
\def\seclabel#1{\global\edef#1{\the\chpno.\the\secno}}
\def\thmlabel#1{\global\edef#1{\the\chpno.\the\thmno}}


\def\dummy{??}
\def\ref#1{\ifx #1\undefined\let#1=\dummy #1\else #1\fi}


\def\eqref#1{\ifx #1\undefined\let#1=\dummy #1\else (#1)\fi}




\input amssym.def
\input amssym.tex

\font\teneusm=eusm10  
\font\seveneusm=eusm7  
\font\fiveeusm=eusm5  
\newfam\eusmfam
\textfont\eusmfam=\teneusm  
\scriptfont\eusmfam=\seveneusm
\scriptscriptfont\eusmfam=\fiveeusm
\def\script#1{{\fam\eusmfam\relax#1}}

\font\sevenex=cmex7
\scriptfont3=\sevenex
\scriptscriptfont3=\sevenex


\font\tencmmib=cmmib10  \skewchar\tencmmib='177
\font\sevencmmib=cmmib7  \skewchar\sevencmmib='177
\font\fivecmmib=cmmib5  \skewchar\fivecmmib='177
\newfam\cmmibfam
\textfont\cmmibfam=\tencmmib  
\scriptfont\cmmibfam=\sevencmmib
\scriptscriptfont\cmmibfam=\fivecmmib

\font\tencmbsy=cmbsy10  \skewchar\tencmbsy='60
\font\sevencmbsy=cmbsy7  \skewchar\sevencmbsy='60
\font\fivecmbsy=cmbsy5  \skewchar\fivecmbsy='60
\newfam\cmbsyfam
\textfont\cmbsyfam=\tencmbsy  
\scriptfont\cmbsyfam=\sevencmbsy
\scriptscriptfont\cmbsyfam=\fivecmbsy




\def\Biglongrightarrow{\hbox to 20pt{\rightarrowfill}}
\def\mapright#1{\smash{
    \mathop{\ \Biglongrightarrow\ }\limits^{#1}}}



\def\RR{{\Bbb R}}
\def\SS{{\Bbb S}}
\def\ZZ{{\Bbb Z}}



\def\cA{{\cal A}}
\def\cB{{\cal B}}

\def\cG{{\cal G}}

\def\today{\ifcase\month\or
  January\or February\or March\or April\or May\or June\or
  July\or August\or September\or October\or November\or
  December\fi
  \space\number\day, \number\year}


\def\deriv#1#2{{{d#1}\over{d#2}}}

\def\pd#1#2{{{\partial#1}\over{\partial#2}}}

\def\rd{\partial}
\def\]{\nabla}
\def\cov{\nabla}


\def\half{{\textstyle {1\over 2}}}
\def\quarter{{\textstyle {1\over 4}}}


\def\a{\alpha}
\def\al{\alpha}
\def\b{\beta}
\def\be{\beta}
\def\c{\gamma}
\def\ga{\gamma}
\def\C{\Gamma}
\def\Ga{\Gamma}
\def\d{\delta}
\def\de{\delta}
\def\rd{\partial}
\def\D{\Delta}
\def\De{\Delta}
\def\e{\varepsilon}
\def\eps{\varepsilon}
\def\z{\zeta}

\def\vth{\vartheta}

\def\k{\kappa}
\def\ka{\kappa}
\def\l{\lambda}
\def\la{\lambda}

\def\r{\rho}

\def\s{\sigma}
\def\si{\sigma}

\def\w{\omega}
\def\W{\Omega}
\def\om{\omega}
\def\Om{\Omega}



\def\SO{{\rm SO}}

\def\SU{{\rm SU}}


\def\less{\setminus}
\def\ol{\overline}

\def\too{\ \longrightarrow\ }

\def\8{\infty}
\def\<{\langle}
\def\>{\rangle}
\def\sumprime_#1{\setbox0=\hbox{$\scriptstyle{#1}$}
  \setbox2=\hbox{$\displaystyle{\sum}$}
  \setbox4=\hbox{${}'\mathsurround=0pt$}
  \dimen0=.5\wd0 \advance\dimen0 by-.5\wd2
  \ifdim\dimen0>0pt
  \ifdim\dimen0>\wd4 \kern\wd4 \else\kern\dimen0\fi\fi
 \mathop{{\sum}'}_{\kern-\wd4 #1}}


\def\cp{{\Bbb C}{\Bbb P}}
\def\hp{{\Bbb H}{\Bbb P}}


\def\ad{{\rm ad}\thinspace}
\def\Ad{{\rm Ad}\thinspace}

\def\Aut{{\rm Aut} \thinspace}
\def\Centre{{\rm Centre}\thinspace}
\def\Conf{{\rm Conf}}
\def\coker{{\rm coker}\thinspace}

\def\dist{{\rm dist}\thinspace}

\def\Gl{{\rm Gl}}

\def\Hom{{\rm Hom}}

\def\im{{\rm im}\thinspace}

\def\ker{{\rm ker}\thinspace}

\def\loc{{\rm loc}}

\def\Scale{{\rm Scale}}

\def\supp{{\rm supp\thinspace}}
\def\tr{{\rm tr}}



\def\sI{\script{I}}
\def\sJ{\script{J}}
\def\sL{\script{L}}

\def\sT{\script{T}}
\def\sU{\script{U}}
\def\sV{\script{V}}



\def\rad{{\bf r}}
\def\trans{{\bf p}}


\singlespace
\numberson
\pageno=1
\chpno=0
\eqnno=0
\secno=0
\thmno=0


\centerline{\bf GEOMETRY OF THE ENDS OF THE MODULI SPACE OF}
\centerline{\bf ANTI-SELF-DUAL CONNECTIONS}
\medskip
\centerline{Paul M. N. Feehan}
\footnote{}{Research was supported in part by a National Science Foundation Mathematical 
Sciences Postdoctoral Fellowship under grant DMS 9306061.}
\footnote{}{{\it Version:\/} September 26, 1994. }
\footnote{}{{\it Journal reference:\/} Journal of Differential Geometry {\bf 42} (1995), 465--553.}
\medskip
\chapter{Introduction}
Let $X_0$ be a closed, oriented, $C^\8$
four-manifold and let $M_{X_0,P}(g_0)$ be the moduli space of
$g_0$-anti-self-dual  
connections on a principal $G$ bundle $P$ over $X_0$. The subspace
$M_{X_0,P}^*(g_0)$, obtained by excluding the reducible connections is then
a finite-dimensional, usually non-compact, $C^\8$ manifold. 
The moduli space $M_{X_0,P}^*(g_0)$ is naturally endowed with 
a metric ${\bf g}$ of Weil-Petersson type, called the $L^2$ metric, and our
purpose in this article is to study the geometry of the moduli space ends. 

\subsection{(a) Main results} It has been conjectured by D. Groisser and T.
Parker in [G-P87], [G-P89] and by S. K. Donaldson in [D90a] that the moduli
space of anti-self-dual connections, endowed with the $L^2$ metric, has
finite volume and diameter. The goal of this article
is to prove this conjecture under the hypotheses described below.

\proclaim{Theorem}{\thmlabel{\thmFinVolDiam}
Let $X_0$ be a closed, connected, oriented, simply-connected, $C^\8$
four-manifold with generic metric $g_0$ and let $P$ be
a principal $G$ bundle over $X_0$ such that either (1) $G=\SU(2)$ or
$\SO(3)$ and $b^+(X_0)=0$, or 
(2) $G=\SO(3)$ and $w_2(P)\ne 0$, where $w_2(P)$ is the second
Stiefel-Whitney class of $P$. 
Then the moduli space $M_{X_0,P}^*(g_0)$ of irreducible
$g_0$-anti-self-dual connections on $P$ has finite volume and diameter with
respect to the $L^2$ metric ${\bf g}$ defined by $g_0$.}  

We plan to discuss the case of $G=\SU(2)$ and $b^+(X_0)>0$
in a subsequent article. 
Note that when $G=\SO(3)$ and $w_2(P)\ne 0$, the trivial
(product) connection $\Theta$ does not appear in the Uhlenbeck
compactification $\ol{M}_{X,P}(g_0)$.
By `diameter' we mean the sum of the diameters of the connected
components of  
$M_{X_0,P}^*(g_0)$; the hypotheses imply that $M_{X_0,P}^*(g_0)$ has
finitely many path components. In [D89] Donaldson conjectured
that the $L^2$-metric completion of the moduli space coincides with the
Uhlenbeck compactification [D86], [D-K]. We announce here the following
result whose proof is included in [F94].  

\proclaim{Theorem}{\thmlabel{\thmMetCompletion} Under the hypotheses of
{\rm Theorem \ref{\thmFinVolDiam}}, the completion of $M_{X_0,P}^*(g_0)$
with respect to the $L^2$ metric ${\bf g}$ is homeomorphic to the Uhlenbeck
compactification $\ol{M}^u_{X_0,P}(g_0)$.}    

The requirement that $X_0$ be simply-connected implies that the moduli
space of flat connections consists of a single point representing the
product connection over $X_0$. This assumption simplifies the description
of the ends of the moduli spaces 
$M_{X_0,P}^*(g_0)$, but is not important in the derivation of bounds for the
components of ${\bf g}$. We assume $G=\SU(2)$ or $\SO(3)$ in order to
appeal to the generic metric theorems of Freed and Uhlenbeck which ensure
that the moduli space is a $C^\8$ manifold: otherwise, the
bounds for ${\bf g}$ obtained in Chapter 5 hold for any compact Lie group. 
For the sake of clarity, we assume $G=\SU(2)$ for the remainder of the
article and denote $M_{X_0,P}(g_0)$ by $M_{X_0,k}(g_0)$, where $c_2(P)=k\ge
0$ is the second Chern class.
 
\subsection{(b) History} The properties of the $L^2$ metric have been
investigated by many authors in recent years, but most extensively by
Groisser and Parker. In particular, they have conducted detailed studies of
its behaviour at the boundary of certain $k=1$ moduli spaces. Explicit
formulas for the  
components of ${\bf g}$ have been found by Doi, Matsumoto, and Matumoto
[D-M-M], Groisser and Parker [G-P87], and Habermann [Hab] when $k=1$ and
$X_0$ is the four-sphere ${\Bbb S}^4$ with its standard round metric $g_1$.
Groisser conducted a similar study when $X_0$ is the complex projective
space 
$\ol{\cp}^2$, equipped with the Fubini-Study metric $g_{\rm FS}$ [G90]. 
Their formulas imply that these $k=1$ moduli spaces have finite ${\bf
g}$-volume and 
${\bf g}$-diameter. More generally, Groisser and Parker have established
Theorem \ref{\thmFinVolDiam} in the special case $k=1$ [G-P89]. They also
obtained $C^0$ bounds for ${\bf g}$ in neighbourhoods of the reducible
connections, the `conical ends', for any $k\ge 1$. In [G92], Groisser
refined some of the $k=1$ results obtained in [G-P89]. It is worth
recalling that the $L^2$ metric is not invariant with respect to
conformal changes in the metric $g_0$ on $X_0$.

The approach of [G-P89] does not appear to readily generalise to the case
$k>1$, since their method relies on Donaldson's collar map which gives a
diffeomorphism from the `bubbling end' of $M_{X_0,1}^*(g_0)$ to the collar
$X_0\times (0,\la_0)$. For this reason we adopt a quite different method
which uses the gluing 
techniques of Taubes and Donaldson to construct a system of local
coordinate charts covering the `ends' of the moduli space. We then estimate the components of ${\bf g}$ with respect to these coordinates. 
In the case of the Weil-Petersson metric on Teichm\"uller space, estimates
of this type have been obtained by Masur [Mas]. In [F92], the author
proved Theorem \ref{\thmFinVolDiam}, when $X_0=\SS^4$ and $k=2$, using
the ADHM correspondence [D-K]. After the present work was submitted,
a preprint was received from Peng giving $L^2$ estimates
for the derivatives with respect to moduli parameters
of the family of anti-self-dual connections $A$ on the
connected sum $X_0\#_\la\SS^4$ constructed in \S7.2.2 of [D-K],
with $H^2_{A_i}=0$ [Pe]. His $L^2$ estimates are defined
with respect to a family of metrics $g_\la$ 
which are conformally equivalent to $g_0$ and 
pinch the neck of the connected sum as
$\la\to 0$, away from the neck coinciding with $g_0$ on $X_0$ and 
converging in $C^2$ to the standard round metric on the unit sphere $\SS^4$.

\subsection{(c) Outline and strategy} It remains to summarise the methods
used in the proofs of our main results. Let us first recall the definition
of the $L^2$ metric. 
The tangent space $T_A M_{X_0,k}^*(g_0)$
is identified with the cohomology group $H_A^1=\ker d_A^{+,g_0}/ 
\im d_A^{*,g_0}$. Given tangent vectors $[a]$, $[b]$, the $L^2$ metric
${\bf g}$ is defined by 
$$
{\bf g}_{[A]}([a],[b])=(\pi_Aa,\pi_Ab)_{L^2(X_0,g_0)}, \eqlabel
$$
where $\pi_A=1-d_A(d_A^{*,g_0}d_A)^{-1}d_A^{*,g_0}$ is the $L^2$ orthogonal
projection from 
$L^2\Om^1(X_0,\ad P)$ to the subspace $\ker d_A^{*,g_0}$. Clearly, ${\bf
g}([a],[b])$ is bounded above by $\|a\|_{L^2}\|b\|_{L^2}$, and so a
reasonable strategy is to seek upper bounds for ${\bf g}$ over the moduli
space ends. This will suffice for our present application. 

\subsubsection{(i) Moduli space ends and the bubble tree compactification}
Our first task is to describe useful models for the ends of the moduli
space of anti-self-dual connections.  
Let $(A_0,x_1,\dots,x_{m_0})$ be a point in the stratum $\ol{M}^u_{X_0,k}(g_0)\cap (M_{X_0,k_0}(g_0)\times s^{k-k_0}(X_0))$ of the Uhlenbeck compactification 
(see \S 4.1) which lies away from the diagonals of the symmetric product,
so that $m_0=k-k_0$ and each point $x_i$ has multiplicity $1$. 
Then every point $[A]\in M_{X_0,k}(g_0)$ which is close enough to
$(A_0,x_1,\dots,x_{m_0})$ in the Uhlenbeck topology can be shown to lie in
a neighbourhood constructible by gluing or `gluing neighbourhood' [D86],
[D-K]. Thus, suppose $[A_\al]$ is a sequence in 
$M_{X_0,k}(g_0)$ which converges weakly to $(A_0,x_1,\dots,x_{m_0})$. As
described in \S 4.2, the sequence of connections $[A_\al]$ produces
sequences of local mass centres $x_{i\al}$ converging to the points $x_i$
and sequences of local scales $\la_{i\al}$ converging to zero. Using the
scales $\la_{i\al}$, one now dilates the metric $g_0$ around the points
$x_{i\al}$ and produces a sequence of conformally equivalent, $C^\8$
metrics $g_\al$ on a connected sum $X\equiv X_0\#_{i=1}^{m_0}\SS^4$. As the
scales $\la_{i\al}$ tend to zero, the corresponding neck is pinched and the
connected-sum metrics $g_\al$ converge in $C^\8$ on compact subsets away
from the neck regions to the metric $g_0$ on $X_0$ and the standard round
metric $g_1$ (of radius $1$) on each copy of $\SS^4$. This `conformal
blow-up' procedure gives a sequence of $g_\al$-anti-self-dual connections
$[\check A_\al]$ which converges {\it strongly} (in the sense of [D-K]) to
a limit $(A_0,I_1,\dots,I_{m_0})$ over the join $X_0\vee_{i=1}^{m_0}\SS^4$,
where the $I_i$ are the standard one-instantons over $X_i=\SS^4$ with
centre at the north pole $n$ and scale $1$. Here, strong convergence means
$C^\8$ convergence on compact sets away from the necks and such that
$c_2(A_0)+\sum_{i=1}^{m_0}c_2(I_i)=k$: there are no singular points and
there is no curvature loss over the necks. One obtains an
open neighbourhood in $\ol{M}^u_{X_0,k_0}(g_0)$ of the boundary point
$(A_0,x_1,\dots,x_{m_0})$ by gluing up the limit $(A_0,I_1,\dots,I_{m_0})$.

On the other hand, if the set $Z_0\equiv (x_1,\dots,x_{m_0})$ lies in the
diagonal of the symmetric product $s^{k-k_0}(X_0)$, the limiting behaviour
of the sequence $[A_\al]$ may be rather more complicated. Suppose $[\check
A_\al]$ is the corresponding sequence of $g_\al$-anti-self-dual connections
over $X=X_0\#_{i=1}^{m_0}\SS^4$ produced by conformal blow-ups. The
sequence $\check A_\al$ converges in $C^\8$ on compact subsets of $X_0\less
Z_0$ to a $g_0$-ASD connection $A_0$ over $X_0$, but in general only
converges weakly to an Uhlenbeck limit $(A_i,Z_i)$ over the four-spheres
$X_i\equiv\SS^4$, where $Z_i=(x_{i1},\dots,x_{im_i})$ is contained in
$X_i\less\{s\}$ and $s$ is the south pole. If the connection $A_i$, $i>0$,
is not flat, then the conformal blow-ups may be chosen so that it is {\it
centred} in the sense of [T88]: its mass centre lies at the north pole and
it has scale (essentially its `standard deviation') equal to $1$ (see \S
4.2). 

Unless all the singular sets $Z_i$ are empty, one can no longer produce an
open subset of the moduli space $M_{X_0,k}(g_0)$ simply by gluing up the
connections $(A_i)_{i=0}^{m_0}$: because of the nature of the convergence
process, some of the required moduli parameters have been lost in the
limit. 

Instead, the conformal blow-up process above must be iterated. The idea of
iterating conformal blow-ups has been suggested by Sacks and Uhlenbeck in
the context of harmonic maps of $\SS^2$ [S-U]. Taubes described an
iterative scheme of this type which is used to analyse the limiting
behaviour of sequences of connections with uniformly bounded Yang-Mills
functional and functional gradient tending to zero [T88]. Parker and
Wolfson described a bubble tree compactification for pseudoholomorphic
maps of Riemann surfaces into symplectic manifolds and noted that their
method should apply to the case of Yang-Mills connections over
four-manifolds [P-W].  

For the problem at hand, by repeatedly applying conformal blow-ups, we
obtain a sequence of $g_\al$-anti-self-dual connections $\check A_\al$ over
a large connected sum $X\equiv\#_{I\in\sI}X_I$. Here, $\sI$ is a set of
multi-indices $I$ obtained when the conformal blow-up process is iterated.
Thus, $\sI$  records the tree structure and if $I=0$, then $X_I$ is the
four-manifold $X_0$, while if $I\ne 0$, then $X_I$ is a copy of $\SS^4$. 
The construction of the `conformal blow-up maps' $f_{I\al}$ ensures that
the blow-up process must be repeated at most $k$ times in order to produce
a sequence of connections $[\check A_\al]$ which converge strongly to a
limit $(A_I)_{I\in\sI}$ over a join $\vee_{I\in\sI}X_I$, where $A_0$ is a
$g_0$-anti-self-dual connection over $X_0$ and each $A_I$, for $I\ne 0$, is
a $g_1$-anti-self-dual connection over $X_I=\SS^4$. The sequence of metrics
$g_\al$ converges in $C^\8$ on compact subsets away from the neck regions
to the metric $g_0$ on $X_0$ and the standard round metric $g_1$ on each
sphere $X_I$. This convergence scheme produces the `bubble tree
compactification' $\ol{M}^\tau_{X_0,k_0}(g_0)$ and is described in \S 4.3.

In particular, bubble tree degeneration and gluing are inverse to one
another in a natural way. 
One can now glue up the bubble tree limits $(A_I)_{I\in\sI}$ to form
$g$-anti-self-dual connections $A$ over a connected sum $X\equiv
\#_{I\in\sI}X_I$ using the techniques of [D-K] and construct open subsets
of the moduli space $M_{X,k}(g)$ by small deformations of the limit data.
The gluing procedure gives a collection of conformal maps $f_I$ (from a
small ball in a lower level summand $X_{I_-}$ to the complement in the
sphere $X_I$ of a small ball around the south pole) defined in exactly the
same way as the conformal blow-up maps $f_{I\al}$ above. Here, $g$ is a
$C^\8$ metric on $X$ which is conformally equivalent to the old metric
$g_0$ (via the maps $f_I$) and depends on the choice of gluing sites,
frames in the principle $\SO(4)$ frame bundle  $FX_0$, scales, and the
metric $g_0$ on $X_0$: its construction and properties are discussed in \S
3.5. Similar metrics over connected sums are described in [D86] and [T92]. 
Pulling back via the blow-up maps then gives $g_0$-anti-self-dual
connections $\hat A$ over $X_0$ and hence, produces open subsets of the
moduli space $M_{X_0,k}(g_0)$.  

Generalising the arguments in [D86] and [D-K] and employing the compactness
results of \S 4.3, one then shows that $\ol{M}_{X_0,k}^u(g_0)$ has a finite
cover consisting of gluing neighbourhoods $\ol\sV$. Of course, any
precompact open subset of $M_{X_0,k}(g_0)$ is covered by finitely many
Kuranishi charts and these comprise the `gluing charts' in this case.
Moreover, the $L^2$-metric geometry near the reducible connections, the
conical ends, has already been analysed by Groisser and Parker [G-P89],
so we may confine our attention to the more troublesome bubbling ends.  

\subsubsection{(ii) Upper bounds for the components of the $L^2$ metric} We
now outline a method of computing estimates for the $L^2$ metric ${\bf g}$
over the ends of the moduli space. 
In \S\S 3.3 and 3.4 we apply the techniques of [D86] and [D-K] to first construct approximate gluing maps $\sJ':\sT^0/\Ga\to \cB_{X,k}^*$, $t\to [A'(t)]$. Here, $X$ is the connected sum $\#_{I\in\sI} X_I$ with $C^\8$ metric $g$ conformally equivalent to $g_0$ on $X_0$, and $\sT^0/\Ga$ is a certain parameter space. 
If the $g$-self-dual curvature $F^{+,g}(A')$ is sufficiently small one can
then solve the $g$-anti-self-duality equation $F^{+,g}(A'+a)=0$, or
equivalently 
$$
d_{A'}^{+,g}a+(a\wedge a)^{+,g}=-F^{+,g}(A'), \eqlabel\eqnlabel{\eqIntroASDa}
$$
for $a\in\Om^1(X,\ad P)$. This gives a $C^\8$ family of $g$-anti-self-dual
connections $A\equiv A'+a$ and thus a gluing map $\sJ:\sT^0/\Ga\to
M_{X,k}^*(g)$, $t\to [A(t)]$. The solutions $a$ to Eq. \eqref{\eqIntroASDa}
are expressed in the form $a=P\xi$, where $\xi\in\Om^{+,g}(X,\ad P)$ and
$P$ is a right inverse to the
operator $d_{A'}^{+,g}$ (constructed as in [D-K] by patching together right
inverses 
$P_I$ for the operators $d_{A_I}^{+,g_I}$ over the summands $X_I$.
Therefore, Eq. \eqref{\eqIntroASDa} takes the
shape 
$$
\xi+(P\xi\wedge P\xi)^{+,g}=-F^{+,g}(A'). \eqlabel\eqnlabel{\eqIntroASDXi}
$$
Following [D-K], we assemble the framework required for solving Eq.
\eqref{\eqIntroASDXi} in \S 5.1.  

Now the $L^2$ metric ${\bf g}$ depends on the choice of metric $g_0$, not
just the conformal class $[g_0]$. So, using the conformal maps $f_I$, we
pull back the family of $g$-anti-self-dual connections $A(t)=A'(t)+a(t)$
over $X$ to an equivalent $C^\8$ family of $g_0$-anti-self-dual connections
$\hat A(t)=\hat A'(t)+\hat a(t)$ over $X_0$. Hence, we obtain gluing maps
$\hat\sJ:\sT^0/\Ga\to M_{X_0,k}^*(g_0)$, $t\to [\hat A(t)]$ analogous to
those constructed by Taubes. The properties
of the gluing maps $\sJ$ and $\hat\sJ$ are discussed in \S 5.2. 

The problem then is to estimate the differentials $D\hat\sJ$ and this task
is comprised of two parts. The first part is to bound the derivatives $\rd
A'/\rd t$: this local calculation is the subject of \S\S3.7 to 3.9 and the
main results are summarised in \S 3.10. 
The more difficult part is to bound the derivatives of the correction
terms, $\rd\hat a/\rd t$: this involves bounding the derivatives of global
operators such as $P$ and is described in \S\S 5.3 to 5.5. 
The problem of expressing bounds for derivatives of $\hat a(t)$ in terms of
bounds for derivatives of $a(t)$ is the subject of \S 3.5. Some care is
required here, since the conformal maps $f_I$ vary with the scale and
centre parameters, as does the metric $g$ in Eq. \eqref{\eqIntroASDXi}. 
The required estimates for the derivatives $\rd a/\rd t$ are then computed
in \S\S5.3 to 5.5 in terms of bounds for $\rd P/\rd t$ and $\rd\xi/\rd t$;
the estimates for $\rd\xi/\rd t$ are obtained implicitly from Eq.
\eqref{\eqIntroASDXi}. For the special case of a
neighbourhood of a point $(A_0,A_1)$ (with $H^2_{A_0}=0$), 
$L^2$ estimates for the derivatives $\rd A/\rd t$ were later obtained
independently by Peng using similar methods [Pe].

It is the estimates for derivatives with respect to the scales $\la_I$
which require the most care. For example, difficulties arise when bounding
the derivatives $\rd\hat A'/\rd \la_I$ because of the dependence on $\la_I$
of the conformal maps $f_I$ and the cut-off functions required to
patch the connections $A_I$ together over the connected sum. These
derivatives are ill-behaved as $\la_I\to 0$ and the necks of the connected
sum $X$ are pinched. Problems also occur when one attempts to bound $\rd
a/\rd\la_I$, since $a=P\xi$ and the construction of
$P$ involves cut-off functions 
with badly behaved derivatives with respect to $\la_I$ as $\la_I\to 0$. 
The final estimates for the differentials $D\hat\sJ$ and the corresponding
bounds for the $L^2$ metric ${\bf g}$ are sumarised in \S5.6. The constants
appearing in the bounds for ${\bf g}$ depend only on the gluing
neighbourhood. Theorem \ref{\thmFinVolDiam} then follows immediately from these estimates. 

\medskip
\noindent{\bf Acknowledgements.} The author would like to 
thank J. W. Morgan, T. Mrowka, D. H. Phong, and C. H. Taubes for 
helpful conversations. The current work was begun at
Columbia University and completed while the author was a postdoctoral
fellow at the Mathematical Sciences Research Institute in 1992--93.
He would like to thank these institutions and the Department of
Mathematics at Harvard University for their
generous support and hospitality.  

\endchapter
\eqnno=0
\secno=0
\thmno=0


\chapter{Preliminaries} In this Chapter we establish our notation
and define the $L^2$ metric. Unless stated otherwise, we adhere to the
standard conventions of [D-K]. 
For further details concerning gauge theory, we refer to
[D-K] or [F-U] and the references therein, while for details
concerning the 
$L^2$ metric, we refer to [G-P87], [G-P89]. 

Let $X$ be a closed, connected,
oriented, $C^\8$ four-manifold with Riemannian metric $g$ and let $P$ be a
principal $G$ bundle over $X$ with Lie algebra $\frak{g}$. As noted in the
Introduction, we will generally confine our attention in this article to
the case $G=\SU(2)$ for the sake of clarity.
We let $\Omega^l(P,\frak{g})$ denote the space of $C^\8$ $\frak{g}$-valued
$l$-forms, let $\ad P=P\times_\Ad\frak{g}$ be the adjoint bundle,
and let
$\Omega^l(X,{\rm ad}\,P)$ be the space of $C^\8$ $\ad P$-valued $l$-forms
on $X$. Let ${\cal A}_P$ be the affine subspace in $\Omega^1(P,\frak{g})$
of $C^\8$ connection 1-forms on $P$. 
For a connection $A$ on $P$, we let $\]_A$
be the corresponding covariant derivative, let $d_A$ be the
exterior covariant derivative, and let 
$F_A\in\Omega^2(X,\ad P)$ denote the curvature.

Let $\cG_P$ be the group of $C^\8$ bundle automorphisms or gauge
transformations. Recall that the isotropy group $\C_A\subset\cG_P$ of
a connection $A$ on $P$ is isomorphic to the centraliser of the holonomy
group of $A$ in $G$ and the centre $Z$ of the bundle
structure group $G$ is isomorphic to the centre of ${\cal
G}_P$. Thus $\C_A\supset Z$ and we let
${\cal A}_P^*$ be the dense open subset of connections
$A\in{\cal A}_P$ with $\C_A=Z$, so that $\cA_P^*$ is the space of
irreducible connections on $P$ when $G=\SU(2)$ or $\SO(3)$. 

The bundles $\Lambda^l T^*X\otimes\ad P$ have fibre metrics
$\<\ ,\ \>$ 
induced by the Riemannian metric $g$ on $X$ and the inner product on the
Lie algebra $\frak{g}$ 
given by $-1$ times the Cartan-Killing form: if $\xi_1,\xi_2\in\frak{g}$,
then $\<\xi_1,\xi_2\>=-\tr(\xi_1\xi_2)$. In particular, we may define Sobolev
spaces $L^p_n\Om^l(X,\ad P)$ in the usual way and consider the action of the
$L^2_{n-1}$ gauge transformations $\cG$ on the space of
$L^2_n$ connections $\cA_P$ (for $n>2$) with quotient $\cB_P=\cA_P/\cG_P$,
omitting the explicit Sobolev notation when no confusion can arise.

The tangent space $T_A{\cal A}^*_P$ is equal to $\W^1(X,\ad P)$ while
the tangent space to the ${\cal G}$-orbit through $A\in{\cal A}_P^*$
is $\im d_A\subset\W^1(X,\ad P)$. This induces an $L^2$-orthogonal
decomposition $T_A{\cal A}^*_P=\ker d_A^*\oplus \im d_A$,
where $\ker d_A^*\subset\W^1(X,\ad P)$.
There is an associated horizontal projection operator 
$\pi_A:T_A{\cal A}_P^*\to\ker d_A^*$, with
$\pi_A = 1 - d_A G^0_A d^*_A$, where $G^0_A$ is the Green's operator for
the Laplacian $\Delta^0_A=d_A^*d_A$. To identify the tangent space
$T_{[A]}{\cal B}^*_P$, introduce $C^\8$ paths $A(t)$ in ${\cal A}_P^*$ and
$u(t)$ in ${\cal G}_P$, $u(0)=1$. If $A^u(t)\equiv u_t(A_t)$, then  
$$
\deriv{A^u}{t} = \Ad(u^{-1})\deriv{A}{t} +
d_{A^u}\left(u^{-1}\deriv{u}{t}\right). \eqlabel 
$$
Thus $ dA/dt(0)$ defines an element of $\W^1(X,\ad P)/\im d_A$ and therefore the tangent space
$T_{[A]}{\cal B}_P^*$ is given by $\W^1(X,\ad P)/\im d_A \simeq \ker d_A^*$.  

Let $M_P(g)$ be the moduli space of $g$-anti-self-dual connections on the
$G$ bundle $P$ over $X$, that is
$\{[A]\in {\cal B}_P: F^{+,g}(A)= 0\}$, and let $M^*_P(g)$ be the
dense open subset  
$M_P(g)\cap \cB_P^*$. If $A(t)$ is a $C^\8$ path in ${\cal
A}_P$ satisfying $F^{+,g}(A(t))= 0$, then  
$d A/dt(0)$ defines an element of $\ker d_A^*/\im d_A^{+,g}$.
The $g$-anti-self-dual condition $F^{+,g}(A)= 0$ is equivalent to
$d^{+,g}_A\circ d_A=0$ and so we have the elliptic deformation complex  
$$
\W^0(X,\ad P)\mapright{d_A}\W^1(X,\ad P)
\mapright{d^{+,g}_A}\W^{+,g}(X,\ad P) \eqlabel\eqnlabel{\eqEllDefCplx}
$$
with associated cohomology groups $H^*_A$, where $H^0_A$ is the Lie
algebra of $\C_A$, the group $H^1_A = \ker d^{+,g}_A/\im d_A$ is just the
tangent space $T_{[A]}M_P(g)$, and $H^2_A = \coker
d^{+,g}_A$. By Hodge theory there are natural isomorphisms
$H_A^0\simeq\ker\D_A^0$,  
$H_A^1\simeq\ker d_A^*\cap \ker d_A^{+,g}$, and  
$H_A^2\simeq\ker\D_A^{+,g}$, where the Laplacian 
$\D_A^{+,g}$ is equal to $d_A^{+,g}(d_A^{+,g})^*$.
 
If $[A]$ is an irreducible point of $M_P(g)$, then $H_A^0=0$, and an
irreducible point $[A]$ is regular if $H_A^2=0$. The moduli space
$M_P(g)$ is regular if all its irreducible points are regular points,
and in that case, $M^*_P(g)$ is a $C^\8$ manifold of dimension  
$$
\dim M_P(g)=8k(P)-3(1-b_1(X)+b^+(X)), \eqlabel\eqnlabel{\eqDimMod}
$$
with tangent space $T_{[A]}M^*_P(g)=H_A^1$ at the point $[A]$. 

According to the Freed-Uhlenbeck theorems, the anti-self-dual
moduli spaces $M^*_P(g)$ are smooth manifolds when $g$ is generic.
More precisely, if 
$b^+(X)>0$, $P$ is any $\SU(2)$ or $\SO(3)$ bundle $P$ over $X$, and the
metric $g$ on $X$ is generic, then 
(1) $M^*_P(g)$ contains no points $[A]$ with $H_A^2\ne 0$;
(2) If $b^+(X)>0$ and $l>0$ then 
$M_P(g)$ contains no points $[A]$ with $H_A^0\ne 0$ for any bundle
$P$ with $0<k(P)\le l$;
(3) If $b^+(X)=0$ and $P$ is non-trivial, then 
the cohomology groups $H^2_A$ are zero for 
all the reducible $g$-anti-self-dual connections $A$ on $P$, and a
neighbourhood of point 
$[A]\in M_P(g)$ with $H_A^0\ne 0$ is homeomorphic to a cone over
$\cp^{4k-2}$ and diffeomorphic away from the cone point $[A]$. 

It remains to define the $L^2$ metric.
The quotient space $\cB_P^*$ inherits a (weak) Riemannian $L^2$ metric ${\bf
g}$ by requiring that the projection map for the principal $\cG_P/Z$ bundle
$\cA_P^*\to\cB_P^*$ be a Riemannian 
submersion: if $[a], [b]$ are tangent vectors in $T_{[A]}\cB_P^*$, then
$$
{\bf g}_{[A]}([a],[b])\equiv \int_X\<\pi_A a,\pi_A b\>\, dV_g, \eqlabel
\eqnlabel{\L2metric}
$$
and this restricts to give a $C^\8$ Riemannian metric
${\bf g}$ on the moduli space $M^*_P(g)$. 

\endchapter
\eqnno=0
\secno=0
\thmno=0


\chapter{Differentials of the approximate gluing maps}
Our purpose in this Chapter is to construct the approximate gluing
maps $\script{J}':\sT/\Ga\to\cB^*_{X,k}$ and 
$\hat\script{J}':\sT/\Ga\to\cB^*_{X_0,k}$, and to estimate the
differentials $D\script{J}'$, and especially $D\hat\script{J}'$. 
The construction of $\script{J}'$ uses the method employed by Donaldson in
[D86], [D-K]. The induced maps $\hat\script{J}'$  are
essentially the approximate gluing maps described by Taubes in [T82],
[T84a], [T88]. In the former case, we obtain an almost
$g$-anti-self-dual 
connection $A'$ over a connected sum $X=X_0\#_{I\in\sI}\SS^4$ with metric
$g$ conformally equivalent to $g_0$ on $X_0$, while in the latter case we
obtain an almost $g_0$-anti-self-dual connection $\hat A'$ over $X_0$ with
its fixed 
metric $g_0$. In Chapter 5, we obtain a system of coordinate charts
$\hat\script{J}:\sT/\Ga\to M^*_{X_0,k}(g_0)$
covering the moduli space by perturbing the maps $\hat\script{J}'$
using the techniques of [D-K] for solving the anti-self-dual equation. 

\section{Preliminary estimates for connections and curvature}
We describe some pointwise estimates for local connection one-forms and
curvature two-forms.  
We first consider estimates for connection one-forms in radial gauge on a
$C^\8$ manifold $X$ with $C^\8$ metric $g$. Suppose $P\to X$ is a principal
$G$ bundle, $A$ is a $C^\8$ connection 
on $P$, and $B$ is an open geodesic ball centred at $x_0\in X$ with radius
$\varrho/2$, where $\varrho$ is the injectivity radius of $(X,g)$.
Define a $C^\8$ local section $\si: B\to P$ by parallel transport of 
a point in the fibre $P|_{x_0}$
along radial geodesics through $x_0$. If $\ga$ is a radial geodesic in $B$ with $\ga(0)=x_0$ and $\dot\ga(t)=\xi_t$, then
$\si^*A(x_0)=0$ and
$\iota_{\xi_t}\si^*A(\ga(t))=0$, $t>0$.
If $\phi^{-1}:B\to {\Bbb R}^4$ is a geodesic normal coordinate system
centred at $x_0$ and we define a geodesic $\ga$ by $\ga(t)=\phi(tx)$, $x\in
B$, $t\in[0,1]$, 
then $\ga^\mu(t)=tx^\mu$, $\dot\ga={\bf x}$, and
$\iota_{\bf x}\si^*A=x^\mu (\si^*A)_\mu$.
We recall the following estimates for local connection one-forms in radial
gauge. 

\proclaim{Lemma}{\thmlabel{\lemRadConnEst}
{\rm [U82, p. 14]} Let $A$ be a $C^\8$ connection on a principal $G$
bundle $P\to X$, where $X$ is a $C^\8$ manifold with $C^\8$ metric $g$, $B$
is a geodesic ball of radius $\varrho/2$ centred 
at $x_0\in X$, $\si:B\to P$ is a local section such that $\si^*A$ is in
radial gauge centred at $x_0$, and $\phi^{-1}:B\to {\Bbb R}^n$ is a
geodesic normal coordinate system centred at $x_0$. If
$K=\|F_A\|_{L^\8(B,g)}$, then 
$|\phi^*\si^*A|_g(x)\le K|x|,\qquad |x|<\varrho/2$.}

Let $\hp^1$
be the right quaternionic projective space,
with the standard identifications 
${\Bbb H}\simeq{\Bbb R}^4$ and 
$\hp^1\simeq{\Bbb S}^4$.
Coordinate patches for ${\Bbb S}^4$ may then be defined by 
$U_n=\{[x,y]: y\ne 0\}={\Bbb S}^4\less\{s\}$ and 
$U_s=\{[x,y]:x\ne 0\}={\Bbb S}^4\less\{n\}$
covering the north pole $n=[0,1]$ and 
south pole $s=[1,0]$, respectively. We let 
$\phi_n^{-1}:U_n\to{\Bbb R}^4$, $[x,y]\mapsto xy^{-1}$ and
$\phi_s^{-1}:U_s\to{\Bbb R}^4$, $[x,y]\mapsto yx^{-1}$
denote the standard local coordinate charts. If $g_1$ is 
the standard round metric of radius
1 on $\SS^4$, then
$$
(\phi_\a^*g_1)_{\mu\nu}(x)=h_1^2(x)\d_{\mu\nu}
={4\over (1+|x|^2)^2}\d_{\mu\nu}, \qquad x\in\RR^4, \eqlabel
$$
for $\a=n,s$; the standard flat metric on 
${\Bbb R}^4$ is denoted by $\d$. 

Let $A$ be a $C^\8$ connection on a principal $G$ bundle $P\to{\Bbb S}^4$,
where ${\Bbb S}^4$ has its standard metric $g_1$. We define a system of
local sections $\s_\a: U_\a\to P$, $\al=n,s$, by parallel transport 
of points in the fibres $P|_\a$ along radial geodesics
through the north or south poles. The estimates below follow easily since
$A$ is smooth over $\SS^4$ with metric $g_1$:

\proclaim{Lemma}{\thmlabel{\lemFSphereEst}
Let $A$ be a $C^\8$ connection on a principal $G$ bundle 
$P\to {\Bbb S}^4$, where ${\Bbb S}^4$ has metric $g_1$ and $K=\|F_A\|_{L^\8({\Bbb S}^4,g_1)}$. Then, for $\a,\b\in\{n,s\}$, 
$$
|\phi_\b^*F(\s_\a^*A)|_\d(x) \le 4K{1\over(1+|x|^2)^2}\quad\hbox{for }
\cases{x\in {\Bbb R}^4 &if $\a=\b$ \cr
       x\in {\Bbb R}^4\less\{0\} &if $\a\ne\b$.\cr}
$$}

\proclaim{Lemma}{\thmlabel{\lemRadConnSphEst} Given the hypotheses of\/ {\rm Lemma \ref{\lemFSphereEst}}, if  
the local connection one-forms $\s_\a^*A$ are in radial gauge, then
$|\phi_\a^*\s_\a^*A|_{g_1}(x) \le K|x|$, for $x\in {\Bbb R}^4$ and $\al=n,s$.}


\section{Connections over the four-sphere and conformal diffeomorphisms} 
\seclabel{\secConnSphere}
Recall that the group of conformal diffeomorphisms of $\SS^4$ acts on the
space $\cA_P$ of $C^\8$ connections on a $G$ bundle $P$ over $\SS^4$. The
group ${\bf D}\times{\bf T}$ of dilations and translations of $\RR^4$ may be
identified with a
subgroup of the conformal group of $\SS^4$. Hence, in this section we
discuss some aspects of the induced action of $\RR^+\times \RR^4$ on
the space $\cA_P$. For related material we refer to [D83], [F-U], [G-P87], 
[G-P89], and [T88].

Let $P$ be a $G$ bundle with $C^\8$ connection $A\in\Om^1(P,\frak{g})$ over
a $C^\8$ manifold $X$ and 
suppose $\varphi_t$ is a $C^\8$ one-parameter group of diffeomorphisms of
$X$ generating a vector field $\xi\in C^\8(TX)$. Let $\tilde\xi\in
C^\8(TP)$ be the 
horizontal vector field covering $\xi$ and let $\tilde\varphi_t$ be the
one-parameter group of diffeomorphisms of $P$ generated by $\tilde\xi$.
Then $\tilde\varphi_t$ commutes with right $G$ multiplication and covers
$\varphi$.  
Fixing $\Om\in\Om^1(P,\frak{g})$, we obtain a $C^\8$ one-parameter
family of $C^\8$ one-forms $\tilde\varphi_t^*\Om$ on $P$ with 
$$
\deriv{\tilde\varphi_t^*\Om}{t}\Bigr|_{t=0}=\sL_{\tilde\xi}\Om,\eqlabel
$$
where $\sL_{\tilde\xi}\Om\in\Om^1(P,\frak{g})$ denotes the Lie derivative
of $\Om$ with respect to $\tilde\xi$; in particular,  
$\tilde\varphi_t^*A$ is a $C^\8$
one-parameter family of $C^\8$ connection one-forms on $P$. 

\proclaim{Lemma}{\thmlabel{\lemLieDerConn} Let $P$ be a $G$ bundle with
connection $A\in\Om^1(P,\frak{g})$ over a manifold $X$. Given a vector
field $\xi\in C^\8(TX)$, let $\tilde\xi\in C^\8(TP)$ be its horizontal lift. If
$F_A\in\Om^2(P,\frak{g})$ is the curvature of $A$, then  
$\sL_{\tilde\xi}A=\iota_{\tilde\xi}F_A$.}

\pf Since $\tilde\xi$ is horizontal, then $A(\tilde\xi)=0$ and so
for any vector field $\eta\in TP$, we have
$(\sL_{\tilde\xi}A)(\eta)
=(\iota_{\tilde\xi} dA+d\iota_{\tilde\xi}A)(\eta)
=dA(\eta,\tilde\xi)$.
But $F_A(\eta,\tilde\xi)=dA(\eta,\tilde\xi)+\half[A(\eta),A(\tilde\xi)]$ and so the result follows.\qed

We also need to consider Lie derivatives of $\ad P$-valued one-forms.
Recall that if $\pi:P\to X$ is the bundle projection, there is an
injective map $\pi^*:\Om^1(X,\ad P)\hookrightarrow \Om^1(P,\frak{g})$.  
The one-forms $\Om$ in the image of $\pi^*$ are characterised by the
properties (a) $R_u^*\Om=\Ad(u^{-1})\Om$, for all $u\in G$, and (b)
$\Om(\eta)=0$ if $\eta\in TP$ is vertical. Hence, the action of
$\tilde\varphi_t$ on $\Om^1(P,\frak{g})$ induces an action on $\Om^1(X,\ad
P)=\Ga(T^*X\otimes \ad P)$.
Thus, if $\om\in\Om^1(X,\ad P)$, we obtain a $C^\8$ one-parameter family of $C^\8$ $\ad P$-valued one-forms $\tilde\varphi_t^*\om$ on $X$ with
$$
\deriv{\tilde\varphi_t^*\om}{t}\Bigr|_{t=0}=\sL_{\tilde\xi}\om, \eqlabel
$$
where $\sL_{\tilde\xi}\om\in\Om^1(X,\ad P)$ denotes the Lie derivative of $\om$ with respect to $\tilde\xi$.

For the purposes of calculation, it is useful to phrase the preceding discussion in terms of local one-forms on $X$. It is convenient to choose a system of local sections $\si_\al:U_\al\to P$ which are {\it parallel} with respect to the connection $A$ and vector field $\xi$, in the sense that $A(\si_{\al*}\xi)=0$. 
For example, one can try to construct $\si_\al$ by first choosing a section $\si_\al|_{V_\al}$, where $V_\al$ is a submanifold of $U_\al$ transverse to the vector field $\xi$, and then extend by parallel translation along integral curves of $\xi$ to construct a section $\si_\al$ over a tubular neighbourhood $U_\al$ of $V_\al$. Local sections of this type are described in [U82a, pp. 14-15] and [F-U, pp. 146-147]. 

Given a system of $(A,\xi)$-parallel local sections $\si_\al$, we have  $\tilde\xi=\si_{\al*}\xi$ and $\varphi_t=\si_\al^*\tilde\varphi_t$ over $U_\al$. Hence, for $\om\in\Om^1(X,\ad P)$ we see that
$\si_\al^*\tilde\varphi_t^*\om=\varphi_t^*\si_\al^*\om$ and
$\si_\al^*\sL_{\tilde\xi}\om=\sL_\xi\si_\al^*\om$ on $U_\al$,
and similarly for $A\in\Om^1(P,\frak{g})$.
Indeed, one can see that the transition functions $\{u_{\al\be}\}$ are constant along the vector field $\xi$. For if $\si_\be=\si_\al u_{\al\be}$, then
$\si_{\be*}\xi=\si_{\al*}\xi\cdot u_{\al\be}+\si_\al\cdot u_{\al\be*}\xi$,
which gives
$A(\si_{\be*}\xi)=\Ad(u_{\al\be}^{-1})A(\si_{\al*}\xi)+A(\si_\al\cdot u_{\al\be*}\xi)$,
and thus 
$du_{\al\be}(\xi)=0$,
since $A(\si_{\al*}\xi)=A(\si_{\be*}\xi)=0$ and
$A(\si_\al\cdot u_{\al\be*}\xi)=u_{\al\be*}\xi$. Here, $\si_\al\cdot u_{\al\be*}\xi$ is the vector field on $P|_{U_\al}$ obtained by differentiating the maps $G\to P$ given by $u\mapsto \si_\al(x)u$.
When computing Lie derivatives of local connection one-forms or $\ad
P$-valued one-forms with respect to a vector field $\xi$, we shall always
require that the local sections $\si_\al$ be $(A,\xi)$-parallel.

It is often useful to express $\sL_{\tilde\xi}\om$ in terms of covariant
derivatives. Suppose $X$ has a $C^\8$ metric $g$. 
We have $\sL_{\xi}\om=\iota_\xi d\om+d\iota_\xi\om$, or in local coordinates,
$(\sL_\xi\w)_\mu=\xi^\nu\rd{\w_\mu}/\rd{x^\nu}+\w_\nu\rd{\xi^\nu}/\rd{x^\mu}$.
We find that 
$$
\sL_{\tilde\xi}\om=\cov_\xi^{A,g}\om+\om(\cov^g\xi), \eqlabel
$$
using normal geodesic coordinates $\{x^\mu\}$ and $(A,\xi)$-parallel local
sections $\{\si_\al\}$.
In the sequel, we omit the ``tildes'' to indicate lifts of vector fields or
diffeomorphisms on the base to the total space of a principal bundle --
this being understood from the context. 
Remark that if $\Phi:X\to X$ is a diffeomorphism and $\om\in\W^1(X,\ad P)$, then we have $\sL_\xi\Phi^*\w=\Phi^*\sL_{\Phi_*\xi}\w$. 

Let $A$ be a $C^\8$ connection on a $G$ bundle $P$ over $\SS^4$ and let $\om\in\W^1(\SS^4,\ad P)$. 
For any $t\in(-\8,\8)$, let $\de_t$ be the dilation 
of ${\Bbb R}^4$ given by $x\mapsto e^tx$ and for any $p\in\RR^4$, let $\tau_p$ be the translation of ${\Bbb R}^4$ defined by $\tau_p: x\mapsto x-p$. If
$\de_t$ and $\tau_p$ again denote the conformal diffeomorphisms of $\SS^4$ induced by the chart $x=\phi_n^{-1}$, then the group ${\bf C}=\SO(4)\times{\bf D}\times {\bf T}$ of rotations, dilations, and translations of ${\Bbb R}^4$ is identified with the subgroup in $\Conf(\SS^4,g_1)$ of diffeomorphisms which fix the south pole $s\in \SS^4$. Setting $\varphi_t=\de_t$ or $\tau_{tp}$, these diffeomorphisms are generated by the vector fields
$$
\rad\equiv x^\mu\pd{}{x^\mu}
\quad\hbox{and}\quad-\trans\equiv -p^\mu\pd{}{x^\mu}.\eqlabel
$$
We always choose $p\in\RR^4$ with $|p|\le 1$.
We next describe the construction of $(A,\xi)$-parallel local sections $\si_\al$ for $\xi=\rad$ or $\trans$. 

Considering the group of dilations ${\bf D}$, let $\s_n,\s_s$ be the local sections formed by choosing points in the fibres $P|_n,P|_s$ and then
parallel translating along radial directions from the poles. The transition function $u$ will be constant along the radial directions, $du(\rad)=0$, and the local connection one-forms $\si_\al^*A$ are in {\it radial} gauge. 
On the other hand, considering the group of translations ${\bf T}$, suppose first that $\trans=\rd/\rd x^4$ and let $\s_n|_{\SS^3},\s_s|_{\SS^3}$ be the local sections formed by parallel translation from the north and south poles of the three-sphere $\SS^3\subset \SS^4$ defined by the image of the $x^1x^2x^3$-plane under the map $\phi_n:\RR^4\to\SS^4\less\{s\}$. We obtain local sections $\s_n,\s_s$ by parallel translation along the $x^4$-axis. The transition function $u$ will now be constant along the $x^4$-axis, $du(\trans)=0$, and the local connection one-forms $\si_\al^*A$ are in a {\it transverse} gauge. By a linear change of coordinates, the same argument applies to arbitrary translations. 

For the dilations, we have
$$
\deriv{\de_t^*\om}{t}\Bigr|_{t=0}
=\sL_{\rad}\om=\iota_\rad d\om+\om, \eqlabel\eqnlabel{\eqLieRadOm}
$$
using $\sL_{\xi}\om=\iota_\xi d\om+d\iota_\xi\om$, or in local coordinates,
$(\sL_\xi\w)_\mu
=\xi^\nu\rd{\w_\mu}/\rd{x^\nu}+\w_\nu\rd{\xi^\nu}/\rd{x^\mu}$. 
Similarly, for the translations we have
$$
\deriv{\tau_{tp}^*\om}{t}\Bigr|_{t=0}=-\sL_{\trans}\om=-\iota_\trans d\om, 
\eqlabel\eqnlabel{\eqLieTransOm}
$$
where $\trans=p^\mu\rd/\rd x^\mu$.

For any $\la\in(0,\8)$, let $c_\l$ be the diffeomorphism of ${\Bbb S}^4$ defined by the chart $x=\phi_n^{-1}$ and the dilation $c_\la$ of $\RR^4$ given by 
$x\mapsto x/\l$. Then $c_\la=\de_t$ with $t=-\log\la$, and so from Eq. \eqref{\eqLieRadOm} we have
$\pd{}{\l} c_\l^*\w=-{1\over\la}c_\l^*\sL_\rad\w$.
Similarly, for the translations $\tau_q$, $q\in \RR^4$, we see that Eq. \eqref{\eqLieTransOm} gives
$\pd{}{p}\tau_q^*\om=-\tau_q^*\sL_\trans\om$,
where $\rd/\rd p\equiv p^\mu\rd/\rd{q^\mu}$ on the left-hand side and using $\tau_{q+tp}=\varphi_{tp}\circ\tau_q$ on the right. 
Combining these actions, we find that
$$
\pd{}{\l}\tau_q^*c_\l^*\om=-{1\over\la}\tau_q^*c_\l^*\sL_\rad\om
\quad\hbox{and}\quad
\pd{}{p}\tau_q^*c_\l^*\om=
-{1\over\la}\tau_q^*c_\l^*\sL_\trans\om,
\eqlabel\eqnlabel{\eqDctauStarOmDLaq}
$$
Similarly, considering the action of the dilations $c_\la$ and translations $\tau_q$ on connection one-forms, we have
$$
\pd{}{\l}\tau_q^*c_\l^*A=-{1\over\la}\tau_q^*c_\l^*\iota_\rad F_A
\quad\hbox{and}\quad
\pd{}{p}\tau_q^*c_\l^*A=
-{1\over\la}\tau_q^*c_\l^*\iota_\trans F_A.
\eqlabel\eqnlabel{\eqDctauStarADLaq}
$$
These derivative formulas play a significant role in the sequel. 

It is convenient at this point to recall Taubes' definition of a centred connection over the four-sphere [T88, p. 343].  
Let $A$ be an $g_1$-ASD connection on a $G$ bundle $P$ with $c_2(P)=k$ over $\SS^4$ with its standard metric $g_1$. Pulling back via the chart $x=\phi_n^{-1}:\SS^4\less\{s\}\to\RR^4$, we obtain a $\de$-ASD connection $A$ on a $G$ bundle $P$ over $\RR^4$ with its standard metric $\de$. Let $\Theta$ denote the flat connection on the product bundle.
Suppose $A\ne\Theta$: then the {\it mass centre} $q$ and {\it scale} $\la$ are defined by
$$\leqalignno
{q&=\Centre[A]\equiv
{1\over 8\pi^2k}\int_{{\Bbb R}^4}|F_A|^2\, d^4x,
&\numeqn\eqnlabel{\eqCentScalSph}\cr 
\la^2&=\Scale^2[A]\equiv{1\over 8\pi^2k}
\int_{{\Bbb R}^4}|x-q|^2|F_A|^2\, d^4x. &\cr}
$$
If $A=\Theta$, we set $\Centre[A]=0$ and $\Scale[A]=0$.
The connection $A$ is called
{\it centred} if $\Centre[A]=0$ and $\Scale[A]=1$. Eq. \eqref{\eqCentScalSph} leads to the following {\it Tchebychev ineqality}:
$$
\int_{|x-q|\ge R\la}|F_A|^2\, d^4x \le 8\pi^2k R^{-2},\qquad R\ge 1.
\eqlabel\eqnlabel{\eqTchebySphere}
$$
Hence, the ball $B(q,R\la)$ contains $A$-energy greater than or equal to
$8\pi^2k(1-R^{-2})$.  

Setting $f_{\la,q}=c_\la\circ\tau_q$, we see that $\Centre[(f_{\la,q}^{-1})^*A]=0$ and $\Scale[(f_{\la,q}^{-1})^*A]=1$.
Let $M_k$ denote the moduli space of $g_1$-ASD connections on the bundle
$P$ over $\SS^4$ and let $M_k^0$ denote the {\it moduli space of centred
$g_1$-ASD connections}. Note that $M_1^0$ consists of a single point
representing the standard one-instanton over $\SS^4$. More generally, the
relationship between $M_k$ and $M_k^0$ is explained below:

\proclaim{Proposition}{\thmlabel{\propCentModDiff}
For any $k>0$, the space $M_k^0$ is a smooth submanifold of $M_k$. Moreover, $M_k$ is diffeomorphic to 
$M_k^0\times{\Bbb R}^4\times (0,\8)$.}

\pf One argues as in [T88, pp. 343-344] and [T84b, pp. 365-367]. Given $[A]\in M_k$ with $\Centre[A]=q$ and $\Scale[A]=\la$, set $f_{\la,q}=c_\la\circ\tau_q$.
The map $[A]\to ([(f_{\la,q}^{-1})^*A],q,\la)$ then
gives the required diffeomorphism. \qed


\section{Gluing construction of approximately anti-self-dual connections}
\seclabel{\secGluCon}
We describe the approximate gluing constructions of Donaldson [D86], [D-K],
and Taubes [T82], [T84a], [T88], adapted to the case of ``bubble trees''.
For clarity, we first discuss the construction of approximately
anti-self-dual connections over single 
connected sums.
Let $X_0$ be our closed, smooth four-manifold with metric
$g_0$ and injectivity radius $\varrho_0$, and let $X_1={\Bbb S}^4$ with its
standard round metric $g_1$ of radius $1$. Let $x_1$ be a point in $X_0$
and let $x_{1n}, x_{1s}$ denote the north and south poles of $X_1$. 
Let $P_i\to X_i$ be principal $G$ bundles with $c_2(P_i)=k_i$, $i=0,1$. Let
$FX_0$ be the principle $\SO(4)$ bundle of oriented, orthnormal frames over
$X_0$.  

A choice of frame $v_1\in FX_0|_{x_1}$
defines a geodesic normal coordinate system
$\phi_1^{-1}=\exp_{v_1}^{-1}: B_1(\varrho_0)\to{\Bbb R}^4$.
Denote $\phi_{1\al}=\phi_\al$, $\a=s,n$, where $\phi_\a^{-1}: 
U_\a={\Bbb S}^4\less\{\a\}\to {\Bbb R}^4$ are the standard coordinate charts on the four-sphere. Let
$B_1(r)=B(x_1,r)$ be the open geodesic ball in $X_0$ with centre $x_1$ and radius $r$,
and let $B_{1s}(r)=\phi_{1s}(\{x\in\RR^4:|x|<r\})$, an open ball in $X_1$ with centre $x_{1s}$.
Let $\W_1(r,R)=\W(x_1,r,R)$ be the open annulus $B_1(R)\less\ol B_1(r)$ centred at $x_1\in X_0$, with inner radius $r$ and outer radius $R$; similarly, let
$\W_{1s}(r,R)=\W(x_{1s},r,R)$ be the open annulus $B_{1s}(R)\less\ol B_{1s}(r)$ in $X_1$. 

Let $N>4$ be a large parameter, to be fixed later, and let 
$\l_1>0$ be a small scale parameter such that $\l_1^{1/2}N\ll 1$. We define
open sets  
$X_0'=X_0\less \ol{B}_1(N^{-1}\l_1^{1/2})$, 
$X_0''=X_0\less \ol{B}_1(\half\l_1^{1/2})$, and 
$X_0'''=X_0\less \ol{B}_1(2N\l_1^{1/2})$ ---
the complements in $X_0$ of small balls around the point $x_1$. Likewise, define open sets $X_1'$, $X_1''$, and $X_1'''$ in the sphere $X_1$. 
Let $\W_1$ denote the annulus $\Om_1(N^{-1}\l_1^{1/2},N\l_1^{1/2})$ in $X_0$ and let $\Om_{1s}=\Om_{1s}(N^{-1}\l_1^{1/2},N\l_1^{1/2})$ 
be the corresponding annulus in $X_1$.
Let $c_1$ be the dilation map on
${\Bbb R}^4$ defined by $x\mapsto x/\l_1$.
Define balls $B_1'=B_1(N\l_1^{1/2})$ and
$B_1''=B_1(2\la_1^{1/2})$ centred at $x_1$ in $X_0$ and a diffeomorphism
$$
f_1= \phi_{1n}\circ c_1\circ \phi_1^{-1}: B_1'\too X_1'. \eqlabel
$$
Hence, $f_1$ identifies the small balls $B_1'$ and $B_1''$ in $X_0$ with
the open sets $X_1'$ and $X_1''$ in $X_1$, and restricts to a
diffeomorphism $f_1:\W_1\to\W_{1s}$. 

We let $X$ be the connected sum $X_0\#_{f_1}X_1$. In \S 3.5 we define a smooth metric $g$ on $X$ which closely approximates the metrics $g_i$ on each summand $X_i'$ and such that the map $f_1: B_1'\to X_1'$ is conformal. Thus, $(X,g)$ is conformally equivalent to $(X_0,g_0)$.

Let $A_i$ be $g_i$-anti-self-dual connections on the bundles $P_i\to X_i$, $i=0,1$.
The connections $A_0, A_1$, together with a choice of points in the fibres $P_0|_{x_1},P_1|_{x_{1s}}$, define local sections $\s_1:B_1(\varrho_0)\to P_0$ and 
$\s_{1s}:X_1\less \{x_{1n}\}\to P_1$ by
parallel transport along radial geodesics
through $x_1$, $x_{1s}$. Hence, we obtain local
trivialisations $P_0|_{B_1}\simeq B_1\times G$ and $P_1|_{B_{1s}}\simeq B_{1s}\times G$. 

Let $b_1\ge 4N\la_1^{1/2}$ be a small parameter, $b_1< \quarter\min\{1,\varrho_0\}$: we will eventually set $b_1=4N\la_1^{1/2}$. Choose
cutoff functions $\psi_i$ on $X_i$ such that $0\le \psi_i\le 1$, with
$\psi_0=1$ on $X_0\less B_1(b_1)$,
$\psi_0=0$ on $B_1(b_1/2)$, and similarly for $\psi_1$ on $X_1$.
We let $A_0'=\psi_0A_0$ be the $C^\8$ connection on 
the bundle $\pi_0:P_0\to X_0$ defined by  
$$
A_0'=\cases{A_0& on $P_0|_{X_0\less B_1(b_1)}$\cr
           \pi_0^*(\psi_0\s_1^*A_0)& on $P_1|_{B_1(b_1)}$.\cr} \eqlabel
$$
Of course, we have the analogous definition for the $C^\8$ connection $A_1'$ over $X_1$ and we obtain almost anti-self-dual connections which are flat on the balls $B_1'$, $B_{1s}'$.

To construct the cutoff functions $\psi_i$, choose a $C^\8$ bump function $\z$ on ${\Bbb R}^1$ such that
$\z(t)=1$ for $t\ge 1$ and $\z(t)=0$
for $t\le 1/2$. 
Define a $C^\8$ cutoff function $\psi_b$ on ${\Bbb R}^4$ by
$\psi_b(x)=\zeta(|x|/b)$, for any $b>0$. Set $\psi_0=(\phi_1^{-1})^*\psi_{b_1}$ and extend by $1$ on $X_0\less B(x_1,b_1)$ and by zero on $B(x_1,b_1/2)$ to give $\psi_0\in C^\8(X_0)$; likewise, set $\psi_1= (\phi_{1s}^{-1})^*\psi_{b_1}$ and extend to give $\psi_1\in C^\8(X_1)$. Each $\psi_i$ extends by zero to give a $C^\8$ cutoff function on the connected sum $X$.

Choose a $G$-isomorphism $\r_1\in \Gl_{x_1}$, where  
$\Gl_{x_1}\equiv\Hom_G(P_0|_{x_1}, P_1|_{x_{1s}})\simeq G$ is  
the space of ``gluing parameters''. Using the connections $A_i$ over the
small $\half b_1$-balls, spread out the fibre isomorphism $\r_1$ to give a
bundle isomorphism $\tilde\rho_1:P_0|_{\W_1}\to P_1|_{\W_{1s}}$ covering
the diffeomorphism $f_1:\W_0\to\W_1$. Thus,
$\si_1\tilde\rho_1=f_1^*\si_{1s}$ on $\Om_1$. We define  
the smooth connected-sum bundle $P\to X$ with second Chern class
$c_2(P)=k=k_0+k_1$ by setting   
$P|_{X_0'} = P_0|_{X_0'}$ and $P|_{X_1'} = P_1|_{X_1'}$. Note that the bundle $P$ is defined by transition functions independent of the scale $\la_1$.
We define a smooth connection $A'=A_0'\#A_1'$ on
$P\to X$ by setting $A'=A_i'$ on each summand $X_i'$. 

If $\C_{A_i}$ are the isotropy
groups of the connections $A_i$ and 
$\C=\C_{A_0}\times\C_{A_1}$, then we recall that the gluing construction
gives a bijection between the gauge equivalence classes $[A'(\r_1)]$ in
$\cB_{X,k}$ and $\Gl_{x_1}/\Ga$ [D-K, p. 286].  

Using the diffeomorphism $f_1:B_1'\to X_1'$, we pull back the bundle $P$ over $X$ to a bundle $\hat P$ over $X_0$, given by
$\hat P|_{X_0'} = P_0|_{X_0'}$ and
$\hat P|_{B_1'} = f_1^*P_1|_{B_1'}$. 
We have an induced system of local sections of $\hat P|_{B_0'}$ given near $x_1$ by
$\hat\s_{1n}=f_1^*\si_{1n}:B_1'\to P$,
$\hat\s_{1s}=f_1^*\s_{1s}:B_1'\less\{x_1\}\to P$, 
and $\hat\s_1=\s_1:\Om_1(N^{-1}\la_1^{1/2},\varrho_0)\to P$.
The corresponding transition functions  
$\hat u_1= f_1^*u_1:B_1'\less\{x_0\}\to G$
and $\tilde\rho_1:\W_1\to G$ 
are determined by $\hat\s_{1s}=\hat\si_{1n}\hat u_1$ on $B_1'\less\{x_0\}$
and $\si_1\tilde\rho_1=f_1^*\hat\s_{1s}$ on $\Om_1$. 

On the pull-back bundle $\hat P\to X_0$ we define the corresponding smooth pull-back connection $\hat A'$ by
setting $\hat A'=A_0'$ on $\hat P|_{X_0'}$ and
$\hat A'=f_1^*A_1'$ on $\hat P|_{B_1'}$. 
We obtain local connection 1-forms for $\hat A'$ over $X_0$ given by 
$\hat\si_{1n}^*\hat A'=f_1^*\si_{1n}^*A_1'$ on the ball $B_1'$,
$\hat\s_1^*\hat A'=\s_1^*A_0'$ on the annulus $\Om_1(N^{-1}\la_1^{1/2},\varrho_0)$, and
$\hat\s_{1s}^*A'=f_1^*\s_{1s}^*A_1'$ on the punctured ball $B_1'\less\{x_1\}$.

On the annulus $\W_1$ we have
$\hat\s_{1s}^*\hat A'=\s_1^*\hat A'=0$, and since
$$
\hat\s_{1s}^*\hat A'=\tilde\rho_1^{-1}\hat\s_1^*\hat A_0'\tilde\rho_1 + \tilde\rho_1^{-1}d\tilde\rho_1
\quad\hbox{on }\W_1, \eqlabel
$$
we see that $d\tilde\rho_1=0$ on $\W_1$ and so $\tilde\rho_1$ is constant on $\Om_1$. The transition function $\hat u_1$ on $B_1'\less\{x_0\}$ is
independent of $\l_1$, since $u_1$ on $X_1\less\{x_{1n},x_{1s}\}$ is
constant along geodesics connecting the north and south
poles. Thus, the bundle $\hat P$ is defined by transition functions which
are constant with respect to $\la_1$. 

We now generalise the preceding discussion to give a construction of
approximately  anti-self-dual connections over multiple connected sums. The
description we give here is closely related to Taubes' iterated gluing
construction [T88, \S 4]. 
The construction parallels the description of the ends of the bubble tree
compactification $M_{X_0,k}(g_0)$ described in Chapter 4. 

It is convenient at this point to introduce some terminology.
Let $I=(i_1,\dots,i_r)$ denote a multi-index of positive integers. The
{\it length} of $I$ is $r$; we regard $0$ as a multi-index of length
zero. 
Given $I=(i_1,\dots,i_r)$, we let $I_-=(i_1,\dots,i_{r-1})$; we
will often denote a multi-index of the form
$(i_1,\dots,i_{r+1})$ by $I_+$ or if we wish to be more specific, by
$Ij$, where $j=i_{r+1}>0$ or $s,n$ (indicating north or south poles of
$\SS^4$), with a slight abuse of notation. 
Let $\script{I}$ be an oriented tree with a finite set of vertices $\{I\}$,
including a base vertex $0$, and a set of edges $\{(I,I_+)\}$. If
$I=(i_1,\dots,i_r)$ and $I=(j_1,\dots,j_t)$, then we say $I<J$ if $r<t$ and
$J=(i_1,\dots,i_r,j_{r+1},\dots,j_t)$. The valence of each vertex $I$ is
the number of edges emanating from that vertex. The height of the
tree $\sI$ is the number of levels -- the length of the longest multi-index
minus one. 
With respect to a given vertex $I$, the edge $(I_-,I)$ is called ${\it
incoming}$ and  
the edge $(I,I_+)$ is called ${\it outgoing}$.

The construction of a $C^\8$, approximately $g$-anti-self-dual connection
$A'$ of second Chern class $k\ge 1$, associated with a tree $\script{I}$,
requires the following data:  
 
\data\thmlabel{\dataGluingData} {\it Gluing data for approximately
anti-self-dual connections.} 
\item{(1)} To each {\it vertex} $I$, we associate a $g_I$-anti-self-dual
connection $A_I$ on a $G$ bundle $P_I\to X_I$ with $c_2(P_I)=k_I\ge 0$. If
$I=0$, then $X_0$ is the base four-manifold with metric $g_0$, while if
$I>0$, then $X_I=\SS^4$ with its standard round metric $g_I\equiv g_1$ of
radius $1$.  

\item{(2)} To each {\it edge} $(I_-,I)$, we associate the data
$(b_I,\la_I,\rho_I,x_I,v_I)$ given by the 
\itemitem{(i)} Connection cutoff parameter $b_I$;
\itemitem{(ii)} Scale parameter $\la_I$;
\itemitem{(iii)} Bundle gluing parameter 
$\rho_I\in\Gl_{x_I}$, where $\Gl_{x_I}=\Hom(P_{I_-}|_{x_I},P_I|_{x_{Is}})$;
\itemitem{(iv)} Centre or gluing site $x_I\in X_{I_-}$;
\itemitem{(v)} Frame $v_I\in FX_0|_{x_I}$ if $I_-=0$.
\item{(3)} Constants $b_0$, $d_0$, $\la_0$, $N$.
\medskip

For convenience, if $I_+=Is$, we denote $b_{Is}=b_I$, $\la_{Is}=\la_I$, $N_{Is}=N$, and $\rho_{Is}=\rho_I$. We let $x_{In}$, $x_{Is}$ denote the north and south poles of the spheres $X_I=\SS^4$. If $I_->0$, then $x_I\equiv\phi_{I_-n}(q_I)\in X_I$, where $q_I\in\RR^4$. Define
$$
\ol b=\max_{I\in\script{I}} b_I\quad\hbox{and}\quad
\ol\la=\max_{I\in\script{I}}\la_I. \eqlabel
$$
The gluing data should satisfy the following constraints:

\cond\thmlabel{\condGluConstraint} {\it Gluing data constraints.} 
\item{(1)} Scales: $4N\la_I^{1/2}\le b_I<\quarter\min\{1,\varrho_0,d_0\}$, $4<N_0\le N$, and $0<\la_I\le\la_0$;
\item{(2)} Separation of centres: Suppose $x_I,x_{I'}\in X_{I_-}$.
\itemitem{(i)} If $I_-=0$, then $\dist_{g_0}(x_I, x_{I'})>4(b_I+b_{I'})$,
\itemitem{(ii)} If $I_->0$, then $|q_I-q_{I'}|>4(b_I+b_{I'})$;
\item{(3)} Topology: $\sum_{I\in\script{I}}k_I=k$ and $k_I>0$ for some $I>0$.
 \medskip

\rmk Definition \ref{\dataGluingData}, together with the constraints of
Condition \ref{\condGluConstraint} should be compared with the definition
of ``bubble tree ideal'' connections in \S 4.3.  
The requirements on the scales and separation of centres are in place
simply to ensure that the different gluing regions do not interfere with
one another. 
\medskip

The gluing procedure now generalises to give a $C^\8$ family of
approximately $g$-anti-self-dual connections $A'=\#_{I\in\script{I}}A_I'$
on a bundle $P$ over a multiple connected sum $X=\#_{I\in\script{I}}X_I'$.
First, consider the definition of coordinate charts, open balls, and annuli in $X_0$. If $I_-=0$, 
let $\phi_I^{-1}=\exp_{v_I}^{-1}: B(x_I,\varrho_0)\to{\Bbb R}^4$ be a geodesic normal coordinate chart
defined by a point $v_I$ in the oriented frame bundle fibre $FX_0|_{x_I}$.
Let $B_I(r)=B(x_I,r)$ be the open geodesic ball in $X_0$ with centre $x_I$ and radius $r$. 

Turning to the four-spheres $X_I$, for any $I>0$,
let $\phi_{I\al}=\phi_\al$, $\a=s,n$ be the standard inverse coordinate charts 
on $X_I$. Define open neighbourhoods in $X_I$ by 
$$\leqalignno
{B_{Is}(r)&=B(x_{Is},r)=\phi_{Is}\left(\{x\in\RR^4:|x|<r\}\right)\quad\hbox{and} &\numeqn\cr
B_{I_+}(r)&=B(x_{I_+},r)=\phi_{In}\left(\{x\in\RR^4:|x-q_{I_+}|<r\}\right). &\cr}
$$
Let $\W_I(r,R)=\W(x_I,r,R)$ be the open annulus $B_I(R)\less\ol B_I(r)$ centred at $x_I\in X_{I_-}$, with inner radius $r$ and outer radius $R$.

Define small balls $B_I'=B(x_I,N\la_I^{1/2})$ and annuli $\Om_I=\Om(x_I,N\la_I^{1/2}, N\la_I^{1/2})$ in $X_{I_-}$, $I>0$. The open subset $X_{I_-}'$ is the complement in $X_{I_-}$ of the balls $\ol{B}_I(N^{-1}\la_I^{1/2})$, the open subset $X_{I_-}''$ is the complement in $X_{I_-}$ of the balls $\ol{B}_I(\half\la_I^{1/2})$, and the open subset $X_{I_-}'''$ is the complement in $X_{I_-}$ of the balls $\ol{B}_I(2N\la_I^{1/2})$
We define identification maps $f_I$ by
$$
f_I= \phi_{In}\circ c_I\circ \phi_I^{-1}: B_I'\too X_I',
\eqlabel\eqnlabel{\eqConfBlwUp}
$$
where $c_I$ is the dilation $x\to x/\l_I$ on
${\Bbb R}^4$. The maps $\phi_I$ above are local coordinate charts on $X_{I_-}$ given by 
$$
\phi_I^{-1}=\cases{\exp_{v_I}^{-1} &if $I_-=0$ \cr
\tau_I\circ\phi_{I_-n}^{-1} &if $I_->0$, \cr} \eqlabel
$$ 
where $\tau_I$ is the translation $x\to x-q_I$ on $\RR^4$. 
The charts $\phi_I^{-1}=\exp_{v_I}^{-1}$ may be replaced by $\bar\phi_I^{-1}=\tau_{p_I}\circ\exp_{v_I}^{-1}$, $|p_I|\ll\varrho_0$, if we wish to compute derivatives with respect to the centres $x_I$ in $X_0$.
For notational consistency, we let $f_0$ denote the identity map on $X_0$.  

Using the diffeomorphisms $f_I:\Om_I\to\Om_{Is}$ we obtain a connected sum $X=\#_{I\in\script{I}}X_I'$. 
We again defer to \S 3.5 for the precise definition of a metric $g$ on $X$ closely approximating the metrics $g_I$ on the summands $X_I'$ and such that the maps $f_I:B_I'\to X_I'$ are conformal. With this choice of metric, the connected sum $(X,g)$ is conformally equivalent to $(X_0,g_0)$.

We have a local section $\si_I$ of $P_{I_-}$ defined by 
a choice of point in the fibre $P_{I_-}|_{x_I}$ and 
$A_{I_-}$-parallel translation from $x_I$; similarly, we have local sections $\si_{In}$, $\si_{Is}$ of $P_I$ defined by 
a choice of points in the fibres $P_I|_{x_{In}}$, $P_I|_{x_{Is}}$ and
$A_I$-parallel translation from $x_{In},x_{Is}$. These sections provide local trivialisations $P_{I_-}|_{B_I(\varrho_0)}\simeq B_I(\varrho_0)\times G$ and 
$P_I|_{X_I\less\{x_{In}\}}\simeq X_I\less\{x_{In}\}\times G$. 
Define $C^\8$ cutoff functions $\psi_I$ on each summand $X_I$ by setting
$$
\psi_I\equiv(\phi_{Is}^{-1})^*\psi_{b_I}
\prod_{I_+}(\phi_{I_+}^{-1})^*\psi_{I_+}
\quad\hbox{on }X_I, \eqlabel\eqnlabel{\eqPsiI}
$$ 
where the factor $(\phi_{Is}^{-1})^*\psi_{b_I}$ is omitted when $I=0$. Note that $\psi_I=0$ on the balls $B_{Is}(b_I/2)$ and $B_{I_+}(b_{I_+}/2)$ in $X_I$ and smoothly extends by $1$ on the complement of the balls $B_{Is}(b_I)$ and $B_{I_+}(b_{I_+})$ in $X_I$. Lastly, extend each $\psi_I$ by zero to give a $C^\8$ cutoff function on the connected sum $X$.
Setting $A_{I_-}'=\psi_{I_-}A_{I_-}$, $A_I'=\psi_IA_I$,
we obtain $C^\8$ almost anti-self-dual connections $A_{I_-}',A_I'$ which are flat on the balls $B_I(b_I/2), B_{Is}(b_I/2)$.

The gluing parameter $\rho_I$ provides an isomorphism of the fibres
$:P_{I_-}|_{x_I}\simeq P_I|_{x_{Is}}$. Using the connections
$A_{I_-},A_I$, this identification is extended to give a bundle
isomorphism $\tilde\rho_I:P_{I_-}|_{\Om_I}\to P_I|_{\Om_{Is}}$ covering
$f_I$. Using these identification maps we obtain a connected-sum $G$ bundle
$P\to X$ with $c_2(P)=k$ and whose transition functions are constant with
respect to the scales $\la_I$. 
The cutoff connections $A_I'$ on $P_I$ patch together to give a $C^\8$
connection $A'$ on $P$. As before, the connection $A'$ on the connected-sum
bundle $P$ over $X$ pull back via the maps $f_I$ to give a connection $\hat
A'$ on a bundle $\hat P$ over $X_0$.

Lastly, we record some estimates for the connections $A'$ when restricted
to a summand $X_I'$. For this and later purposes, we define the following
Sobolev norms:
Let $\]^{g_I}$ denote the Levi-Civita connection on $TX_I$ defined by the
metric $g_I$, so that if $f\in C^\8(X_I)$, then
$$
\|f\|_{L^p_n(X_I,g_I)}=\sum^n_{i=0}\|(\]^{g_I})^if\|_{L^p(X_I,g_I)}, \eqlabel
$$
for any $1\le p\le\8$ and integer $n\ge 0$. Similarly, if $\al\in\W^l(X_I,\ad P_I)$, then
$$
\|\al\|_{L^p_n(X_I,A_I,g_I)}=\sum^n_{i=0}\|(\]^{A_I,g_I})^i\al\|_{L^p(X_I,g_I)}.\eqlabel
$$
It is important to note that these norms will depend only on a set of {\it
fixed} connections, $\{A_I\}_{I\in\sI}$, and a set of {\it fixed} metrics
$\{g_I\}_{I\in\sI}$. 

Recalling that $A_I'=\psi_IA_I$,
define one-forms $a_I\in\Om^1(X_I,\ad P_I)$ by setting
$A_I=A_I'+a_I$. Thus
$$
a_I=\cases{(1-\psi_I)\si_{I_+}^*A_I& on $B_{I_+}(b_{I_+})$, \cr 
           0& on $X_I\less \bigcup_{I_+}B_{I_+}(b_{I_+})$. \cr}
$$
With the aid of bounds for the derivatives of the cutoff functions
$\psi_J$ for $C=C(g_J)$ and $J=I_-$ or $I$,
$$\leqalignno
{|d\psi_J|_{g_J}&\le Cb_J^{-1}\quad\hbox{on }\Om_J(b_J/2,b_J)
\quad\hbox{and}\quad\|d\psi_J\|_{L^4(X_J,g_J)}\le C, 
&\numeqn\eqnlabel{\eqdPsiIEst}\cr}
$$
standard arguments then give the following estimates:

\proclaim{Lemma}{\thmlabel{\lemAIprimeEst}
Let $1\le p<\8$. Then there exists a constant $C=C(A_I,g_I,p)$ such that
\item{\rm (a)} $\|a_I\|_{L^\8(X_I,g_I)} \le C{\ol b}$ and
$\|a_I\|_{L^p(X_I,g_I)} \le C{\ol b}^{4/p+1}$,
\item{\rm (b)} $\|F(A'_I)\|_{L^\8(X_I,g_I)} \le C$ and
$\|F^{+,g_I}(A'_I)\|_{L^p(X_I,g_I)} \le C{\ol b}^{4/p}$.}


\section{Approximate gluing maps} Adopting a more global perspective, the
construction of \S\ref{\secGluCon} yields a family of ``approximate gluing
maps'', $\sJ':\sT/\Ga\to\cB^*_{X,k}$ and $\hat\sJ':\sT/\Ga\to\cB^*_{X_0,k}$,
which we describe in this section. 
We first recall that the standard Kuranishi models give the required
parametrisations for neighbourhoods of points $[A_I]$ in
$M_{X_I,k_I}(g_I)$.  
Let $A_I$ be a $g_I$-anti-self-dual connection over $X_I$, with isotropy group $\Ga_{A_I}$
and $H_{A_I}^2=0$. For a small enough open neighbourhood $T_{A_I}$ of $0\in H^1_{A_I}$, we have smooth $\Ga_{A_I}$-equivariant maps
$$\leqalignno
{&\al_I: T_{A_I}\too \ker d_{A_I}^{*,g_I}\subset\Om^1(X_I,\ad P_I)
&\numeqn\cr}
$$
solving the $g_I$-anti-self-dual equation
$F^{+,g_I}(A_I+\al_I(t_I))=0$, $t_I\in T_{A_I}$.
Setting $A_I(t_I)=A_I+\al_I(t_I)$, we obtain a homeomorphism 
$$
\vth_I:T_{A_I}/\Ga_{A_I}\too U_{A_I},\qquad t_I\longmapsto [A_I(t_I)],
\eqlabel
$$
onto an open neighbourhood $U_{A_I}$ of $[A_I]\in M_{X_I,k_I}$. If $A_I$ is
the product connection, $\Theta$, then 
$\Ga_{A_I}=\SU(2)$ and so $H_{A_I}^0\ne 0$, while $H_{A_I}^1=0$. 
If $A_I$ is a non-trivial reducible connection, then
$\Ga_{A_I}={\Bbb S}^1$ and $H_{A_I}^0\ne 0$: we have a homeomorphism $\vth_I:T_{A_I}/\Ga_{A_I}\to U_{A_I}$ and a diffeomorphism
$\vth_I:(T_{A_I}\less\{0\})/\Ga_{A_I}\to U_{A_I}\less [A_I]$. Finally, if
$A_I$ is irreducible, then $\Ga_{A_I}=(\pm 1)$ and
$H_{A_I}^0=0$: in this case we have a diffeomorphism $\vth_I:T_{A_I}/\Ga_{A_I}\to U_{A_I}$.

We now dispose of the construction of neighbourhoods of reducible connections in $M_{X_0,k}(g_0)$.
Recall that the reducible connections in $M_{X_0,k}(g_0)$ are in one-to-one correspondence with
pairs $\{\pm c\}$, where $c\in H^2(X_0,\ZZ)$ satisfies $c^2=k$. In
particular, there are only finitely many and so to describe a neighbourhood
of any such reducible connection $[A]\in M_{X_0,k}(g_0)$, we may employ the
Kuranishi model $\vth_A:T_A/\Ga_A\to U_A$. 

We now describe the approximate gluing maps $\sJ'$ and $\hat\sJ'$,
beginning with the parameter spaces $\sT/\Ga$.
First, with the centres $\{x_I\}$ and scales $\{\la_I\}$ held fixed, the
parameter spaces $T_{A_I}$ and $\Gl_{x_I}$ combine to give a $C^\8$
manifold 
$$\leqalignno 
{T&\equiv 
T_{A_0}\times\prod_{I\in\script{I}}\left(T_{A_I}\times \Gl_{x_I}\right),
&\numeqn\eqnlabel{\eqT}\cr}
$$
parametrising a ``small'' family of approximately anti-self-dual connections. 
Then 
$$
\Ga\equiv\Ga_{A_0}\times\prod_{I\in\script{I}}\Ga_{A_I}
\eqlabel\eqnlabel{\eqBigStab}
$$ 
acts freely on $T$ and $T/\Ga$ is a $C^\8$ manifold. 
If we allow the centres, now denoted $y_I$, to move over disjoint balls
$B(x_I,r_0)\subset X_{I_-}$ and allow the scales $\l_I$ to vary in the
interval $(0,\la_0)$, the parameter space of Eq. \eqref{\eqT} is augmented
to give a $C^\8$ manifold 
$$\leqalignno
{\script{T}
&\equiv T_{A_0}\times\prod_{I\in\script{I}}
\left(T_{A_I}\times \Gl_{x_I}\times B(x_I,r_0)\times (0,\la_0)\right),
&\numeqn\eqnlabel{\eqBigT}\cr}
$$
parametrising a ``large'' family of approximately $g$-anti-self-dual connections.
Again, $\Ga$ acts freely on $\sT$ and $\sT/\Ga$ is a $C^\8$ manifold. 
We fix local trivialisations of the frame bundle $FX_0$ over the balls
$B(x_I,r_0)$ and these provide smooth families of geodesic normal
coordinate charts on $X_0$. 

We note that the almost anti-self-dual connections $A'$
produced by \S\ref{\secGluCon} are indeed irreducible:

\proclaim{Lemma}{\thmlabel{\lemAprimeIrred}
Let $A'$ be a connection on the $G$ bundle $P$ over $X$ defined by {\rm
Data \ref{\dataGluingData}} and {\rm Condition \ref{\condGluConstraint}}.
Then $A'$ is irreducible, that is $H_{A'}^0=0$, for small enough $b_0$ and
large enough $N_0$.} 

The Lemma follows from Aronszajn's unique continuation principle for
solutions to $\De_{A'}\eta=0$ via standard methods, so the proof is
omitted. Hence, 
the approximate gluing construction of \S\ref{\secGluCon} gives a $C^\8$ map
$$
\script{J}':\sT/\Ga\too \cB_{X,k}^*,\qquad t\longmapsto [A'(t)], 
\eqlabel\eqnlabel{\eqJprime}
$$
where $\cB_{X,k}^*$ has the structure of an $L^2_n$ Hilbert manifold, $n\ge
3$. Moreover, $\sJ'$ is a $C^\8$ submersion
onto its image; see \S 5.2. We refer to $\script{J}'$ as an {\it
approximate gluing map} over $X$ and its image $\sU'\subset\cB_{X,k}^*$ as
an {\it approximate gluing neighbourhood}. 

The dimension of the parameter space $\sT/\Ga$ is given by 
$$
\dim \sT/\Ga = \dim H^1_{A_0}-\dim H^0_{A_0}
+\sum_{I>0}(\dim H^1_{A_I}-\dim H^0_{A_0}+8),
\eqlabel
$$
since each factor $\Gl_{x_I}\times B(x_I,r_0)\times (0,\la_0)$ has
dimension $8$, $\dim H^0_{A_I}=\dim\Ga_{A_I}$, and $H^2_{A_I}=0$ for all
$I\ge 0$ by hypothesis. Families of centred
$g_I$-anti-self-dual connections $A_I\in M^0_{X_I,k_I}(g_I)$ are parametrised by small
balls $T^0_{A_I}$ and thus, we obtain a $C^\8$ parameter
space 
$$\leqalignno
{\script{T}^0
&\equiv T_{A_0}\times\prod_{I\in\script{I}}
\left(T^0_{A_I}\times \Gl_{x_I}\times B(x_I,r_0)\times (0,\la_0)\right),
&\numeqn\eqnlabel{\eqBigTZero}\cr}
$$
with $C^\8$ quotient $\sT^0/\Ga$ of dimension
equal to $\dim M_{X,k}(g)$. The map $\sJ':\sT^0/\Ga\to \cB_{X,k}^*$ is a
$C^\8$ embedding; see \S 5.2. 
 
Lastly, using  
the conformal  diffeomorphisms $f_I$, the bundle $P$ over $X$ pulls back to
a bundle $\hat P$ over $X_0$. 
The gluing construction now produces an approximately $g_0$-anti-self-dual connection
$\hat A'$ in $\cB_{X_0,k}^*$. The map $\script{J}'$ of Eq.
\eqref{\eqJprime} pulls back to a $C^\8$ map 
$$
\hat\script{J}':\sT/\Ga\too \cB_{X_0,k}^*,\qquad t\longmapsto [\hat A'(t)].
\eqlabel\eqnlabel{\eqHatJprime}
$$
Again, $\hat\script{J}'$ is a $C^\8$ submersion onto its image and and is a
$C^\8$  embedding when the parameter space $\sT/\Ga$ is replaced by the
smaller parameter space $\sT^0/\Ga$; see \S5.2. As before, the image $\sV'$
of $\hat\script{J}'$ in $\cB_{X_0,k}^*$ is called an approximate gluing
neighbourhood. 


\section{Metrics on connected sums}
\seclabel{\secMetBubbleTree}
In this section we define a conformal structure $[g]$ on the 
connected sum $X=\#_{I\in\script{I}}X_I$. This is accomplished by replacing
the standard round metric $g_I$ on each spherical summand $X_I'$ by a
quasi-conformally 
equivalent metric $\tilde g_I$ so that the identification maps $f_I:
B_I'\to X_I'$ are conformal. We then construct a $C^\8$ metric $g$ on $X$
in the conformal class $[g]=[g_0]$ and compare the resulting $L^p$ norms
for the different possible metrics on each summand $X_I'$. Our construction
is modelled after the constructions of Donaldson and Taubes for metrics on
connected sums -- see [D86, p. 322], [D-K, p. 293], [D-S], and [T92]. The
metric $g$ depends on the choice of fixed base metric $g_0$, fixed neck
width parameter $N$, scales $\l_I$, centres $x_I$, and frames $v_I$. We
also obtain bounds for the derivatives of $g$ with respect to $\l_I$ and
$x_I$. 

With respect to a geodesic normal coordinate system $x=\phi_{i_1}^{-1}$ on $B_{i_1}(\varrho_0)\subset X_0$, the covariant components of $g_0$ satisfy
$$\leqalignno
{(\phi_{i_1}^*g_0)_{\mu\nu}(0)=\d_{\mu\nu}&\quad\hbox{and}\quad 
\pd{(\phi_{i_1}^*g_0)_{\mu\nu}}{x^\a}(0)=0, &\numeqn\eqnlabel{\eqMetNrmGeoCrd}\cr 
|(\phi_{i_1}^*g_0)_{\mu\nu}-\d_{\mu\nu}|(x)\le c|x|^2 &\quad\hbox{and}\quad 
\left|\pd{(\phi_{i_1}^*g_0)_{\mu\nu}}{x^\a}\right|(x)\le c|x|, 
\qquad |x|<\varrho_0/2, &\cr}
$$
for some $c=c(g_0)$ [K-N].
The analogous relations hold for the contravariant components of $g_0$. 
We now define a conformal structure $[g]$ on $X$:

\defn\thmlabel{\defnAlmostRndMet} The conformal structure $[g]$ on $X$ is defined by the $C^\8$ metric $g_0$ on $X_0'$ and a choice 
of $C^\8$ metric $\tilde g_I$ on 
each summand $X_I'$, $I>0$, given by 
$$
(\phi_{In}^*\tilde g_I)_{\mu\nu}(x)
\equiv\cases{h_1^2(x)(\phi_I^*g_0)_{\mu\nu}(\l_Ix) &if $I_-=0$\cr
h_1^2(x)h_1^{-2}(\l_Ix+q_I)(\phi_I^*\tilde g_{I_-})_{\mu\nu}(\l_Ix)
&if $I_->0$,\cr}
$$
where $|x|<N\l_I^{-1/2}$. 
For convenience, we let $g_I\equiv g_1$ denote the standard metric on $X_I$ and let $\tilde g_0\equiv g_0$ denote the metric on $X_0$.

Definition \ref{\defnAlmostRndMet} provides the following expression for $\tilde g_I$:
$$\leqalignno
{(\phi_{In}^*\tilde g_I)_{\mu\nu}(x)
&=h_1^2(x)(\phi_{i_1}^*g_0)_{\mu\nu}(y(x)),\qquad |x|<N\l_I^{-1/2}, &\numeqn\cr}
$$
where 
$$\leqalignno
{y(x)&=\phi_{i_1}^{-1}\circ f_{i_1}^{-1}\circ\cdots\circ f_I^{-1}\circ\phi_{In}(x) &\numeqn\cr
&=\l_{i_1}(\l_{i_1i_2}(\cdots(\l_{I_-}(\l_Ix+q_I)+q_{I_-})\cdots)+q_{i_1i_2}).
&\cr}
$$
The map $f_I:B_I'\to X_I'$ is now conformal with respect 
to the metrics $\tilde g_{I_-}$ on $B_I'\subset X_{I_-}'$ and 
$\tilde g_I$ on $X_I'$:
$$\leqalignno
{(\phi_I^*f_I^*\tilde g_I)_{\mu\nu}(x)
&=\cases{\l_I^{-2}h_1^2(x/\l_I)(\phi_I^*\tilde g_0)_{\mu\nu}(x) &if $I_-=0$\cr
\l_I^{-2}h_1^2(x/\l_I)h_1^{-2}(x+q_I)(\phi_I^*\tilde g_{I_-})_{\mu\nu}(x)
&if $I_->0$,\cr}
&\cr}
$$
where $|x|<N\l_I^{1/2}$. Thus, $f_I^*\tilde g_I$ is conformally equivalent to the metric $g_{I_-}$ on $\W_I$ and so we obtain a conformal structure $[g]$ on $X=\#_{I\in\script{I}} X_I$.  

We must verify that $\tilde g_I$ is a good approximation to the standard round metric $g_I$ on $X_I'$ for $\la_{i_1}$ small. 

\proclaim{Lemma}{\thmlabel{\lemMetXIEst}
For any $I>0$,
the metric $\tilde g_I$ converges to $g_I$ in $C^\8$ on compact subsets of $X_I\less\{x_{Is}\}$ as $\la_{i_1}\to 0$. Moreover, we have the following bounds:
\item{\rm (a)} For any integer $l\ge 0$, there is a constant $c=c(g_0,l)$ such that 
$$
\left|{\rd^l (\phi_{In}^*\tilde g_I)_{\mu\nu}\over 
\rd x^{\al_1}\cdots \rd x^{\al_l}}
-{\rd^l (\phi_{In}^*g_I)_{\mu\nu}\over 
\rd x^{\al_1}\cdots \rd x^{\al_l}}\right|\le cN^2\la_{i_1}h_1^2(x),
\qquad |x|<N\l_I^{-1/2}.
$$
The analogous bounds hold for the contravariant components $(\phi_{In}^*\tilde g_I)^{\mu\nu}$, provided $h_1^2(x)$ is replaced by $h_1^{-2}(x)$.
\item{\rm (b)} Let $*_{\tilde g_I}$ denote the Hodge star operator for $\tilde g_I$. Then there is a constant $c=c(g_0)$ such that 
$$
\|*_{\tilde g_I}\z-*_{g_I}\z\|_{L^\8(X_I',g_I)}
\le cN^2\l_{i_1}\|\z\|_{L^\8(X_I',g_I)}, \qquad \zeta\in\Om^2(X_I',\ad P_I).
$$}

\pf (a) This follows easily from Eq. \eqref{\eqMetNrmGeoCrd} and Definition \ref{\defnAlmostRndMet}. 
\noindent (b) This follows immediately from (a) and the definition of the Hodge star operator. \qed

We will also require bounds for the derivatives of $\tilde g_J$ with respect to
the scales $\l_I$ and centres $x_I$. The following estimates will suffice
for our application. 

\proclaim{Lemma}{\thmlabel{\lemDMetDLaqEst} If $0<I\le J$, there is a constant $c=c(g_0,J)$ such that the following bounds hold.
\item{\rm (a)} For any $|x|<N\l_J^{-1/2}$,
$$
\left|\pd{(\phi_{Jn}^*\tilde g_J)_{\mu\nu}}{\l_I}\right|(x)
\le \cases{c\la_{i_1}h_1^2(x)& if $I<J$ and $|I|=1$ \cr
c\la_{i_1}^2h_1^2(x)& if $I<J$ and $|J|\ge 2$ \cr
cN^2h_1^2(x)& if $I=J$ and $|J|=1$ \cr
c\la_{i_1}^2|x|h_1^2(x) \le cN\la_{i_1}^2\l_I^{-1/2}h_1^2(x) 
& if $I=J$ and $|J|\ge 2$. \cr}
$$
\item{\rm (b)} If $\rd/\rd p_I\equiv p_I^\al\rd/\rd q_I^\al$, then for any $|x|<N\l_J^{-1/2}$,
$$
\left|\pd{(\phi_{Jn}^*\tilde g_J)_{\mu\nu}}{p_I}\right|(x)
\le \cases{cN\la_{i_1}^{1/2}h_1^2(x) & if $I=J$ and $|J|=1$ \cr
c\la_{i_1}^2h_1^2(x)& if $I\le J$ and $|J|\ge 2$. \cr}
$$
The analogous bounds in (a) and (b) hold for the contravariant components of $\tilde g_J$, if $h_1^2(x)$ is replaced by $h_1^{-2}(x)$.
\item{\rm(c)} For any $\zeta\in\Om^2(X_J',\ad P_J)$, then
$$\leqalignno
{\left\|\pd{*_{\tilde g_J}}{\l_I}\zeta\right\|_{L^\8(X_J',g_J)}
&\le \cases{cN\l_I^{-1/2}\|\z\|_{L^\8(X_J',g_J)}&if $I=J$ and $|J|\ge 2$\cr
cN^2\|\z\|_{L^\8(X_J',g_J)}& otherwise, \cr} &\cr
\left\|\pd{*_{\tilde g_J}}{p_I}\zeta\right\|_{L^\8(X_J',g_J)}
&\le cN\l_{i_1}^{1/2}\|\z\|_{L^\8(X_J',g_J)}. &\cr}
$$}

\pf (a) The inequalities follow from Eq. \eqref{\eqMetNrmGeoCrd} and Definition \ref{\defnAlmostRndMet}. 

\noindent (b) The proof is similar. When $|I|=1$, we recall that the normal geodesic chart $\phi_{i_1}\equiv\exp_{v_{i_1}}$ is replaced by $\bar\phi_{i_1}\equiv \exp_{v_{i_1}}\circ\tau_{q_{i_1}}$ in order to compute the required derivative at $q_{i_1}=0$ (corresponding to $x_{i_1}=\phi_{i_1}(0)$). 
The estimates follow immediately from (a) and (b). \qed

We next define an honest $C^\8$ metric $g$ on $X$. Consider a neck
$\Om_I=f_I^{-1}(\Om_{Is})$ labelled by the multi-index $I$. 
We replace the metric $\tilde g_{I_-}$ on the annulus $\W_I$ and
replace the metric $\tilde g_I$ on the annulus $\W_{Is}$ by conformally equivalent metrics 
$m_{I_-}\tilde g_{I_-}$ and
$m_I\tilde g_I$ so that
$$
m_{I_-}g_{I_-}=f_I^*(m_Ig_I)\quad\hbox{on }\W_I. \eqlabel
$$ 
Hence, the metrics $m_{I_-}\tilde g_{I_-}$ and $m_I\tilde g_I$ agree on the neck and patch together to give a $C^\8$ metric, say $g$, on a neighbourhood of the neck in the connected sum $X_{I_-}\# X_I$. On the annulus 
$\W_I=\phi_I(\{x\in\RR^4:N^{-1}\l_I^{1/2}<|x|<N\l_I^{1/2}\})$ we have 
$$\leqalignno
{(\phi_I^*f_I^*\tilde g_I)_{\mu\nu}(x)
&=\cases{4\l_I^2(\l_I^2+|x|^2)^{-2}
(\phi_I^*\tilde g_0)_{\mu\nu}(x) &if $I_-=0$\cr
4\l_I^2(\l_I^2+|x|^2)^{-2}
h_1^{-2}(x+q_I)(\phi_I^*\tilde g_{I_-})_{\mu\nu}(x) &if $I_->0$.\cr} 
&\numeqn\cr}
$$ 
By comparing $f_I^*\tilde g_I$ and $g_{I_-}$ on $\Om_I$, a little experimentation reveals that the $C^\8$ conformal factors $m_{I_-}$ and $m_I$ can be chosen so that
$$\leqalignno
{\k^{-1}\le m_{I_-}\le \k N^4
&\quad\hbox{on }\W_I(N^{-1}\l_I^{1/2},N\l_I^{1/2}), &\numeqn\cr
\k^{-1}\le m_{I_-}\le \k 
&\quad\hbox{on }\W_I(\half\l_I^{1/2},N\l_I^{1/2}), &\cr
m_{I_-}=1 
&\quad\hbox{on }\W_I(2\l_I^{1/2},4N\l_I^{1/2}), &\cr}
$$
and likewise for $m_I$ on $\Om_{Is}$, for some constant $\ka=\ka(g_0)$.
For each summand $X_I$, we smoothly extend the $m_I$ to $X_I'$ by setting $m_I\equiv 1$ away from the neck regions. This gives a $C^\8$ metric $g$ on $X=\#_{I\in\sI}X_I$ by setting 
$$
g\equiv m_I\tilde g_I\quad\hbox{on } X_I',\quad\hbox{for all }I\in\sI. \eqlabel
$$
The construction ensures that each $m_I$ obeys
$$
\k^{-1}\le m_I\le \k N^4\quad\hbox{on }X_I',
\quad \k^{-1}\le m_I\le \k\quad\hbox{on }X_I'',
\quad\hbox{and }m_I=1\quad\hbox{on }X_I'''. \eqlabel\eqnlabel{\eqMIBound}
$$
Thus, the metrics $\tilde g_I$ and $g$ 
are equivalent on $X_I''$ with constants independent of $N$, and equivalent over $X_I'$ with constants now depending on $N$.

The Hodge star operator $*_g:\W^2(X,\ad P)\to \W^2(X,\ad P)$ only depends on the conformal class $[g]$ of $g$ and so over each summand $X_I'$ of $X$ we have $*_g=*_{m_I\tilde g_I}=*_{\tilde g_I}$. From Lemma \ref{\lemDMetDLaqEst}, we obtain:

\proclaim{Lemma}{\thmlabel{\lemDHodgeStarDLaq} There is a constant $c=c(g_0)$ such that for any $\zeta\in\Om^2(X,\ad P)$, then
\item{\rm (a)}
$\|\pd{*_g}{\l_I}\zeta\|_{L^\8(X,g)}
\le cN\l_I^{-1/2}\|\z\|_{L^\8(X,g)}$, 
\item{\rm (b)}
$\|\pd{*_g}{p_I}\zeta\|_{L^\8(X,g)}
\le cN\l_{i_1}^{1/2}\|\z\|_{L^\8(X,g)}$.}

We will often need to compare $L^p$ norms defined by the different metrics
$g_I$, $\tilde g_I$, and $g$ over $X_I'\subset X$. The required ``comparison
estimates'' are below follow in a straightforward way from Lemma
\ref{\lemMetXIEst} and Eq. \eqref{\eqMIBound}, and
similar inequalities may be found in [D-K, p. 294].

\proclaim{Lemma}{\thmlabel{\lemMetComparEst} For any $I\ge 0$, the
following holds.  
\item{\rm (a)} If $2\le p<\8$ and $4\le q<\8$, there is a constant $c=c(g_0,k,p,q)$, $1\le c<\8$, such that for any $\w\in\W^1(X_I',\ad P_I)$ and $\z\in\W^2(X_I',\ad P_I)$, then
$$\leqalignno
{\|\w\|_{L^q(X_I',g)}\le c\|\w\|_{L^q(X_I',g_I)}&\quad\hbox{and}\quad 
\|\z\|_{L^p(X_I',g)}\le c\|\z\|_{L^p(X_I',g_I)}, &\cr
\|\w\|_{L^q(X_I'',g_I)}\le c^{-1}\|\w\|_{L^q(X_I'',g)}&\quad\hbox{and}\quad
\|\z\|_{L^p(X_I'',g_I)}\le c^{-1}\|\z\|_{L^p(X_I'',g)}.&\cr}
$$
\item{\rm (b)} If $1\le p<\8$, $n\ge 1$, and $b_0$ is sufficiently small, there is a constant $c=c(g_0,k,n,N,p)$, $1\le c<\8$, such that for any $\a\in\W^n(X_I',\ad P_I)$, then
$$
c^{-1}\|\a\|_{L^p(X_I',g_I)}\le 
\|\a\|_{L^p(X_I',\tilde g_I)},\quad \|\a\|_{L^p(X_I',g)}
\le c\|\a\|_{L^p(X_I',g_I)}.
$$}

Lastly, having defined the conformal structure $[g]$ of $X$, we apply
the estimates for $d\psi_I$ in Eq.
\eqref{\eqdPsiIEst}, the estimates for $A_I'$ and $F(A_I')^{+,g_I}$ in
Lemma \ref{\lemAIprimeEst}, and the estimates for $*_g-*_{g_I}$ in Lemma
\ref{\lemMetXIEst} to obtain a bound for the 
$L^p$-norm of the $g$-self-dual curvature 
$F^{+,g}(A')=\half(1+*_g)F(A')$ of the connection $A'$ on the connected sum
bundle $P$ over $X$. Similar estimates have been given by Taubes and Donaldson.

\proclaim{Proposition}{\thmlabel{\propFAprimeEst}
For $1\le p<\8$ and sufficently small $b_0$, there exists a constant
$C=C(g_0,p,\sT)$ such that for any $t\in \sT$ one has
$\|F^{+,g}(A')\|_{L^p(X,g)}\le C{\ol b}^{4/p}$.}


\section{Estimates over connected sums and conformal vector fields}
The goal of this section is to obtain $L^2$ estimates for the derivatives
with respect to the scales $\l_I$ and centres $x_I$ of $\ad\hat P$-valued
one-forms $\hat\om$ over the base manifold $X_0$ obtained by pulling back
$\ad P$-valued one-forms $\om$ over the connected sum $X$. 

Following Taubes
[T84b], [T88], let us begin by defining some useful Sobolev norms on
$\Om^1(\SS^4,\ad P)$ and examine their behaviour under conformal
diffeomorphisms.   
Suppose $A$ is a $C^\8$ connection on a $G$ bundle $P$ over $\SS^4$. Let  
$g_1$ be the standard round metric on $\SS^4$ and let $\de$ be the flat
metric on $\SS^4\less\{s\}$ obtained via the conformal identification
$\phi_n^{-1}:\SS^4\less\{s\}\to \RR^4$. Let $\]^{A,g_1}$ denote the
covariant derivative on $\Om^1(\SS^4,\ad P)$ defined by the connection $A$
and metric $g_1$, while $\]^{A,\de}$ denotes the covariant derivative on
$\Om^1(\SS^4\less\{s\},\ad P)$ defined by $A$ and $\de$. 
Define an $L^2_1$ norm on $\Om^1(\SS^4,\ad P)$ by
$$
\|\om\|_{L^2_1(\SS^4,A,g_1)}\equiv\|\om\|_{L^2(\SS^4,g_1)}
+\|\]^{A,g_1}\om\|_{L^2_1(\SS^4,g_1)}. \eqlabel
$$
Similarly, if $\om$ has compact support in $\SS^4\less\{s\}$, define 
$$\leqalignno
{|\om|_A&\equiv\|\]^{A,\de}\om\|_{L^2(\SS^4,\de)} &\cr
\|\om\|_{L^2_1(\SS^4,A,\de)}&\equiv\|\om\|_{L^2(\SS^4,\de)}
+\|\]^{A,\de}\om\|_{L^2(\SS^4,\de)}. &\numeqn\cr}
$$
The properties of $|\cdot|_A$ and $\|\cdot\|_{L^2_1(\SS^4,A,\de)}$ are
described by the following result of [T84b]. Recall that ${\bf C}={\bf
D}\times{\bf T}\times\SO(4)$ is identified, using
$\phi_n:\RR^4\to\SS^4\less\{s\}$, with the subgroup of conformal
diffeomorphisms of $(\SS^4,g_1)$ which fix the south pole. 

\proclaim{Lemma}{{\rm [T84b, Proposition 2.4]}\thmlabel{\lemTaubesCEst} Given an $L^2_1$ connection $A$ on a $G$ bundle $P$ over $\SS^4$, then the following holds. 
\item{\rm (a)} $|\cdot|_A$ extends to a continuous norm on $L^2_1\Om^1(\SS^4,\ad P)$. 
\item{\rm (b)} The norm $|\cdot|_A$ is ${\bf C}$-invariant: for any $f\in{\bf C}$, $|f^*\om|_{f^*A}=|\om|_A$.
\item{\rm (c)} There exists a constant $1\le z<\8$, which is independent of $P$, $A$, $f$, and $\om\in\Om^1(\SS^4,\ad P)$, such that 
$$\leqalignno
{z^{-1}\|\om\|_{L^2_1(\SS^4,A,g_1)}
&\le |\om|_A\le z\|\om\|_{L^2_1(\SS^4,A,g_1)}, &\cr
z^{-1}\|\om\|_{L^2_1(\SS^4,A,g_1)}
&\le \|\om\|_{L^2_1(\SS^4,A,\de)}\le z\|\om\|_{L^2_1(\SS^4,A,g_1)}. &\cr}
$$}

\proclaim{Lemma}{{\rm [T88, Lemma 3.1]}\thmlabel{\lemTaubesConfEst}
Let $A$ be a $C^\8$ connection on a $G$ bundle $P$ over $\SS^4$ with its standard metric $g_1$ and let $f:\SS^4\to\SS^4$ be a conformal diffeomorphism. Then there exists a constant $1\le z<\8$, which is independent of $P$, $A$, $f$, and $\om\in\Om^1(\SS^4,\ad P)$, with the following significance:
$$
z^{-1}\|\om\|_{L^2_1(\SS^4,A,g_1)}
\le \|f^*\om\|_{L^2_1(\SS^4,f^*A,g_1)}\le z\|\om\|_{L^2_1(\SS^4,A,g_1)}.
$$}

Recall that $c_\la$ denotes both the dilation $x\mapsto x/\la$ of $\RR^4$
and the conformal diffeomorphism of $(\SS^4,g_1)$ induced by $\phi_n$. 
A straightforward application of H\"older's inequality yields the following
``transfer estimates'' for the maps $c_\la$.  

\proclaim{Lemma}{\thmlabel{\lemXferEst}
Let $2\le p\le p_1\le 4$, let $\l\in(0,1]$, and let $U$ be an open subset of $\SS^4\less B(s,N\la^{1/2})$. Let $P$ be a $G$ bundle over $\SS^4$.
Then there is a constant $C=C(N)$ such that the following holds.
\item{\rm (a)} 
If $\w\in\W^1(U,\ad P)$, then
$\|c_\la^*\w\|_{L^p(c_\la^{-1}(U),g_1)}
\le C\la^{2/p-2/p_1}\|\w\|_{L^{p_1}(U,g_1)}$.
\item{\rm (b)} If $\zeta\in\W^2(U,\ad P)$, then
$\|c_\la^*\z\|_{L^2(c_\la^{-1}(U),g_1)}
\le C\|\z\|_{L^2(U,g_1)}$.}

We next consider the action of the conformal group on $\Om^1(\SS^4,\ad P)$.
Let $f_{\la,q}$ denote the lift to $\SS^4$, via the chart $\phi_n$, of the
conformal diffeomorphism $c_\la\circ\tau_q$ on $\RR^4$. Let $P$ be a $G$
bundle over $\SS^4$ and suppose 
$\w\in\W^1(\SS^4,\ad P)$. Then Eq. \eqref{\eqDctauStarOmDLaq} gives
$$
\pd{f_{\la,q}^*\w}{\l} = -{1\over\la}f_{\la,q}^*\sL_\rad\w\quad\hbox{and}\quad 
\pd{f_{\la,q}^*\w}{p} = -{1\over\la}f_{\la,q}^*\sL_\trans\w,
\eqlabel\eqnlabel{\eqDfStarOmDLaq}
$$
where $\rd/\rd p\equiv p^\mu\rd/\rd q^\mu$.
It will be convenient to express the above Lie derivatives in terms of covariant derivatives. If $A$ is a $C^\8$ connection on $P$, then Eqs. \eqref{\eqLieRadOm} and \eqref{\eqLieTransOm} imply 
$$
\sL_\rad\w=\om+\]^{A,\d}_\rad\w\quad\hbox{and}\quad
\sL_\trans\w=\]^{A,\de}_\trans\w. 
\eqlabel\eqnlabel{\eqLieOmCov}
$$
This leads to the following estimates for the derivatives of $f_{\la,q}^*\w$ with respect to $\l$ and $q$.

\proclaim{Lemma}{\thmlabel{\lemDfstarOmDLaq} Let $A$ be a $C^\8$ connection on a $G$ bundle $P$ over $\SS^4$, let $U\subset \SS^4\less B(s,N\la^{1/2})$ be an open subset, and let $\w\in\W^1(U,\ad P)$. Let $\rd/\rd p= p^\mu\rd/\rd q^\mu$, $|p|\le 1$.
Then there is a constant $C=C(q,N)$ such that the following bounds hold. 
\item{\rm (a)}
$\|\rd{f_{\la,q}^*\w}/\rd{\la}\|_{L^2(f_{\la,q}^{-1}(U),g_1)} 
\le C\la^{-1/2}\|\w\|_{L^2_1(U,A,g_1)}$;
\item{\rm (b)}
$\|\rd{f_{\la,q}^*\w}/\rd{p}\|_{L^2(f_{\la,q}^{-1}(U),g_1)} 
\le C\|\w\|_{L^2_1(U,A,g_1)}$.}

\pf (a) Observe that $U=\phi_n(B(0,N\la^{-1/2}))$ and 
$f_{\la,q}^{-1}(U)=\phi_n(B(q,N\la^{1/2}))$. 
From Eqs. \eqref{\eqDfStarOmDLaq} and \eqref{\eqLieOmCov}, we have
$$
\pd{f_{\la,q}^*\w}{\la}=-\la^{-1}f_{\la,q}^*\sL_\rad\w\quad\hbox{and}\quad
f_{\la,q}^*\sL_\rad\w=f_{\la,q,*}^{-1}\rad\lrcorner f_{\la,q}^*\]^{A,\de}\w
+f_{\la,q}^*\w\quad\hbox{on }f_{\la,q}^{-1}(U),
$$
where $\rad=y^\mu\rd/\rd y^\mu$ and $f_{\la,q,*}^{-1}\rad=x^\mu\rd/\rd x^\mu$ with respect to the coordinates $y=\phi_n^{-1}$ on $U$ and $x=\tau_q\circ\phi_n^{-1}$ on $f_{\la,q}^{-1}(U)$. Since $|f_{\la,q,*}^{-1}\rad|_{g_1}\le C\la^{1/2}$ on $f_{\la,q}^{-1}(U)$, Lemma \ref{\lemXferEst} implies 
$$\leqalignno
{\|f_{\la,q}^*\sL_\rad\w\|_{L^2(f_{\la,q}^{-1}(U),g_1)}
&\le \|f_{\la,q}^*\w\|_{L^2(f_{\la,q}^{-1}(U),g_1)}
+C\la^{1/2}\|f_{\la,q}^*\]^{A,\de}\w\|_{L^2(f_{\la,q}^{-1}(U),g_1)} &\cr
&\le C\la^{1/2}\|\w\|_{L^4(U,g_1)}
+C\la^{1/2}\|\]^{A,\de}\w\|_{L^2(U,g_1)} &\cr
&= C\la^{1/2}\|\w\|_{L^2_1(U,A,\de)}, &\cr}
$$
the last step following by conformal invariance.
Lemma \ref{\lemTaubesCEst} then gives (a). 

\noindent (b) From Eqs. \eqref{\eqDfStarOmDLaq} and \eqref{\eqLieOmCov}, we have
$$
\pd{f_{\la,q}^*\w}{p_I}=-\la^{-1}f_{\la,q}^*\sL_\trans\w\quad\hbox{and}\quad
f_{\la,q}^*\sL_\trans\w=f_{\la,q,*}^{-1}\trans\lrcorner f_{\la,q}^*\]^{A,\de}\w
\quad\hbox{on }f_{\la,q}^{-1}(U),
$$
where $\trans=p^\mu\rd/\rd y^\mu$ on $U$ and $f_{\la,q,*}^{-1}\trans=\la p^\mu\rd/\rd x^\mu$ on $f_{\la,q}^{-1}(U)$. Since $|f_{\la,q,*}^{-1}\trans|_{g_1}\le C\la$ on $f_{\la,q}^{-1}(U)$, Lemma \ref{\lemXferEst} implies 
$$\leqalignno
{\|f_{\la,q}^*\sL_\trans\w\|_{L^2(f_{\la,q}^{-1}(U),g_1)}
&\le C\la\|f_{\la,q}^*\]^{A,\de}\w\|_{L^2(f_{\la,q}^{-1}(U),g_1)} &\cr
&\le C\la\|\]^{A,\de}\w\|_{L^2(U,\de)}
\le C\la\|\w\|_{L^2_1(U,A,\de)}. &\cr}
$$
Lemma \ref{\lemTaubesCEst} then gives (b). \qed

We will frequently need to compute estimates for
families of one-forms $\om$ over connected sums $X$ and to this end, it
will be useful to define suitable Sobolev norms which depend only on the
fixed connections $A_I$ and, in particular, the fixed metrics $g_I$ on each
summand $X_I$ rather than varying metric $g$ on $X$.
Let $P$ be the $G$ bundle over the connected sum $X=\#_{I\in\sI}X_I$
defined in \S\ref{\secGluCon}. Then we may view any $\om\in\Om^1(X,\ad P)$
as a collection of $\om_I\in\Om^1(X_I',\ad P_I)$ which agree over the necks
$\Om_I=f_I^{-1}(\Om_{Is})$ connecting each pair $X_{I_-}$ and $X_I$: 
$$
\si_I^*\om_{I_-}=\Ad(\rho_I^{-1})f_I^*\si_{Is}^*\om_I
\quad\hbox{on }\Om_I,
$$
where $f_I:\Om_I\to\Om_{Is}$ is the identification map. 

From \S\ref{\secMetBubbleTree}, we recall that there is a $C^\8$ metric $g$
on $X$ which agrees, modulo the conformal factors $m_I$, with the metrics
$g_0$ on the base $X_0'$ and $\tilde g_I\simeq g_1$ on the four-spheres
$X_I'$. Moreover, the $L^q$ norms on $\Om^1(X_I',\ad P_I)$, $4\le q<\8$,
and $L^p$ norms on $\Om^2(X_I',\ad P_I)$, $2\le p<\8$, compare uniformly
when defined with the metrics $g_I$, $\tilde g_I$, or $g=m_I\tilde g_I$ on
$X_I'$. The constants involved in these norm comparisons are independent of
the scale parameters $\la_J$ for forms supported on $X_I'$ and independent
of both the $\la_J$ and $N$ for forms supported on $X_I''$. Thus, we may
conveniently define $L^q$ norms on $\Om^1(X,\ad P)$, $4\le q<\8$, and $L^p$
norms on $\Om^2(X,\ad P)$, $2\le p<\8$, using the metric $g$ on $X$. 

In Chapter 5, we will need to bound the $L^2_1$ norms of solutions $\om\in \Om^1(X,\ad P)$ to the $g$-anti-self-dual equation $F^{+,g}(A'+\om)=0$ over $X$.
Unfortunately, since the conformal factors $m_I$ have badly behaved derivatives over the neck regions, the norm comparisons described above do not hold for $L^2_n$ Sobolev norms if $n\ge 1$. Of course, problems of this type are encountered in [D86], [D-K], and [T92]. 
So, given such an $\om\in\Om^1(X,\ad P)$, with $\om=\{\om_I\}_{I\in\sI}$ as above,  and $1\le p<\8$, define 
$$
\|\om\|_{\sL^p_1(X)}\equiv \sum_{I\in\sI}\|\om_I\|_{L^p_1(X_I,A_I,g_I)},
\eqlabel
$$
by analogy with Eq. (6.25) in [T92]. 

Recall that a one-form $\w\in\W^1(X,\ad P)$ pulls back to a one-form
$\hat\w\in\W^1(X_0,\ad\hat P)$ defined by 
$$
\hat\w=f_0^*\cdots f_J^*\w\quad\hbox{on }
f_0^{-1}\cdots f_J^{-1}(X_J')\subset X_0, \eqlabel
$$
for each $J\in\script{I}$. We will need
estimates for the derivatives of $\hat\w$ with respect to the scales
$\la_I$ and centres $x_I$. To begin, we need suitable expressions for these
derivatives: 

\proclaim{Lemma}{\thmlabel{\lemDffStarOmDLaq}
Let $\w\in\W^1(X_J',\ad P_J)$, let $0<I\le J$, and let $\rd/\rd p_I=
p_I^\mu\rd/\rd q_I^\mu$. Then: 
\item{\rm (a)}
$\pd{}{\l_I} f_0^*\cdots f_J^*\w
=f_0^*\cdots f_J^*\pd{\w}{\l_I}$, for $J<I$;
\item{\rm (b)}
$\pd{}{\l_I}f_0^*\cdots f_I^*\w
= -\l_I^{-1}f_0^*\cdots f_I^*\sL_\rad\w$, for $J=I$;
\item{\rm (c)}
$\pd{}{\l_I} f_0^*\cdots f_J^*\w
= -\l_I^{-1}f_0^*\cdots f_I^*\sL_\rad f_{I_+}^*\cdots f_J^*\w$, 
for $J>I$;
\item{\rm (d)}
$\pd{}{p_I}f_0^*\cdots f_J^*\w
=f_0^*\cdots f_J^*\pd{\w}{p_I}$, for $J<I$;
\item{\rm (e)}
$\pd{}{p_I}f_0^*\cdots f_I^*\w
= -\l_I^{-1}f_0^*\cdots f_I^*\sL_\trans\w$, for $J=I$;
\item{\rm (f)}
$\pd{}{p_I} f_0^*\cdots f_J^*\w
= -\l_I^{-1}f_0^*\cdots f_I^*\sL_\trans f_{I_+}^*\cdots f_J^*\w$,
for $J>I$.}

\rmk When $I_-=0$, then $\rd/\rd p_I=p_I^\mu\rd/\rd p_I^\mu$ and
$f_I=\phi_{In}\circ c_I\circ\phi_I^{-1}$ is replaced by ${\ol
f}_I=\phi_{In}\circ c_I\circ\tau_{p_I}\circ\phi_I^{-1}$ in order to compute
the derivative at $p_I=0$.  
\medskip

These expressions lead to the following bounds for the derivatives with
respect to the scales $\l_I$ and centres $x_I$ of the pull-backs
$f_0^*\cdots f_J^*\w$. 

\proclaim{Lemma}{\thmlabel{\lemDffstarOmDLaqEst}
Let $\w\in\W^1(X_J',\ad P_J)$, let $U= f_0^{-1}\circ\cdots\circ
f_J^{-1}(X_J')\subset X_0$, let $0<I\le J$, and let $\rd/\rd p_I=
p_I^\mu\rd/\rd q_I^\mu$ with $|p_I|\le 1$. Then there is a constant
$C=C(g_0,N)$ such that the following holds. 
\item{\rm (a)}
$\left\|\pd{}{\l_I}f_0^*\cdots f_J^*\w\right\|_{L^2(U,g_0)}
\le C\left\|\pd{\w}{\l_I}\right\|_{L^2(X_J',g_J)}$,
for $J<I$;
\item{\rm (b)}
$\left\|\pd{}{\l_I}f_0^*\cdots f_J^*\w\right\|_{L^2(U,g_0)} 
\le C\l_I^{-1/2}\|\w\|_{L^2_1(X_J',A_J,g_J)}$,
for $J\ge I$;
\item{\rm (c)}
$\left\|\pd{}{p_I}f_0^*\cdots f_J^*\w\right\|_{L^2(U,g_0)}
\le C\left\|\pd{\w}{p_I}\right\|_{L^2(X_J',g_J)}$,
for $J<I$;
\item{\rm (d)}
$\left\|\pd{}{p_I}f_0^*\cdots f_J^*\w\right\|_{L^2(U,g_0)} 
\le C\|\w\|_{L^2_1(X_J',A_J,g_J)}$,
for $J\ge I$.}

\pf (a) By repeatedly applying Lemma \ref{\lemXferEst}, we find that
$$\leqalignno
{\left\|\pd{}{\l_I}f_0^*\cdots f_J^*\w\right\|_{L^2(U,g_0)}
&=\left\|f_0^*\cdots f_J^*\pd{\w}{\l_I}\right\|_{L^2(U,g_0)} 
\le C\left\|\pd{\w}{\l_I}\right\|_{L^2(X_J',g_J)}, &\cr}
$$
as required for (a).
For $J=I$ and $U=f_0^{-1}\circ\cdots\circ f_I^{-1}(X_I')\subset X_0$,
Lemmas \ref{\lemXferEst} and \ref{\lemDfstarOmDLaq} show that 
$$\leqalignno
{\left\|\pd{}{\l_I}f_0^*\cdots f_I^*\w\right\|_{L^2(U,g_0)}
&=\left\|f_0^*\cdots f_{I_-}^*\pd{f_I^*\w}{\l_I}\right\|_{L^2(U,g_0)} &\cr
&\le C\left\|\pd{f_I^*\w}{\l_I}\right\|_{L^2(f_I^{-1}(X_I'),g_{I_-})} 
\le C\l_I^{-1/2}\|\w\|_{L^2_1(X_I',A_I,g_I)}. &\cr}
$$
Let $V=f_{I_+}^{-1}\circ\cdots\circ f_J^{-1}(X_J')\subset X_I'$, so that
$U=f_0^{-1}\circ\cdots\circ f_I^{-1}(V)\subset X_0$. Then
for $J>I$, we have
$$\leqalignno
{\left\|\pd{}{\l_I}f_0^*\cdots f_J^*\w\right\|_{L^2(U,g_0)}
&=\left\|f_0^*\cdots f_{I_-}^*\pd{}{\l_I}
f_I^*\cdots f_J^*\w\right\|_{L^2(U,g_0)} &\cr
&\le C\left\|\pd{}{\l_I}f_I^*f_{I_+}^*\cdots 
f_J^*\w\right\|_{L^2(f_I^{-1}(V),g_{I_-})} &\cr
&\le C\l_I^{-1/2} 
\|f_{I_+}^*\cdots f_J^*\w\|_{L^2_1(V,f_{I_+}^*\cdots f_J^*A_J,g_I)} &\cr
&\le C\|\w\|_{L^2_1(X_J',A_J,g_J)}. &\cr}
$$
by repeatedly applying Lemma \ref{\lemTaubesConfEst} in the last step. This gives (b); the proofs of (c) and (d) are similar.\qed

Finally, we obtain our estimate for the derivatives of $\hat\w$ with respect to the scales $\l_I$ and centres $x_I$.

\proclaim{Proposition}{\thmlabel{\propDhatOmDLaq}
There is a constant $C=(g_0,\sT)$ such that for any 
$\w\in\W^1(X,\ad P)$ and $t\in\sT$, the following bounds hold.
\item{\rm (a)}
$\|\rd{\hat\w}/\rd{\l_I}\|_{L^2(X_0,g_0)} 
\le C(\|\rd{\w}/\rd{\l_I}\|_{L^2(X,g)} 
+\l_I^{-1/2}\|\w\|_{\sL^2_1(X)})$;
\item{\rm (b)}
$\|\rd{\hat\w}/\rd{p_I}\|_{L^2(X_0,g_0)} 
\le C(\|\rd{\w}/\rd{p_I}\|_{L^2(X,g)} 
+\|\w\|_{\sL^2_1(X)})$.}

\pf From Lemma \ref{\lemDffstarOmDLaqEst} we have
$$\leqalignno
{\left\|\pd{\hat\w}{\l_I}\right\|_{L^2(X_0,g_0)} 
&\le C\sum_{J<I}\left\|\pd{\w}{\l_I}\right\|_{L^2(X_J',g_J)} 
+ C\l_I^{-1/2}\sum_{J\ge I}\|\w\|_{L^2_1(X_J',A_J,g_J)}, &\cr}
$$
and so (a) follows from Lemma \ref{\lemMetComparEst}. Similarly, Lemma \ref{\lemDffstarOmDLaqEst} gives
$$\leqalignno
{\left\|\pd{\hat\w}{p_I}\right\|_{L^2(X_0,g_0)} 
&\le C\sum_{J<I}\left\|\pd{\w}{p_I}\right\|_{L^2(X_J',g_J)} 
+ C\sum_{J\ge I}\|\w\|_{L^2_1(X_J',A_J,g_J)}. &\cr}
$$
and likewise, (b) follows from Lemma \ref{\lemMetComparEst}.\qed


\section{Derivatives with respect to scales and centres} 
\seclabel{\secDerScaleCentre}
We obtain $L^p$ estimates for the derivatives of the connections $A'$ and $\hat A'$ and of the $g$-self-dual curvature $F^{+,g}(A')$ with respect to the scales $\l_I$ and centres $x_I$. 

Throughout this section we require that $b_J=4N\l_J^{1/2}$ for all $J$. Let us first record 
the following bounds for the derivatives of the cutoff functions $\psi_J$ for $J=I_-$ or $I$:
$$\leqalignno
{\left|\rd{\psi_J}/\rd{\la_I}\right|_{g_J}&\le CN^{-1}\la_I^{-1}
\quad\hbox{and}\quad
\left|\rd{d\psi_J}/\rd{\la_I}\right|_{g_J}\le CN^{-2}\la_I^{-3/2}
\quad\hbox{on }X_J', &\numeqn\eqnlabel{\eqdPsiIEstdLaq}\cr
\left|\rd{\psi_J}/\rd{p_I}\right|_{g_J}&\le CN^{-1}\la_I^{-1/2}
\quad\hbox{and}\quad
\left|\rd{d\psi_J}/\rd{p_I}\right|_{g_J}\le CN^{-2}\la_I^{-1}
\quad\hbox{on }X_J', &\cr}
$$
where $\rd/\rd p_I\equiv p_I^\mu\rd/\rd q_I^\mu$ and $|p_I|\le 1$. The constant $C$ depends only on $g_J$.
We now begin with the $L^p$ estimates for derivatives of the connections $A'$.

\proclaim{Proposition}{\thmlabel{\propDAprimeDLaqEst}
Suppose $1\le p<\8$ and $I>0$. Then for sufficiently small $\la_0$, there is a constant $C=C(g_0,p,\sT)$ such that for any $t\in \sT$,
\item{\rm (a)} $\|\rd{A'}/\rd{\l_I}\|_{L^p(X,g)}\le C\l_I^{2/p-1/2}$,
\item{\rm (b)} $\|\rd{A'}/\rd{p_I}\|_{L^p(X,g)}\le C\l_I^{2/p}$.}

\pf (a) Observe that $\rd A'/\rd\l_I$ is non-zero only on the supports of
$\rd\psi_{I_-}/\rd\l_I$ and $\rd\psi_I/\rd\l_I$, given by the annuli 
$\W_I(\half b_I,b_I)$ in $X_{I_-}'$ and $\W_{Is}(\half b_I,b_I)$ in $X_I'$.

\step{1.} {\it Estimate of $\rd A'/\rd\l_I$ over $X_{I_-}'$.}
Recall that $\psi_{I_-}=1$ on the complement of the balls $B_I(b_I)$ in $X_{I_-}$, while $0<\psi_{I_-}<1$ on $\W_I(\half b_I,b_I)$, and 
$\psi_{I_-}=0$ on $B_I(\half b_I)$. We have 
$\s_I^*A'=\psi_{I_-}\s_I^*A_{I_-}$ on 
$\W_I(\half b_I,b_I)$ and thus
$\s_I^*\pd{A'}{\l_I}=
\pd{\psi_{I_-}}{\l_I}\s_I^*A_{I_-}$ on $X_{I_-}'$.
Since $|\rd\psi_I/\rd\l_I|\le C\l_I^{-1}$ by Eq. \eqref{\eqdPsiIEstdLaq} and 
$|\s_I^*A_{I_-}|_{g_{I_-}}\le C\l_I^{1/2}$ on 
$\W_I(\half b_I,b_I)$ by Lemmas \ref{\lemRadConnEst} and \ref{\lemRadConnSphEst}, we obtain the pointwise bound
$$
\left|\pd{A'}{\l_I}\right|_{g_{I_-}}
\le \cases{C\l_I^{-1/2}& on $\W_I(\half b_I,b_I)$ \cr
           0 & on $X_{I_-}'\less\W_I(\half b_I,b_I)$. \cr}
$$
Hence, we get the integral estimate
$$
\int_{X_{I_-}'}\left|\pd{A'}{\l_I}\right|_g^p\,dV_g
\le C\l_I^{2-p/2}, \eqlabel
$$
noting that $g=\tilde g_{I_-}$ on $X_{I_-}\less B_I(\half b_I)$ and appealing to Lemma \ref{\lemMetXIEst}.

\step{2.} {\it Estimate of $\rd A'/\rd\l_I$ over $X_I'$.}
A similar argument shows that
$$
\int_{X_I'}\left|\pd{A'}{\l_I}\right|_g^p\,dV_g
\le C\l_I^{2-p/2}, \eqlabel
$$
and combining the integral bounds from Steps 1 and 2 gives (a). 
For (b) we use the pointwise estimates $|\rd\psi_J/\rd p_I|\le C\l_I^{-1/2}$, $J=I_-,I$. The same argument as in (a) then gives the required bound. \qed
  
Our next task is to obtain a $L^p$ estimates for the derivatives of the $g$-self-dual curvature $F^{+,g}(A')$.

\proclaim{Proposition}{\thmlabel{\propDFAprimeDLaq}
Suppose $1\le p<4$ and $I>0$.
Then for sufficiently small $\la_0$, there exists a constant $C=C(g_0,p,\sT)$ such that for any $t\in \sT$,
\item{\rm (a)} $\|\rd{F^{+,g}(A')}/\rd{\l_I}\|_{L^p(X,g)}
\le C\l_I^{2/p-1}$, 
\item{\rm (b)} $\|\rd{F^{+,g}(A')}/\rd{p_I}\|_{L^p(X,g)}
\le C(\l_I^{2/p-1/2}+\ol\la^{1/2p})$.}

\pf (a) We note that $F^{+,g}(A')=F^{+,\tilde g_J}(\psi_JA_J)$ on $X_J'$ and so
$\rd F^{+,g}(A')/\rd\l_I$ is supported on 
$\cup_{J\ge I_-}X_J'$. It is convenient to obtain estimates separately over the regions $X_{I_-}'$, $X_I'$, and $X_J'$, $J>I$.

\step{1.} {\it Estimate of $\rd F^{+,g}(A')/\rd\l_I$ over $X_{I_-}'$.} 
On the annulus $\W_I(\half b_I,b_I)$ we have
$F^{+,g}(A')={1\over 2}(1+*_{\tilde g_{I_-}})F(\psi_{I_-}A_{I_-})$ and
$$\leqalignno
{F(\psi_{I_-}A_{I_-})
&=\psi_{I_-}F(A_{I_-})+d\psi_{I_-}\wedge\s_I^*A_{I_-}
+(\psi_{I_-}^2-\psi_{I_-})\s_I^*A_{I_-}\wedge\s_I^*A_{I_-}. &\cr}
$$
Therefore, we see that
$$\leqalignno
{\pd{F^{+,g}(A')}{\l_I}
&={1\over 2}(1+*_{\tilde g_{I_-}})\pd{F(\psi_{I_-}A_{I_-})}{\l_I}, &\cr
\pd{F(\psi_{I_-}A_{I_-})}{\l_I}
&=\pd{\psi_{I_-}}{\l_I}F(A_{I_-})
+\pd{d\psi_{I_-}}{\l_I}\wedge\s_I^*A_{I_-} 
+(2\psi_{I_-}-1)\pd{\psi_{I_-}}{\l_I}\s_I^*A_{I_-}\wedge\s_I^*A_{I_-}. &\cr}
$$
on $X_{I_-}'$.
The metric $\tilde g_{I_-}$ is independent of $\l_I$ 
and so applying the pointwise estimates of Lemmas \ref{\lemRadConnEst}, \ref{\lemRadConnSphEst}, and Eq. \eqref{\eqdPsiIEstdLaq}, we find that
$$
\left|\pd{F^{+,g}(A')}{\l_I}\right|_{g_{I_-}} 
\le \cases{C\l_I^{-1}& on $\W_I(\half b_I,b_I)$\cr
           0& on $X_{I_-}\less\W_I(\half b_I,b_I)$.\cr}
$$
Consequently, we obtain
$$\leqalignno
{\int_{X_{I_-}'}\left|\pd{F^{+,g}(A')}{\l_I}\right|_g^p\,dV_g
&\le C\l_I^{2-p}, &\numeqn\cr}
$$
where we observe that $g=\tilde g_{I_-}$ on $\W_I(\half b_I,b_I)$.

\step{2.} {\it Estimate of $\rd F^{+,g}(A')/\rd\l_I$ over $X_I'$.} We have
$F^{+,g}(A')={1\over 2}(1+*_{\tilde g_I})F(\psi_IA_I)$ and
$F(\psi_IA_I)=\psi_IF(A_I)
+d\psi_I\wedge\s_{Is}^*A_I
+(\psi_I^2-\psi_I)\s_{Is}^*A_I\wedge\s_{Is}^*A_I$ on $X_I'$. Thus,
$$\leqalignno
{\pd{F^{+,g}(A')}{\l_I}
&={1\over 2}\pd{*_{\tilde g_I}}{\l_I}F(\psi_IA_I)
+{1\over 2}(1+*_{\tilde g_I})\pd{F(\psi_IA_I)}{\l_I},
&\cr
\pd{F(\psi_IA_I)}{\l_I}
&=\pd{\psi_I}{\l_I}F(A_I)
+\pd{d\psi_I}{\l_I}\wedge\s_{Is}^*A_I 
+(2\psi_I-1)\pd{\psi_I}{\l_I}\s_{Is}^*A_I\wedge\s_{Is}^*A_I &\cr}
$$
on $X_I'$. Applying the pointwise estimates of Lemmas \ref{\lemRadConnSphEst}, \ref{\lemMetXIEst}, \ref{\lemDMetDLaqEst}, and Eq. \eqref{\eqdPsiIEstdLaq}, we find that
$$
\left|\phi_{In}^*\pd{F^{+,g}(A')}{\l_I}\right|_{g_I}(x)
\le \cases{0 & on $B_{Is}(\half b_I)$\cr
           C\l_I^{-1}& on $\W_{Is}(\half b_I,b_I)$\cr
           C|x|& on $X_I\less B_{Is}(b_I)$.\cr}
$$
Now $g=\tilde g_I$ on $X_I\less B_{Is}(\half b_I)$, and so applying the above estimates and H\"older's inequality gives
$$\leqalignno
{\int_{X_I'}\left|\pd{F^{+,g}(A')}{\l_I}\right|_g^p\,dV_g
&\le C\l_I^{2-p}, &\numeqn\cr}
$$
completing Step 2.

\step{3.} {\it Estimate of $\rd F^{+,g}(A')/\rd\l_I$ over $X_J'$, $J>I$.} 
We have
$$\leqalignno
{\pd{F^{+,g}(\psi_JA_J)}{\l_I}
&={1\over 2}\pd{*_{\tilde g_J}}{\l_I}F(\psi_JA_J) 
\quad\hbox{on }X_J',&\cr}
$$
since $F^{+,g}(A')={1\over 2}(1+*_{\tilde g_J})F(\psi_JA_J)$.
The pointwise estimates of Lemmas \ref{\lemAIprimeEst}, \ref{\lemMetXIEst}, and \ref{\lemDMetDLaqEst}
show that
$$
\left|\phi_{Jn}^*\pd{F^{+,g}(A')}{\l_I}\right|_{g_J}(x)
\le \cases{0 & on $B_{Js}(\half b_J)$\cr
           C|x|& on $X_J\less B_{Js}(\half b_J)$.\cr}
$$ 
Again, $g=\tilde g_J$ on $X_J\less B_{Js}(\half b_J)$, and so
$$\leqalignno
{\int_{X_J'}\left|\pd{F^{+,g}(A')}{\l_I}\right|_g^p\,dV_g
&\le C. &\numeqn\cr}
$$
Combining the integral estimates of Steps 1 to 3 then gives (a). 

\noindent (b) The argument is the same, except that we now use the cutoff function estimates $|\rd\psi_J/\rd p_I|\le C\l_I^{-1/2}$, $|\rd d\psi_J/\rd p_I|\le C\l_I^{-1}$, $J=I_-,I$, and metric estimates $|\rd\tilde g_J/\rd p_I|\le CN{\ol\la}^{1/2}$, $J\ge I$. \qed

Lastly, we have $L^2$ estimates of the derivatives of $\hat A'$ with respect to $\la_I$ and $x_I$. 
 
\proclaim{Proposition}{\thmlabel{\propDHatAprimeDLaqEst}
Suppose $I>0$. Then for sufficiently small $\la_0$, there is a constant $C=C(g_0,\sT)$ such that for any $t\in\sT$,
\item{\rm (a)} $\|\rd{\hat A'}/\rd{\l_I}\|_{L^2(X_0,g_0)}\le C$,
\item{\rm (b)} $\|\rd{\hat A'}/\rd{p_I}\|_{L^2(X_0,g_0)}\le C$.}

\pf (a) The connection one-forms over $X_0$ having non-zero derivatives with respect to $\l_I$ are given by
$$
\hat A'=\cases{f_0^*\cdots f_{I_-}^*\psi_{I_-}A_{I_-}
&over $f_0^{-1}\circ\cdots\circ f_{I_-}^{-1}(X_{I_-}')\subset X_0$\cr
f_0^*\cdots f_I^*\hat A_I'
&over $X_I\less B_{Is}(N_I^{-1}\l_I^{1/2})$,\cr}
$$
where $\hat A_I'$ is the $C^\8$ connection over $X_I$, $I>0$, given by
$$
\hat A_I'=\cases{f_{I_+}^*\cdots f_J^*\psi_JA_J& over the regions $f_{I_+}^{-1}\circ\cdots\circ f_J^{-1}(X_J')\subset X_I$ \cr
\psi_IA_I& over the complement of these regions in $X_I$. \cr}
$$
It is convenient to consider the estimates over these different regions of $X$ separately.

\step{1.} {\it Estimate of $\rd f_0^*\cdots f_{I_-}^*\psi_{I_-}A_{I_-}/\rd\l_I$.} We have 
$\hat A'=f_0^*\cdots f_{I_-}^*\psi_{I_-}A_{I_-}$, which is supported on 
$U_1\equiv f_0^{-1}\circ\cdots\circ f_{I_-}^{-1}(X_{I_-}')\subset X_0$, and so
$\pd{}{\l_I}\hat A'=f_0^*\cdots f_{I_-}^*\pd{}{\l_I}\psi_{I_-}A_{I_-}$
on $U_1$.
Lemma \ref{\lemXferEst} implies that
$$
\left\|f_0^*\cdots f_{I_-}^*\pd{\psi_{I_-}A_{I_-}}{\l_I}\right\|_{L^2(U_1,g_0)}
\le C\left\|\pd{\psi_{I_-}A_{I_-}}{\l_I}\right\|_{L^2(X_{I_-}',g_{I_-})}.
$$
We have $\s_I^*\psi_{I_-}A_{I_-}=\psi_{I_-}\s_I^*A_{I_-}$, where
the section $\s_I$ is chosen so that $\s_I^*A_{I_-}$ is in radial gauge, and so 
the pointwise estimates of Lemmas \ref{\lemRadConnEst}, \ref{\lemRadConnSphEst}, and Eq. \eqref{\eqdPsiIEstdLaq} show that
$$
\left|\pd{\psi_{I_-}A_{I_-}}{\l_I}\right|_{g_{I_-}}
\le \cases{C\l_I^{-1/2}& on $\W_I(\half b_I,b_I)$ \cr
           0 & on $X_{I_-}\less B_I(b_I)$. \cr}
$$
Noting that $g=\tilde g_{I_-}$ on $X_{I_-}\less B_I(b_I)$, we obtain the integral bound
$$\leqalignno
{\int_{X_{I_-}'}\left|\pd{\psi_{I_-}A_{I_-}}{\l_I}\right|_g^2\,dV_g
&\le C\l_I, &\cr}
$$
and combining the preceding integral estimates gives
$\|\pd{}{\l_I}\hat A'\|_{L^2(U_1,g_0)}
\le C\l_I^{1/2}$,
completing Step 1. 

\step{2.} {\it Estimate of $\rd f_0^*\cdots f_I^*\psi_IA_I/\rd\l_I$.} We denote $\hat A'=f_0^*\cdots f_I^*\psi_IA_I$, which is supported on
$U_2\equiv f_0^{-1}\circ\cdots\circ f_I^{-1}(X_I')\subset X_0$, and so
$\pd{}{\l_I}\hat A'=f_0^*\cdots f_{I_-}^*\pd{}{\l_I}f_I^*\psi_IA_I$
on $U_2$.
Repeated application of Lemma \ref{\lemXferEst} then gives the integral bound
$$
\left\|f_0^*\cdots f_{I_-}^*\pd{f_I^*\psi_IA_I}{\l_I}\right\|_{L^2(U_2,g_0)}
\le C\left\|\pd{f_I^*\psi_IA_I}{\l_I}\right\|_{L^2(X_{I_-}',g_{I_-})}.
$$  
Recall that Eq. \eqref{\eqDctauStarADLaq} implies
$\pd{}{\l_I}f_I^*\psi_IA_I
=-\l_I^{-1}f_I^*\iota_\rad F(\psi_IA_I)$ on $B_I'$.
The curvature $F(\psi_IA_I)$ is supported on 
$X_I\less B_{Is}(\half b_I)$ and 
$\pd{}{\l_I} f_I^*\psi_IA_I$
is supported on $B_I(\half N_I^{-1}\l_I^{1/2})$. Then, 
$$
|\phi_{In}^*\iota_\rad F(\psi_IA_I)|_{\de_I}(x)
\le K{|x|\over (1+|x|^2)^2},
$$ 
and since
$\phi_I^*f_I^*\iota_\rad F(\psi_IA_I)(x)
=\l_I^{-1}\phi_{In}^*\iota_\rad F(\psi_IA_I)(x/\l_I)$,
we obtain 
$$
\left|\phi_I^*\pd{f_I^*\psi_IA_I}{\l_I}\right|_{\d_{I_-}}(x)\le 
\cases{{\displaystyle 4K{\l_I^2|x|\over (\l_I^2+|x|^2)^2}}
          &if $|x|<\half N_I^{-1}\l_I^{1/2}$ \cr
        0 &if $|x|\ge \half N_I^{-1}\l_I^{1/2}$, \cr}
$$
where $K\equiv \|F(\psi_IA_I)\|_{L^\8(X_I,g_I)}$ is bounded by a constant $C$ independent of $\la_I$ by Lemma \ref{\lemAIprimeEst}.
But $g=\tilde g_{I_-}$ on
$B_I(\half N_I^{-1}\l_I^{1/2})\subset X_{I_-}''$, and moreover, the metrics $\tilde g_{I_-}$, $g_{I_-}$, and $\d_{I_-}$ are equivalent over the ball $B_I(\half N_I^{-1}\l_I^{1/2})$, with constants depending at most on $x_I$. Thus, we obtain the integral estimate
$$\leqalignno
{\int_{X_{I_-}'}\left|\pd{f_I^*A'}{\l_I}\right|_g^2\,dV_g&\le C\la_I^2, &\cr}
$$
and so, combining these bounds, we have
$\|\pd{}{\l_I}\hat A'\|_{L^2(U_2,g_0)}\le C\la_I$, completing Step 2.

\step{3.} {\it Estimate of $\rd f_0^*\cdots f_I^*\hat A_I'/\rd\l_I$.}  We have $\hat A_I'=f_{I_+}^*\cdots f_J^*\psi_JA_J$ over $V_3\equiv f_{I_+}^{-1}\circ\cdots\circ f_J^{-1}(X_J')\subset B_{I_+}'\subset X_I$, with $J>I$.
We denote $\hat A'=f_0^*\cdots f_I^*\hat A_I'$ and observe that
$\pd{}{\l_I}\hat A'
=f_0^*\cdots f_{I_-}^*\pd{}{\l_I}f_I^*\hat A_I'$ over
$U_3\equiv f_0^{-1}\circ\cdots\circ f_I^{-1}(V_3)\subset X_0$. Thus,
$$\leqalignno
{\pd{f_I^*\hat A_I'}{\l_I}
&=-\l_I^{-1}f_I^*\iota_\rad F(\hat A_I') 
=-\l_I^{-1}f_I^*\iota_\rad F(f_{I_+}^*\cdots f_J^*\psi_JA_J) &\cr
&=-\l_I^{-1}\iota_{f_{I*}^{-1}\rad} f_I^*\cdots f_J^*F(\psi_JA_J). 
&\cr}
$$
Note that $\rd\hat A_I'/\rd\l_I$ is supported on
$f_I^{-1}(B_{I_+}')\subset B_I'$.

As $\rad=y^\mu\rd/\rd y^\mu$ with respect to $y=\phi_{In}^{-1}$ on $X_I\less\{x_{Is}\}$, we have 
$f_{I*}^{-1}\rad=x^\mu\rd/\rd x^\mu$ with respect to $x=\phi_I^{-1}$ on $B_I'$.
If $|y|\le R_0$ on $B_{I_+}'$, for some constant $0<R_0<\8$ depending at most on $x_I$, then $|x|\le R_0\l_I$ on $f_I^{-1}(B_{I_+}')$. Thus,
$|f_{I*}^{-1}\rad|_{g_{I_-}}\le R_0\l_I$ on $f_I^{-1}(B_{I_+}')$
and so we have the pointwise bound
$$\leqalignno
{\left|\pd{\hat A_I'}{\l_I}\right|_{g_{I_-}}
&\le R_0|f_I^*\cdots f_J^*F(A_J')|_{g_{I_-}}
\qquad\hbox{on }f_I^{-1}(B_{I_+}'). &\cr}
$$
Therefore, with the aid of repeated applications of Lemma \ref{\lemXferEst}, we find that
$$\leqalignno
{\left\|f_0^*\cdots f_{I_-}^*\pd{\hat A_I'}{\l_I}\right\|_{L^2(U_3,g_0)}
&\le C\left\|\pd{\hat A_I'}{\l_I}\right\|_{L^2(f_I^{-1}(V_3),g_{I_-})} &\cr
&\le C\|f_I^*\cdots f_J^*F(A_J')\|_{L^2(f_I^{-1}(V_3),g_{I_-})} &\cr
&\le C\|F(A_J')\|_{L^2(X_J',g_J)}. &\cr}
$$
and since $\|F(A_J')\|_{L^2(X_J',g_J)}\le C$, this gives
$\left\|\pd{}{\l_I}\hat A'\right\|_{L^2(U_3,g_0)}\le C$, completing Step
3. Combining the results from Steps 1 to 3 then yields (a).  
For (b) we use the cutoff function estimate $|\rd\psi_J/\rd p_I|\le C\l_I^{-1/2}$, $J=I_-,I$. The vector field $\rad$ is replaced by $\trans=p_I^\mu\rd/\rd y^\mu$, with respect to the coordinates $y=\phi_{In}^{-1}$. Then, $f_{I*}^{-1}\trans= \la_Ip_I^\mu\rd/\rd x^\mu$ with respect to the coordinates $x=\phi_I^{-1}$ and we have the vector field estimate $|f_{I*}^{-1}\trans|\le R_0\la_I$ on $f_I^{-1}(B_{I_+}')$. The required bound then follows by an argument similar to that of (a). \qed 


\section{Derivatives with respect to bundle gluing parameters}
\seclabel{\secDAprimeDGlu}
The purpose of this section is to obtain estimates for the derivatives of the almost ASD connections $A'$ and $\hat A'$ with respect to the bundle gluing parameters $\rho_I\in\Gl_I$, $I>0$. These estimates may be extracted from [D-K, \S 7.2] and we include them here for completeness. 

Since we wish to differentiate a family of connections $A'(\r_I)$ on a family of $G$-bundles $P(\r_I)$ with respect to the gluing parameters $\r_I\in\Gl_I$, we first pull this family back to an equivalent family on a {\it fixed} bundle, say $P(\bar\rho_I)$, as described in [D-K, p. 296]. Let $\bar\rho_I\in\Gl$ be a given gluing parameter: then points $\r_I$ in a small neighbourhood of $\bar\rho_I$ in $\Gl_I$ can be written in the form $\r=\bar\rho_I\exp(v)$, where $v\in V_I\equiv\ad P_I|_{x_{Is}}\simeq \frak{g}$. One regards the fibres of $P_{I_-}$ and $P_I$ as being identified by $\bar\rho_I$ and so $v$ may considered as a local section of both $P_{I_-}$ and $P_I$, covariantly constant with respect to the connections $A_{I_-}'$, $A_I'$. 

We digress in order to construct a set of cutoff functions $\{\ga_I\}$ on $X$ such that
$\sum_{I\in\sI}\ga_I=1$. These cutoffs will be needed here and again in \S 5.1 for patching together certain integral operators over the $X_I$ to give an integral operator over $X$.
Choose a bump function $\ga\in C^\8({\Bbb R}^1)$ such that 
$\ga(t)=1$ if $t\ge 2$ and $\ga(t)=0$ if $t\le \half$.  
Define a cut-off function $\ga_\la\in C^\8(\RR^4)$ by 
$$
\c_\la(x)\equiv \ga(|x|/\l^{1/2}),\qquad x\in\RR^4. \eqlabel\eqnlabel{\eqGamma}
$$
Now define $C^\8$ cutoff functions $\ga_I$ on each summand $X_I$ by setting
$$
\ga_I\equiv(\phi_{Is}^{-1})^*(1-\c_{\la_I})
\prod_{I_+}(\phi_{I_+}^{-1})^*\c_{\la_{I_+}}
\quad\hbox{on }X_I, \eqlabel\eqnlabel{\eqGammaI}
$$ 
where the factor $(\phi_{Is}^{-1})^*(1-\c_{\la_I})$ is omitted when $I=0$. Note that $\ga_I=0$ on the balls $B_{Is}(\half\la_I^{1/2})$ and $B_{I+}(\half\la_{I_+}^{1/2})$ in $X_I$. We extend $\ga_I$ to a $C^\8$ cutoff function on $X_I$ by zero on these balls and by $1$ on the complement of the larger balls $B_{Is}(2\la_I^{1/2})$ and $B_{I+}(2\la_{I_+}^{1/2})$ in $X_I$;
then extend by zero outside $X_I''\subset X$ to give $\ga_I\in C^\8(X)$. By construction, we have $\sum_{I\in\sI}\ga_I=1$ on $X$, with a slight abuse of notation. Indeed, note that $f_I$ maps the annulus $\W_I(\half\l_I,2\l_I)$ around the point $x_I$ in $X_{I_-}$ onto the annulus $\W_{Is}(\half\l_I,2\l_I)$ around the south pole $x_{Is}$ in $X_I$. Then, $f_I^*\c_I+\c_{I_-}=1$
on each annulus $\Om_I$.
Lastly, note that there is a constant $C$, depending at most on the metric $g_0$, such that
$$
|d\ga_I|_{g_I}\le C\l_I^{-1/2}\quad\hbox{on }\Om_I, \Om_{Is}\quad\hbox{and}\quad
\|d\ga_I\|_{L^p(X_I,g_I)}\le C\ol\l^{2/p-1/2}.
\eqlabel\eqnlabel{\eqGammaIEst}
$$
for any $1\le p<\8$. 
Define gauge transformations $u_{I_-}(v)$ on $\Aut P_{I_-}|_{X_I'}$ and
$u_I(v)$ on $\Aut P_I|_{X_I'}$ by setting
$$\leqalignno
{u_{I_-}(v)&=\cases{\exp(\ga_Iv) &on $\Om_I$ \cr
           1           &on $X_{I_-}'\less \Om_I$ \cr}, &\numeqn\cr
u_I(v)&=\cases{\exp(-\ga_{I_-}v) &on $\Om_{Is}$ \cr
           1           &on $X_I'\less \Om_{Is}$. \cr} &\cr}
$$
Note that $u_I$ has a natural extension to a gauge transformation of $P_I$ over all of $X_I$ --- equal to $\exp(-v)$ on $B_{Is}(N_I^{-1}\la_I^{1/2})$, the ball enclosed by the annulus $\Om_{Is}$. Similarly for the gauge transformation $u_{I_-}$. After identifying the bundles and base manifolds over $\Om=\Om_I=\Om_{Is}$, we have 
$u_{I_-}u_I^{-1}=\exp((\ga_{I_-}+\ga_I)v)=\exp(v)$.
Hence, relative to the flat connections $A_{I_-}',A_I'$, the gauge transformations $u_I$ differ by a constant bundle automorphism over $\Om$ and so their action on the connection $A'(\bar\rho_I)$ is the same:
$u_{I_-}(A'(\bar\rho_I))|_\Om=u_I(A'(\bar\rho_I))|_\Om$.
Therefore, while the automorphisms $u_I$ do {\it not} patch together to give a global automorphism of $P(\bar\rho_I)$, their actions on the connection $A'(\bar\rho_I)$ do. Indeed, we can define a connection $A'(\bar\rho_I,v)$ on $P(\bar\rho_I)$ by
$$
A'(\bar\rho_I,v)=\cases{u_{I_-}(A'(\bar\rho_I)) &on $X_{I_-}'$ \cr
                  u_I(A'(\bar\rho_I)) &on $X_I'$. \cr} \eqlabel
$$
If $\r_I=\bar\rho_I\exp(v)$, the connections $A'(\bar\rho_I,v)$ and $A'(\r_I)$ are gauge equivalent [D-K, p. 296]. Thus, as desired, we have an equivalent family of connections $A'(\bar\rho_I,v)$ on
the fixed connected sum bundle $P=P(\bar\rho_I)$. Let $L_I\subset\Gl_I$ be a coordinate neighbourhood and suppose $\bar\rho_I\in L_I$. Then 
$$
\frak{g}\supset B_\frak{g}\too L_I\subset\Gl_I, \qquad v\longmapsto \rho_I(v)\equiv \bar\rho_I\exp(v) \eqlabel\eqnlabel{\eqGluChart}
$$
is a coordinate chart centred at $\bar\rho_I$, where $B_\frak{g}$ is the unit ball in $\frak{g}$, and there is a $C^\8$ embedding
$$
\frak{g}\supset B_\frak{g}\too \cA_{X,P}^*, 
\qquad v\longmapsto A'(\bar\rho_I,v).\eqlabel\eqnlabel{\eqAprimeGluChart}
$$
It remains to consider the derivative of the family $A'(\bar\rho_I,v)$ with respect to $v$. 

Recall that if $u=u(s)$ is a one-parameter family of gauge transformations, $B$ is a fixed 
connection, and $B^u(s)$ is the induced one-parameter family of gauge transformed connections, then 
$dB^u/ds(0)= d_{B^u}(u^{-1}\dot u(0))$,
where $u^{-1}\dot u(0)\in\W^0(X,\ad P)$. Although the $u_{I_-}, u_I$ are not globally defined gauge transformations, this differentiation formula still applies to the one-parameter families $u_{I_-}(s)=u_{I_-}(sv)$ and
$u_I(s)=u_I(sv)$.
Therefore, we have
$$
\pd{A'}{v}(\bar\rho_I)\equiv\deriv{}{s}A'(\bar\rho_I,sv)\Bigr|_{s=0}
=\cases{d_{A'}(\ga_Iv) &on $X_{I_-}'\cap\W$ \cr
                              -d_{A'}(\ga_{I_-}v) &on $X_I'\cap\W$ \cr
                               0 &on $X\less\W$. \cr} 
\eqlabel\eqnlabel{\eqDAprimeDv}
$$
This leads to the following estimate for the derivative of the family
$A'(\rho_I)$ with respect to the gluing parameters $\rho_I$; 
a related and more general estimate is given by
Lemma 7.2.49 in [D-K].

\proclaim{Proposition}{\thmlabel{\propDAprimeDvEst} Let $2\le p<4$ and suppose that $4\le q<\8$ is determined by $1/4+1/q=1/p$. 
There is a constant $c=c(g_0,p,\script{T})$ such that
\item{\rm (a)}
$c|v|\la_I^{2/p-1}\le \|\rd{A'}/\rd{v}\|_{L^q(X,g)}
\le c^{-1}|v|\la_I^{2/p-1}$,
\item{\rm (b)}
$c|v|\la_I^{2/p-1/2}\le \|\rd{A'}/\rd{v}\|_{L^p(X,g)}
\le c^{-1}|v|\la_I^{2/p-1/2}$.}

\pf Note that $\ga_{I_-}+\ga_I=1$ on $\Om$ and so $d_{A'}(\ga_Iv)=-d_{A'}(\ga_{I_-}v)$ on $\Om$. Moreover,
$d_{A'}(\ga_Iv)=d\ga_I\wedge v$ on $\Om$,
and so we have
$\|d_{A'}(\ga_Iv)\|_{L^q(X,g)}=|v|\cdot\|d\ga_I\|_{L^q(X,g)}$.
From Eq.
\eqref{\eqGammaIEst} there is a constant $c>0$ independent of $\la_I$ such that
$$
c|v|\la_I^{2/q-1/2}\le \left\|\pd{A'}{v}\right\|_{L^q(X,g)}
\le c^{-1}|v|\la_I^{2/q-1/2},
$$
since $\|\rd A'/\rd v\|_{L^q(X,g)}=\|d_{A'}(\ga_Iv)\|_{L^q(X,g)}$. Then (a) follows since $2/q-1/2=2/p-1$, and likewise for (b). \qed

Using the conformal maps $f_J$, we pull back the family $A'=A'(\bar\rho_I,v)$ on the fixed bundle $P$ over $X$ to a family $\hat A'(\bar\rho_I,v)$ on the fixed bundle $\hat P$ over $X_0$.

\proclaim{Proposition}{\thmlabel{\propDhatAprimeDvEst} If $2\le p<4$, 
there is a constant $C=C(g_0,p,\script{T})$ such that for any $t\in\script{T}$,
$\|\rd{\hat A'}/\rd{v}\|_{L^p(X_0,g_0)}\le C\la_I^{2/p-1/2}$.}

\pf Since $\rd A'/\rd v=0$ outside the annulus $\Om_{Is}\subset X_I'$, Proposition \ref{\propDAprimeDvEst} gives
$$\leqalignno
{\left\|\pd{A'}{v}\right\|_{L^p(X_I',g_I)} &\le C\la_I^{2/p-1/2}. &\cr}
$$
But $\hat A'=f_0^*\cdots f_I^*A'$ on 
$U\equiv f_0^{-1}\circ\cdots\circ f_I^{-1}(X_I')\subset X_0$, and so
Lemma \ref{\lemXferEst} gives
$$
\left\|f_0\cdots f_I^*\pd{A'}{v}\right\|_{L^p(U,g_0)} 
\le C\left\|\pd{A'}{v}\right\|_{L^p(X_I',g_I)}.
$$
Combining these estimates gives the desired bound. \qed


\section{Derivatives with respect to lower moduli}
In this section we obtain $L^p$ estimates for the derivatives of the
connections $A'$, 
$\hat A'$, and the self-dual curvature $F^{+,g}(A')$ with respect to the
``lower moduli parameters'' $t_I$.
The bundle $P_I$ carrying the family of connections $\{A_I(t_I)\}_{t_I\in
T_{A_I}}$ can be assumed to be fixed with respect to the
parameters $t\in T_{A_I}$ since the space $T_{A_I}$ --- an open ball in
$H_{A_I}^1$ centred at $0$ --- is
contractible. However, the local sections $\si_{I+}(t_I)$ are defined by the
connections $A_I(t_I)$ (together a choice of point in $P_I|_{x_{I_+}}$) and
will vary with $t_I$. Thus, the bundle gluing maps for the connected sum
bundle $P$, defined by $\si_{I_+}(t_I)\mapsto\si_{I_+s}\tilde\rho_{I_+}(t_I)$
(suppressing the identification map $f_I:\Om_{I_+}\to \Om_{I_+s}$), will in
general vary with $t_I$. We may suppose that the remaining parameters are fixed
and thus we obtain a family of connections $A'(t_I)$ on a family of bundles
$P(t_I)$. The difficulty, of course, is that unless we have a family of
connections defined on a fixed bundle, we cannot define the derivative $\rd
A'/\rd t_I$. Problems such as these are discussed in [D87, p.
423]. For our purposes, we note the bundles are all isomorphic and as
$T_{A_I}$ is 
contractible, the connections $A'(t_I)$ could be pulled back by bundle
isomorphisms $h_I\in\Hom(P(0),P(t_I))$ to an equivalent family 
$h_I^*A'(t_I)$ on the {\it fixed} bundle $P(0)$ and then we could define
$$
\pd{A'}{t_I}\equiv \pd{h_I^*A'}{t_I}.\eqlabel\eqnlabel{\eqdA'dtIdefn}
$$
Since any two such families $h_I(t_I)$ of bundle isomorphisms would differ
by a family of automorphisms of the fixed bundle $P(0)$, then
$\rd A'/\rd t_I$ would give a well-defined
tangent vector to $\cB_{P(0)}^*$ at $[A'(t_I)]$ using \eqref{\eqdA'dtIdefn}.
Naturally,
the analogous remarks apply to the family of connections $\hat A'(t_I)$ on the
bundles $\hat P(t_I)$.

In our case, a family of isomorphisms $h_J(t_I):P_J(0)\to P_J(t_I)$
may be described quite explicitly, in a manner similar to that of
\S\ref{\secDAprimeDGlu}, and these will give a gauge equivalent family of
connections 
$h_I^*A'(t_I)$, ${\hat h}_I^*\hat A'(t_I)$ on fixed bundles $P(0)$, $\hat P(0)$
respectively, although just as in \S\ref{\secDAprimeDGlu}, the isomorphisms
$h_J(t_I)$ 
will {\it not} patch together to give a global isomorphism of $P(0)$ with
$P(t_I)$ or $\hat P(0)$ with $\hat P(t_I)$. 
Nonetheless Eq. \eqref{\eqdA'dtIdefn} still makes sense and this
allows us to estimate the length of the tangent
vector $\rd A'/\rd t_I$ in terms of derivatives of the local connection
one-forms, as desired. Let $h_I(t_I):P_I(0)\to P_I(t_I)$ be
a family of bundle isomorphisms represented locally by 
$\si_{I_+}(0)\mapsto \si_{I_+}(t)\theta_{I_+}(t_I)$. Then
$h_I^*A_I(t_I)$ is an equivalent family on the fixed bundle $P_I(0)$, with
$$
\si_{I_+}(0)^*h_I(t_I)^*A_I(t_I)
=\theta_{I_+}(t_I)^{-1}\si_{I_+}(t_I)^*A_I(t_I)\theta_{I_+}(t_I) +
\theta_{I_+}(t_I)^{-1}d\theta_{I_+}(t_I). 
$$
Note that while the local connection one-forms
$\si_{I_+}(t_I)^*A_I(t_I)$ 
are in radial gauge, this will {\it not} in general be the case for the
one-forms $\si_{I_+}(0)^*h_I(t_I)^*A_I(t_I)$. 
We next consider the variation in the bundle gluing maps
$\tilde\rho_{I_+}(t_I)$ induced by the variation in $\si_{I_+}(t_I)$ with
$t_I$. Over $X_I$, we replace $\theta_{I_+}(t_I)$ above
by $\theta_{I_+}(t_I)\exp(\ga_{I_+}v_I(t_I))$ and over $X_{I_{+s}}$ 
define $h_I(t_I)$ by right multiplication with $\exp(\ga_Iv_I(t_I))$.
Recalling the notation of \S\ref{\secDAprimeDGlu},
$v_I:T_{A_I}\to\frak{g}$ is a smooth map with $v(0)=0$ defined (for small
enough 
$T_{A_I}$) by the identity $\rho_{I_+}(t_I)=\rho_{I_+}(0)\exp(v_I(t_I))$.
Lastly, for $J\ne I, I_+$, we set $h_J(t_I)=1$.
{\it Then, for the remainder of this article, we require that the derivatives
$\rd A'/\rd t_I$ be defined by \eqref{\eqdA'dtIdefn}}. 

This understood, we obtain the following estimates
for the derivatives with respect to the parameters $t_I$ of the connections
$A'$ and $\hat A'$ and for the $g$-self-dual curvature $F^{+,g}(A')$. The
proofs are straightforward, following the pattern in
\S\ref{\secDerScaleCentre}, 
and so are omitted.

\proclaim{Proposition}{\thmlabel{\propDAprimeDt}
Let $1\le p<\8$. For sufficiently small $b_0$, there exists a
constant $C=C(g_0,p,\script{T})$ such that for any $t\in \script{T}$,  
\item{\rm (a)} $\|\rd{A'}/\rd{t_I}-\rd{A_I}/\rd{t_I}\|_{L^p(X_I'',g_I)}
\le C\la_I^{2/p}$
\item{\rm (b)} $\|\rd{A'}/\rd{t_I}\|_{L^p(X,g)}\le C$.}

\proclaim{Proposition}{\thmlabel{\propDFplusDt}
Let $1\le p<\8$. For sufficiently small $b_0$, there is a constant
$C=C(g_0,p,\script{T})$ such that for any $t\in \script{T}$, 
$\|\rd{F^{+,g}(A')}/\rd{t_I}\|_{L^p(X,g)}\le C\bar\l^{2/p-1/2}$.}

\proclaim{Proposition}{\thmlabel{\propDhatAprimeDt}
For sufficiently small $b_0$, there is a constant $C=C(g_0,N,\script{T})$
such that for any $t\in \script{T}$,
$\|\rd{\hat A'}/\rd{t_I}\|_{L^2(X_0,g_0)}\le C$.}

\pf Let $U\equiv f_0^{-1}\cdots f_I^{-1}(X_I')\subset X_0$ and note
note that $\pd{}{t_I}\hat A'=\pd{}{t_I}f_0^*\cdots f_I^*\psi_IA_I$,
which is $\pd{}{t_I}f_0^*\cdots f_I^*\psi_IA_I
=f_0^*\cdots f_I^*\pd{}{t_I}\psi_IA_I$ on $U$ and zero elsewhere. Now
$$
\left\|f_0^*\cdots f_I^*\pd{\psi_IA_I}{t_I}\right\|_{L^2(U,g_0)}
\le C\left\|\pd{\psi_IA_I}{t_I}\right\|_{L^2(X_I',g_I)}
$$
by Lemma \ref{\lemXferEst} and
so the result follows.
\qed


\section{Differentials of the approximate gluing maps}
We close this Chapter by
summarising the results of the preceding sections and record our bounds
for the differentials of the approximate gluing maps $\script{J}'$ (which
follow by combining Propositions \ref{\propDAprimeDLaqEst},
\ref{\propDAprimeDt}, and 
\ref{\propDAprimeDvEst})) and
$\hat\script{J}'$ (which follow by
combining Propositions \ref{\propDHatAprimeDLaqEst},
\ref{\propDhatAprimeDt}, and \ref{\propDhatAprimeDvEst}).  

\proclaim{Theorem}{\thmlabel{\thmDJprimeEst} Let $\script{J}':\script{T}\to \cB_{X,k}^*$ be the approximate gluing map $t\mapsto [A'(t)]$. 
Assume $b_I=4N_I\l_I^{1/2}$ for all $I$. Then for sufficiently small $\l_0$ and any $t\in \script{T}$, there is a constant $C=C(g_0,\script{T})$ such that the following estimates hold.
\item{\rm (a)} $\|D\script{J}'(\rd/\rd{t_I^\al})\|_{L^2(X,g)} \le C$,
\item{\rm (b)} $\|D\script{J}'(\rd/\rd{\rho_I^\be})\|_{L^2(X,g)} 
\le C\la_I^{1/2}$,
\item{\rm (c)} $\|D\script{J}'(\rd/\rd{x_I^\mu})\|_{L^2(X,g)}\le C$,
\item{\rm (d)} $\|D\script{J}'(\rd/\rd{\la_I})\|_{L^2(X,g)}\le C$.}

\proclaim{Theorem}{\thmlabel{\thmDhatJprimeEst} Let $\hat\script{J}':\script{T}\to \cB_{X_0,k}^*$ be the approximate gluing map $t\mapsto [\hat A'(t)]$. 
Let $b_I=4N_I\l_I^{1/2}$ for all $I$. Then for sufficiently small $\l_0$ and any $t\in \script{T}$, there is a constant $C=C(g_0,\script{T})$ such that the following estimates hold.
\item{\rm (a)} $\|D\hat\script{J}'(\rd/\rd{t_I^\al})\|_{L^2(X_0,g_0)} \le C$,
\item{\rm (b)} $\|D\hat\script{J}'(\rd/\rd{\rho_I^\be})\|_{L^2(X_0,g_0)} 
\le C\la_I^{1/2}$,
\item{\rm (c)} $\|D\hat\script{J}'(\rd/\rd{x_I^\mu})\|_{L^2(X_0,g_0)}\le C$,
\item{\rm (d)} $\|D\hat\script{J}'(\rd/\rd{\la_I})\|_{L^2(X_0,g_0)}\le C$.}

\endchapter
\eqnno=0
\secno=0
\thmno=0


\chapter{Bubble Tree Compactification of the Anti-self-dual Moduli Space}
In order to describe the ends of the moduli space $M_{X_0,k}(g_0)$ one
customarily appeals to the Uhlenbeck compactification
$\ol{M}_{X_0,k}^u(g_0)$. This allows one to give quite
explicit descriptions of the parts of the ends away from the diagonals in
the symmetric products $M_{X_0,k}(g_0)\times s^l(X_0)$ appearing in the
compactification, as for example in [D86, \S V] and [D-K, \S 8.2]. These
examples consider ideal boundary
points of the form $(A_0,x_1,\dots,x_l)$, where the $x_i$ are distinct
points with 
multiplicity $1$ and $A_0$ is a $g_0$-ASD connection
over $X_0$. Open neighbourhoods of
$(A_0,x_1,\dots,x_l)$ in $\ol{M}_{X_0,k}^u(g_0)$ are then constructed 
by gluing standard one-instantons onto $A_0$.

In order to construct open neighbourhoods of ideal boundary points
corresponding to the diagonals of $\ol{M}_{X_0,k}^u(g_0)$ we must employ
the iterated gluing construction of Chapters 3 and 5. This strategy is
mentioned briefly in [D-K, \S8.2].
The construction gives a homeomorphism $\hat\script{J}:\script{T}^0/\Ga\to\ol\sV$, where $\ol\sV$ is an open neighbourhood of a boundary point in $\ol{M}_{X_0,k}^u(g_0)$ --- a `gluing neighbourhood'. In order to use this procedure to describe the ends of $\ol{M}_{X_0,k}^u(g_0)$, we need to show that $\ol{M}_{X_0,k}^u(g_0)$ is covered by finitely many such gluing neighbourhoods. In particular, we need to show that any point in $M_{X_0,k}(g_0)$ which is sufficiently close to the ideal boundary (with respect to the Uhlenbeck topology) lies in the image of a gluing map $\hat\script{J}$. This is accomplished in two steps:
\smallskip

\noindent{\bf Step 1.} We show that any sequence $\{A_\al\}$ of $g_0$-ASD
connections over $X_0$ converging {\it weakly} to a limit
$(A_0,x_1,\dots,x_{m_0})$ determines a sequence of metrics $\{g_\al\}$ and
a sequence $\{\check A_\al\}$ of $g_\al$-ASD connections over a connected
sum $X\equiv\#_{I\in\script{I}}X_{I\al}'$ which converges {\it strongly} to a limit $(A_I)_{I\in\script{I}}$, in the sense of [D-K, \S7.3]. Here, $(X,g_\al)$ is conformally equivalent to $(X_0,g_0)$, for all $\al$, and is defined exactly as in \S\ref{\secGluCon} and \S\ref{\secMetBubbleTree}.
\smallskip

\noindent{\bf Step 2.} We apply an analogue of Theorem 7.3.2 [D-K] to show that the new sequence $\{\check A_\al\}$ is $D_q$-convergent, $q\ge 4$, in the sense of
[D-K, \S7.3]. The appropriate analogue of Theorem 7.2.62 [D-K] then shows that 
the points $[A_\al]\in M_{X,k}(g_\al)$ lie in the image of some
$\script{J}$ for sufficiently large $\al$. Consequently, the points
$[A_\al]\in M_{X_0,k}(g_0)$ lie in the image of the corresponding map
$\hat\script{J}$, for some parameter space $\sT^0/\Ga$. The choice of
parameter space $\script{T}^0/\Ga$ is essentially determined by
$(A_I)_{I\in\script{I}}$, which we
call the {\it strong}\/ or {\it bubble tree limit}\/ of the sequence
$\{A_\al\}$.  
\smallskip

In this Chapter we discuss Step 1 and describe the bubble
tree compactification of the moduli space of anti-self-dual $\SU(2)$ 
connections --- the
extension to the general case of compact, semi-simple Lie groups being straightforward.
Step 2 is discussed in \S\S 5.1 and 5.2 after the necessary analytical
framework has been established. Throughout this Chapter, we suppose only that
$X_0$ is a closed, oriented, simply-connected $C^\8$ four-manifold, $g_0$
is a $C^\8$ 
metric, and $G=\SU(2)$. 

\section{Uhlenbeck compactification} 
We recall the definition of the Uhlenbeck
compactification [D-K] and describe some of the related 
convergence results we will
need for our description of the bubble tree compactification.

\defn\thmlabel{\defIdealConn} An {\it Uhlenbeck ideal 
$g_0$-ASD connection} on a $G$ bundle $P$ over $X_0$ with 
$c_2(P)=k\ge 0$ is a pair $(A_0,Z_0)$,
where $A_0$ is a $g_0$-ASD connection on a $G$ bundle $P_0$ over $X_0$ with 
$c_2(P_0)=k_0\ge 0$ and $Z_0=\{x_i\}_{i=1}^{m_0}$
is a (possibly empty) set of points in $X_0$ with multiplicities $k_i\ge
1$, for
$i=1,\dots,m_0$, such that $\sum_{i=0}^{m_0} k_i=k$. The curvature
density of $(A_0,Z_0)$ is defined to be the Borel measure
$$
\mu(A_0,Z_0)=|F(A_0)|_{g_0}^2+8\pi^2\delta_{Z_0}, \eqlabel
$$
where $\delta_{Z_0}\equiv \sum^{m_0}_{i=1}k_i\delta_{x_i}$,
so that the total mass of $\mu(A_0,Z_0)$ is $8\pi^2k$.
Setting $l=k_1+\cdots+k_m$ and repeating points according to their
multiplicity, one obtains an element $(x_1,\dots,x_l)$ of the
symmetric product $s^l(X_0)$. 

\defn\thmlabel{\defWeakConv} Let $\{A_\al\}_{\a=1}^\8$, be a sequence of
$g_0$-ASD connections on a $G$ bundle $P$ over $X_0$ with $c_2(P)=k\ge 0$ 
and let 
$(A_0,Z_0)$ be an ideal $g_0$-ASD connection on $P$.
Then the sequence $\{A_\al\}$ {\it converges weakly} to $(A_0,Z_0)$ if:
\item{(a)} The sequence $\{\mu_\a\}_{\a=1}^\8$ converges to $\mu(A_0,Z_0)$
in the weak-* topology on measures;
\item{(b)} There is a sequence of $C^\8$ bundle maps 
$\ga_\al:P_0|_{X_0\less Z_0}\to P|_{X_0\less Z_0}$
such that $\ga_\al^*A_\al$ converges in $C^\8$ on compact subsets of
$X_0\less Z_0$ to the connection $A_0$. 
Equivalently, require that for any integer $n\ge 1$, there is a sequence of
$L_{n+1}^2$ bundle maps $\ga_\al$ such that $\ga_\al^*A_\al$ converges in
$L_{n,\loc}^2$ on $X_0\less Z_0$ to $A_0$. 
\medskip

Via the natural extension of Definition \ref{\defWeakConv} to
sequences of ideal connections, the set of all
Uhlenbeck ideal $g_0$-ASD connections of fixed second Chern class $k$, 
$IM_{X_0,k}(g_0)\equiv\coprod_{l=0}^k (M_{X_0,k-l}(g_0)\times s^l(X_0))$,
is endowed with a metrisable topology. 
Let $\ol{M}_{X_0,k}^u(g_0)$ be the closure of $M_{X_0,k}(g_0)$ in
$IM_{X_0,k}(g_0)$. According to [D-K, Theorem 4.4.4], any
infinite sequence in $M_{X_0,k}(g_0)$ has a weakly convergent subsequence 
with limit point in $\ol{M}_{X_0,k}^u(g_0)$, and in particular, the latter
space is compact [D-K, Theorem 4.4.3].

For our description of the bubble tree compactification, we will need the 
following minor extension of the convergence result in Theorem 4.4.4 [D-K] 
and its cousin, Proposition 9.4.2 [D-K], which allows for a sequence of
metrics $\{g_\al\}$ converging to $g_0$ in $C^\8$. The proof employs standard
arguments well described in [D-K, \S 4.4] and is left to the reader.

\proclaim{Proposition}{\thmlabel{\propMetExhstCpt} 
Let $\{U_\al\}_{\al=1}^\8$ be an exhaustion of the punctured manifold
$X_0\less\{p\}$ by an increasing sequence $\{U_\al\}_{\al=1}^\8$ of 
precompact open sets, so that
$U_1\Subset U_2\Subset\cdots\subset X_0\less\{p\}$ and
$\cup_{\al=1}^\8 U_\al=X_0\less\{p\}$.
Let $\{g_\al\}_{\al=1}^\8$ be a sequence of metrics on the subsets $U_\al$ converging in $C^r$ ($r\ge 3$) on compact subsets of $X_0\less\{p\}$ to a $C^r$ metric $g_0$ on $X_0$. 
Let $P$ be a $G$ bundle over $X_0\less\{p\}$ and let 
$\{A_\al\}_{\al=1}^\8$ be a sequence of $g_\al$-ASD connections on the restrictions $P|_{U_\al}$. If  there is a constant $M<\8$ such that 
$$
\int_{U_\al}|F(A_\al)|_{g_\al}^2\,dV_{g_\al}\le M\quad\hbox{for all $\al$},
$$
then there is a set of points $Z_0=\{x_i\}_{i=1}^{m_0}\subset X_0$ and a 
$g_0$-ASD
connection $A_0$ on a $G$ bundle $P_0$ over $X_0$ such that a subsequence
$\{A_\al\}_{\al=1}^\8$ converges weakly to $(A_0,Z_0)$. \qed} 

The mass of the Uhlenbeck limit $(A_0,Z_0)$ in Proposition
\ref{\propMetExhstCpt} is $8\pi^2$ times an integer and may be computed
from the weakly convergent sequence $\{A_\al\}_{\al=1}^\8$ by 
$$
\lim_{n\to\8}\lim_{\al\to\8}
\int_{V_n}|F(A_\al)|_{g_\al}^2\,dV_{g_\al}, \eqlabel
$$
where $\{V_n\}_{n=1}^\8$ is any exhaustion of $X_0\less\{p\}$ by an
increasing sequence of precompact open subsets. 


\section{Conformal blow-ups}
\seclabel{\ConfBlowUp}
Given a sequence of $g_0$-anti-self-dual connections on a $G$ bundle $P$ over $X_0$
with curvature densities concentrating near a set of `singular points' in
$X_0$, we define associated sequences of mass centres and scales. In a
manner analogous to Chapter 3, we then obtain sequences of `conformal
blow-up maps' $f_{I\al}$ (defined exactly as in \S\ref{\secGluCon}) which
resolve these singularities in a sense that will be made precise below and
in \S4.3. As will become evident, the process of applying conformal
blow-ups may need to be iterated before the singularities are completely
`resolved'.

Let us commence by defining the first level conformal blow-ups. 
Suppose $\{A_\al\}_{\a=1}^\8$ is a sequence of $g_0$-anti-self-dual connections over $X_0$ with weak limit $(A_0,Z_0)$.
Let us consider the behaviour of the sequence $\{A_\al\}_{\a=1}^\8$ in $M_{X_0,k}(g_0)$ near the singular set $Z_0=\{x_i\}_{i=1}^{m_0}$ in more detail. If the point $x_i$ has multiplicity $k_i$, then 
$$
\lim_{r\to\8}\lim_{\a\to\8}\int\limits_{B(x_i,r)}
|F(A_\al)|_{g_0}^2\,dV_{g_0}=8\pi^2k_i.
\eqlabel\eqnlabel{\eqMassLimBall}
$$
Choose constants $d_0,r_0$ such that
$$
0<d_0\le\min_{i\ne j}\dist_{g_0}(x_i,x_j)\quad\hbox{and}\quad
0<r_0<\quarter\min\{1,\varrho_0,d_0\}. \eqlabel\eqnlabel{\eqD0R0Choice}
$$ 
We next define mass centres and scales of $g_0$-anti-self-dual connections restricted
to the fixed ball $B(x_i,r_0)\subset X_0$ by appropriately modifying the
previous definitions of mass centres and scales of \S\ref{\secConnSphere}
for $g_1$-anti-self-dual connections over $\SS^4$. First, note that
$$
\lim_{\a\to\8}\int\limits_{B(x_i,r_0)}
\left(|F(A_\al)|_{g_0}^2-|F(A_0)|_{g_0}^2\right)\,dV_{g_0}=8\pi^2k_i.
\eqlabel\eqnlabel{\eqMassDiffLim}
$$
Choose a frame $v_i$ in $FX_0|_{x_i}$ and let
$q=\phi_{x_i}^{-1}$ be the associated geodesic normal coordinate chart. 
For each $i$, define a sequence of {\it mass centres} 
$\{x_{i\al}\}_{\al=1}^\8$ in $B(x_i,r_0)$ by $x_{i\al}\equiv\phi_{x_i}(q_{i\al})$, where $q_{i\al}=\Centre[A_\al|_{B(x_i,r_0)}]\in\RR^4$ and
$$\leqalignno
{\Centre[A_\al|_{B(x_i,r_0)}]
&\equiv{1\over 8\pi^2k_i}\int\limits_{B(x_i,r_0)}
q\left(|F(A_\al)|_{g_0}^2-|F(A_0)|_{g_0}^2\right)\,dV_{g_0}. 
&\numeqn\eqnlabel{\eqMassCentre}\cr}
$$
Define a sequence of {\it scales} $\{\la_{i\al}\}_{\al=1}^\8$ in $(0,\8)$ by  
setting $\la_{i\al}=\Scale[A_\al|_{B(x_i,r_0)}]$, where
$$\leqalignno
{\Scale^2[A_\al|_{B(x_i,r_0)}]
&\equiv{1\over 8\pi^2k_i}\int\limits_{B(x_i,r_0)}
|q-q_{i\al}|^2\left(|F(A_\al)|_{g_0}^2-|F(A_0)|_{g_0}^2\right)\,dV_{g_0}. 
&\numeqn\eqnlabel{\eqScale}\cr}
$$
As in \S3.2., Eq. \eqref{\eqScale} leads to a {\it Tchebychev inequality}:
$$
\int\limits_{B(x_i,r_0)\less B(x_{i\al},R\la_{i\al})}
\left(|F(A_\al)|_{g_0}^2-|F(A_0)|_{g_0}^2\right)\,dV_{g_0}
\le 8\pi^2k_iR^{-2},\qquad R\ge 1.
\eqlabel\eqnlabel{\eqTchebyBall}
$$
Hence, if $R\gg 1$ and $\al$ is sufficiently large, the balls $B(x_{i\al},R\la_{i\al})$ contain most of the $8\pi^2k_i$ quantity of $A_\al$-energy bubbling off at $x_{i\al}$.

\rmk Other choices of scale function are possible. For example, we might
have chosen $\la_{i\al}$ to be the radius of the ball centred at $x_{i\al}$
containing $A_\al$-energy $8\pi^2(k_i-\half)$. As in [D83], a cutoff
function is required in order to regularise this definition.  
\medskip

Thus, we obtain a sequence of scales $\{\la_{i\al}\}_{\al=1}^\8$ associated
to the sequences of mass centres $\{x_{i\al}\}_{\a=1}^\8$ and connections
$\{A_\al\}_{\al=1}^\8$. Moreover, Eq. \eqref{\eqMassLimBall} implies that
the sequence $x_{i\al}$ converges to $x_i$ and that the sequence of scales
$\la_{i\al}$ converges to zero.  
Choose a sequence of frames $v_{i\al}\in FX_0|_{x_{i\al}}$ converging to
the frame $v_i\in FX_0|_{x_i}$ and let $\phi_{x_{i\al}}^{-1}$ be the
corresponding geodesic normal coordinate charts. 
Let $f_{x_{i\al}}\equiv\phi_{in}\circ
c_{\la_{i\al}}\circ\phi_{x_{i\al}}^{-1}$, where $c_{\la_{i\al}}$ is the
dilation of $\RR^4$ given by $x\mapsto x/\la_{i\al}$, let $\tilde g_{i\al}$
be the approximately round metric on $X_{i\al}'$ defined as in
\S\ref{\secMetBubbleTree}, 
let $P_{i\al}=(f_{x_{i\al}}^{-1})^*P$ be the induced $G$ bundle over
$X_{i\al}'$, 
and let $A_{i\al}=(f_{x_{i\al}}^{-1})^*A_\al$ be the induced $\tilde
g_{i\al}$-anti-self-dual connection on $P_{i\al}$. We call the maps
$f_{x_{i\al}}$ {\it conformal blow-ups}.

We obtain a sequence of open subsets $X_{i\al}'$ which exhaust $X_i\less\{x_{is}\}$, a sequence of metrics $\{\tilde g_{i\al}\}_{\a=1}^\8$,
and a sequence of $\tilde g_{i\al}$-anti-self-dual connections $\{A_{i\al}\}_{\a=1}^\8$ over the $X_{i\al}'$.
The sequence $\{\tilde g_{i\al}\}_{\a=1}^\8$ converges in $C^\8$ on compact
subsets of $X_i\less\{x_{is}\}$ to the standard round metric $g_i$ on
$X_i\equiv\SS^4$. Let $\{g_\al\}_{\al=1}^\8$ be the sequence of $C^\8$
metrics, defined as in \S\ref{\secMetBubbleTree}, on the connected sum
$X\equiv \#_{i=0}^{m_0}X_{i\al}'$, defined as in \S\ref{\secGluCon}, and
let $\{\check A_\al\}_{\al=1}^\8$ be the induced sequence of $g_\al$-anti-self-dual
connections over $X$. We call the connected sums $(X,g_\al)$ {\it conformal
blow-ups} of $(X_0,g_0)$.

There is a uniform upper bound on the $L^2$ norms 
$\|F(A_{i\al})\|_{L^2(X_{i\al}',\tilde g_{i\al})}$ since  
$$
\int_{X_{i\al}'}|F(A_{i\al})|_{\tilde g_{i\al}}^2\,dV_{\tilde g_{i\al}}
=\int_{B(x_{i\al},N\la_{i\al}^{1/2})}|F(A_\al)|_{g_0}^2\,dV_{g_0}
\le 8\pi^2(k_i+1/2),
\eqlabel\eqnlabel{\eqCrvSphLaUpBnd}
$$
for sufficiently large $\al$ by Eq. \eqref{\eqMassLimBall}, 
while Eqs. \eqref{\eqMassDiffLim} and \eqref{\eqTchebyBall} give a lower bound
$$
\int_{X_{i\al}'}|F(A_{i\al})|_{\tilde g_{i\al}}^2\,dV_{\tilde g_{i\al}}
=\int_{B(x_{i\al},N\la_{i\al}^{1/2})}|F(A_\al)|_{g_0}^2\,dV_{g_0}
\ge 8\pi^2(k_i-1/2).
\eqlabel\eqnlabel{\eqCrvSphLaLowBnd}
$$
Proposition \ref{\propMetExhstCpt} provides a subsequence $\{A_{i\al}\}_{\a=1}^\8$ which converges weakly to an ideal $g_i$-anti-self-dual connection $(A_i,Z_i)$ over $X_i$, where $Z_i=\{x_{ij}\}_{j=1}^{m_i}$. The energy bound of Eq. \eqref{\eqTchebyBall} ensures that $Z_i\subset X_i\less\{x_{is}\}$. 
Let $\mu_i=\mu(A_i,Z_i)$ be the associated singular measure on $X_i$ and
note that its mass may be computed by
$$
\int_{X_i}\,d\mu_i=\lim_{R\to\8}\lim_{\al\to\8}\int_{B(x_{in},R)}
|F(A_{i\al})|_{\tilde g_{i\al}}^2\,dV_{\tilde g_{i\al}}.
$$
Since this must be $8\pi^2$ times an integer, Eqs. \eqref{\eqCrvSphLaUpBnd} and \eqref{\eqCrvSphLaLowBnd} imply that $\mu_i$ has mass
$8\pi^2k_i$, where $k_i=\sum_{j=0}^{m_i}k_{ij}$, $A_i$ is a $g_i$-anti-self-dual
connection on a bundle $P_i$ over $X_i$ with $c_2(P_i)=k_{i0}$, and each
point $x_{ij}$ has multiplicity $k_{ij}$.

\rmk It is not strictly necessary that we construct a sequence of honest
metrics $g_\al$ over the connected sums $X=\#_{i=0}^{m_0}X_{i\al}'$ above;
a sequence of conformal structures $[g_\al]$ constructed as in
\S\ref{\secMetBubbleTree} would suffice and this would eliminate the need
for the choice of conformal factors over the necks. In any case, the
actual limits obtained are independent of such choices.

The above conformal blow-up construction produces a sequence of $\tilde
g_{x_{i\al}}$-anti-self-dual connections $A_{x_{i\al}}$ on increasing subsets
$X_{i\al}'$ of the four-sphere $X_i$ with weak $g_i$-anti-self-dual limit $(A_i,Z_i)$.
With the inverse process of gluing in mind, we describe a modified choice
of conformal blow-ups which yield {\it centred limits} $(\tilde A_i,\tilde
Z_i)$. First, a technical lemma concerning the variation of geodesic normal coordinate charts with their coordinate centres is required. The proof uses Taylor's theorem and is left to the reader.

\proclaim{Lemma}{\thmlabel{\lemNormCoordVar} Let $X_0$ be a closed $C^\8$ $n$-manifold with metric $g_0$ and injectivity radius $\varrho_0$. 
Let $x_0\in X$, let $v_0\in FX|_{x_0}$, and let $x=\exp_{v_0}^{-1}$ be the geodesic normal coordinate chart on $B(x_0,\varrho_0)$ defined by the frame $v_0$. Suppose $x_1\in B(x_0,\varrho_0/4)$ and $p=\exp_{v_0}^{-1}(x_1)$, so that $\dist_{g_0}(x_1,x_0)=|p|$. We now define two coordinate charts on $B(x_1,\varrho_0/2)$: 
\item{\rm (a)} Let $v_1\in FX|_{x_1}$ be the frame obtained by parallel translating $v_0$ along the geodesic joining $x_0$ to $x_1$, and let $w=\exp_{v_1}^{-1}$ on $B(x_1,\varrho_0/2)$; 
\item{\rm (b)} Let
$\tau_p$ be the translation on $\RR^n$ given by $q\mapsto q-p$, and let 
$\bar w=\tau_p\circ\exp_{v_0}^{-1}$ on $B(x_1,\varrho_0/2)$. Then the coordinates $\bar w$ converge to $w$ in $C^\8$ on $B(x_0,\varrho_0/4)$ as $p\to 0$: $|\bar w^\mu -w^\mu|=O(|w||p|)$, 
$|\rd\bar w^\mu/\rd w^\al -\de^\mu_\al|=O(p)$, and for all $m\ge 2$,
$\rd^m \bar w^\mu/\rd w^{\al_1}\cdots\rd w^{\al_m}= O(p)$.}
 
Next, we define the mass centre and scale of a positive Borel measure $\mu$ on $\RR^4$ by
$$\leqalignno
{p=\Centre[\mu]&\equiv\int_{\RR^4}x\,d\mu\quad\hbox{and}\quad
\la^2=\Scale^2[\mu]\equiv\int_{\RR^4}|x-p|^2\,d\mu. 
&\numeqn\cr}
$$
Let $\Theta$ be the product connection over $X_i$.
The proof of the following lemma describes how to choose conformal blow-ups
which produce centred limits.

\proclaim{Lemma}{\thmlabel{\lemCentredBlwUp}
Let $\{A_\al\}$ be a sequence of $g_0$-anti-self-dual connections over $X_0$ with weak
limit $(A_0,Z_0)$, where $Z_0=\{x_i\}_{i=1}^{m_0}$ is non-empty. Choose
$r_0$ as in\/ {\rm Eq. \eqref{\eqD0R0Choice}}. Then for each $x_i\in Z_0$,
the sequence $\{A_\al\}$ determines a sequence of points $\{w_{i\al}\}$ converging to $x_i$, a sequence of frames $v_{i\al}\in FX_0|_{w_{i\al}}$ converging to a frame $v_i\in FX_0|_{x_i}$, and a sequence of scales $\{\ka_{i\al}\}$ converging to zero such that the following holds. Fix $N>4$,
let $f_{w_{i\al}}$ be the corresponding sequences of conformal blow-ups, and let $A_{w_{i\al}}$ be the induced sequence of $\tilde g_{w_{i\al}}$-anti-self-dual connections  with weak $g_i$-anti-self-dual limit $(\tilde A_i,\tilde Z_i)$ over the four-sphere $X_i$. 
The limit $(\tilde A_i,\tilde Z_i)$ has the following properties:
\item{\rm (a)} If $\tilde A_i\ne \Theta$, then $\tilde A_i$ is centred; 
\item{\rm (b)} If $\tilde A_i=\Theta$, then the corresponding singular measure $\tilde \mu_i$ is centred.}

\pf (a) We begin by defining, exactly as before, a sequence of points $\{x_{i\al}\}$ converging to $x_i$, a sequence of frames $v_{i\al}\in FX_0|_{x_{i\al}}$ converging to a frame $v_i\in FX_0|_{x_i}$, and a sequence of scales $\{\la_{i\al}\}$ converging to zero. Let $f_{x_{i\al}}$ be the corresponding sequences of conformal blow-ups and let $A_{x_{i\al}}$ be the induced sequence of $\tilde g_{x_{i\al}}$-anti-self-dual connections with weak $g_i$-anti-self-dual limit $(A_i,Z_i)$ over $X_i$. Suppose $\Centre[A_i]=p_i$ and $\Scale[A_i]=\nu_i$. 

\case{1.} $Z_i=\emptyset$. Recall that $f_{x_{i\al}}=\phi_{in}\circ c_{\la_{i\al}}\circ\phi_{x_{i\al}}^{-1}$, $A_{x_{i\al}}=(f_{x_{i\al}}^{-1})^*A_\al$, and $\tilde g_{x_{i\al}}=\la_{i\al}^{-2}(f_{x_{i\al}}^{-1})^*g_0$. Define  
$h_i=\phi_{in}\circ c_{\nu_i}\circ\tau_{p_i}\circ\phi_{in}^{-1}$ and set
$\bar f_{w_{i\al}}=h_i\circ f_{x_{i\al}}$. Then 
$$\leqalignno
{\bar f_{w_{i\al}}&=\phi_{in}\circ c_{\la_{i\al}\nu_i}\circ\tau_{p_i\la_{i\al}}\circ\phi_{x_{i\al}}^{-1}
=\phi_{in}\circ c_{\ka_{i\al}}\circ\bar\phi_{w_{i\al}}^{-1}, &\cr}
$$ 
where $w_{i\al}\equiv\phi_{x_{i\al}}(p_i\la_{i\al})$, 
$\ka_{i\al}\equiv\la_{i\al}\nu_i$, and
$\bar\phi_{w_{i\al}}\equiv\phi_{x_{i\al}}\circ\tau_{p_i\la_{i\al}}^{-1}$.
Thus, $\bar f_{w_{i\al}}$ provides a diffeomorphism from the small ball $B(w_{i\al},N\ka_{i\al}^{1/2})$ in $X_0$ to the open subset 
$B(x_{in},N\ka_{i\al}^{-1/2})$ of $X_i$. The sequence of points $\{w_{i\al}\}$ converges to $x_i$ and the sequence of scales $\{\ka_{i\al}\}$ converges to zero. As in \S3.5, define a sequence of metrics on the increasing subsets $B(x_{in},N\ka_{i\al}^{-1/2})$ by
$\bar g_{w_{i\al}}\equiv\ka_{i\al}^{-2}h_1^2(\bar f_{w_{i\al}}^{-1})^*g_0$. 
Then $\bar g_{w_{i\al}}$ converges to the standard metric $g_i$ in $C^\8$ on compact subsets of $X_i\less\{x_{is}\}$. 
Define a sequence of $\bar g_{w_{i\al}}$-anti-self-dual connections over the balls $B(x_{in},N\ka_{i\al}^{-1/2})$ by 
$\bar A_{w_{i\al}}\equiv(\bar f_{w_{i\al}}^{-1})^*A_\al$, 
and observe that 
$\bar A_{w_{i\al}}=(h_i^{-1})^*A_{x_{i\al}}$. The sequence $\{\bar A_{w_{i\al}}\}$ converges to the centred connection $(h_i^{-1})^*A_i$ in $C^\8$ on compact subsets of  $X_i\less\{x_{is}\}$.

It remains to replace the chart 
$\bar w\equiv \bar\phi_{w_{i\al}}^{-1}$ on $B(w_{i\al},\varrho_0/2)$ by a geodesic normal coordinate chart $w\equiv \phi_{w_{i\al}}^{-1}$. Choose a frame $v'_{i\al}\in FX_0|_{w_{i\al}}$ by parallel translating the frame $v_{i\al}\in FX_0|_{x_{i\al}}$ along the geodesic connecting $x_{i\al}$ and $w_{i\al}$, noting that $\dist_{g_0}(x_{i\al},w_{i\al})=|p_i|\la_{i\al}$. Thus, as $\al\to\8$, the coordinate chart $\bar w$ converges in $C^\8$ on $B(x_i,\varrho_0/4)$ to the geodesic normal coordinate chart $w$ in the sense of Lemma \ref{\lemNormCoordVar}. Define a new sequence of conformal blow-up maps by setting
$f_{w_{i\al}}=\phi_{in}\circ c_{\ka_{i\al}}\circ\phi_{w_{i\al}}^{-1}$, 
and define corresponding sequences of connections and metrics on the balls $B(x_{in},N\ka_{i\al}^{-1/2})$ by 
$A_{w_{i\al}}=(f_{w_{i\al}}^{-1})^*A_\al$ and 
$\tilde g_{w_{i\al}}=\ka_{i\al}^{-2}h_1^2(f_{w_{i\al}}^{-1})^*g_0$.
Lemma \ref{\lemNormCoordVar} implies that the 
sequences $\{\tilde g_{w_{i\al}}\}$ and $\{A_{w_{i\al}}\}$
converge in $C^\8$ on compact subsets of $X_i\less\{x_{is}\}$ to the metric $g_i$ and centred $g_i$-anti-self-dual connection $\tilde A_i\equiv (h_i^{-1})^*A_i$. This completes the proof of (a) in Case 1. 

\case{2.} $Z_i\ne\emptyset$. The proof is similar to that of Case 1. Let $\tilde Z_i=h_i^{-1}(Z_i)$. Then the sequences $\{\bar A_{w_{i\al}}\}$ and
$\{A_{w_{i\al}}\}$ converge in $C^\8$ on compact subsets of 
$X_i\less(\tilde Z_i\cup\{x_{is}\})$ to the centred connection $(h_i^{-1})^*A_i$.

\noindent (b) One sets $\Centre[\mu_i]=p_i$, $\Scale[\mu_i]=\nu_i$, and essentially repeats the proof of Part (a) for the sequence of measures $\mu_{x_{i\al}}\equiv |F(A_{x_{i\al}})|_{g_0}^2$. \qed

\rmk In the sequel, we require that the conformal blow-up maps be chosen as
in Lemma \ref{\lemCentredBlwUp}. However, to conserve notation, we will relabel the points $w_{i\al}$ and scales $\ka_{i\al}$ by $x_{i\al}$ and $\la_{i\al}$, respectively, and the limit $(\tilde A_i,\tilde Z_i)$ by $(A_i,Z_i)$.

A technical point that we have not addressed above is that, just as in [P-W],
the weak limit of the sequence $\{A_{i\al}\}$ apparently depends on certain
choices of parameters in the conformal blow-up construction: 

\subsubsection{(1) Neck width parameter $N$} This was only included in this
Chapter for the sake of consistency with the gluing construction of
Chapters 3 and 5: we could just as well have set $N=2$, say. 

\subsubsection{(2) Radius $r_0$} Following [P-W], the dependency is removed
by letting $r_0\to 0$. The conformal blow-up process gives a sequence of
points $\{x_{i\al}(r_0)\}$, scales $\{\la_{i\al}(r_0)\}$, blow-up maps 
$\{f_{x_{i\al}}(r_0)\}$, metrics $\{\tilde g_{i\al}(r_0)\}$, and
connections $\{A_{i\al}(r_0)\}$. The sequence of connections
$\{A_{i\al}(r_0)\}$ converges to an ideal $g_i$-anti-self-dual limit
$(A_i(r_0),Z_i(r_0))$, for any fixed $r_0>0$. We now let $r_0\to 0$ and by
a standard diagonal argument, we obtain a weakly convergent subsequence
$\{A_{i\al}(r_0)\}$ with weak limit $(A_i,Z_i)$, say. 

\subsubsection{(3) Frames $v_{i\al}$ and $v_i$} The construction is
$\SO(4)$ equivariant: Rotating the frames $v_{i\al}\in FX|_{x_{i\al}}$ and
$v_i\in FX_{x_i}$ by elements of $\SO(4)$ induces an $SO(4)$ action on the 
connections $A_{i\al}$ and $A_i$ as described in \S\ref{\secConnSphere}.  

There is one final issue which will be important in our later discussion
of alternative modes of convergence for sequences of anti-self-dual connections: we
must exclude the possibilty that curvature is lost over the necks $\Om_i$
arising in the conformal blow-up process described above. Of course, the
curvature can only bubble off with masses equal to an integer multiple of
$8\pi^2$, so it suffices to show that we can choose the neck parameters to
ensure that the curvature masses over the necks are strictly less than
$8\pi^2$. So, consider again the sequence $\{A_\al\}_{\al=1}^\8$ of $g_0$-anti-self-dual
connections 
over $X_0$ with weak limit $(A_0,Z_0)$, where $Z_0=\{x_i\}_{i=1}^{m_0}$,
and let $\{A_{i\al}\}_{\al=1}^\8$ be the corresponding sequences of  
$\tilde g_{i\al}$-anti-self-dual connections over $X_{i\al}'$ having weak limits
$(A_i,Z_i)$, where $Z_i=\{x_{ij}\}_{j=1}^{m_j}$. Let
$\{\la_{i\al}\}_{\al=1}^\8$ be the sequence of scales associated to the
sequence of connections $\{A_\al\}_{\al=1}^\8$ and the singular point
$x_i\in Z_0$. Given this set-up, standard arguments yield the
following curvature estimates near $x_i$:
 
\proclaim{Lemma}{\thmlabel{\lemAnnFEst} 
Given $\eps>0$, there exist positive constants $R_0$, $r_1$, and $\al_0$ with the following significance.
For large enough $R_0$, small enough $r_1$ and large enough $\al_0$, then $R_0\la_{i\al}<r_1$ for any $\al\ge\al_0$ and the following holds.
\item{\rm (a)} $|\|F(A_{i\al})\|^2_{L^2(B(x_i,R_0),\tilde g_{i\al})}
-8\pi^2k_i|<\eps^2$,
\item{\rm (b)} $|\|F(A_\al)\|^2_{L^2(B(x_i,R_0\la_{i\a}),g_0)}
-8\pi^2k_i|<\eps^2$,
\item{\rm (c)} $\|F(A_{i\al})\|_{L^2(\W(x_i,R_0,r_1\la_{i\a}^{-1}),
\tilde g_{i\al})}< \eps$,
\item{\rm (d)} $\|F(A_\al)\|_{L^2(\W(x_i,R_0\la_{i\a},r_1),g_0)}< \eps$.}

Thus, we have the following curvature estimate which ensures that in the
limit there is no `curvature loss' over the necks $\Om_i$. (In particular, 
if $A_{i\al}$
converges weakly to $(A_i,Z_i)$, then the singular set $Z_i\subset X_i$
does not contain the south pole $x_{is}$.)

\proclaim{Corollary}{\thmlabel{\CorNeckFEst} 
Given $\eps>0$ and $N>4$, there is an $\al_0>0$ with the following
significance. If
$\Om_{i\al}\equiv\Om(x_{i\al},N^{-1}\la_{i\al}^{1/2},N\la_{i\al}^{1/2})$
and, $B_{i\al}'\equiv B(x_{i\al},N\la_{i\al}^{1/2})$, then for any
$\al\ge\al_0$, we have  
\item{\rm (a)} $\|F(A_\al)\|_{L^2(\Om_{i\al},g_0)}< \eps$, and
\item{\rm (b)} $|\|F(A_\al)\|_{L^2(B_{i\al}',g_0)}^2-8\pi^2k_i|< \eps$.}

Lastly, we note that the conformal blow-up process may of course be
iterated if the singular sets $Z_i$ 
are non-empty. In the next section we show that after repeating the
conformal blow-up process at most $k$ times, we obtain a sequence of
$g_\al$-anti-self-dual connections $\{\check A_\al\}$ which is strongly convergent.
Indeed, given the weakly convergent sequence $\{A_{i\a}\}_{\a=1}^\8$ over
the $X_{i\al}'$ near a point $x_{ij}$ with multiplicity $k_{ij}$ in the
singular set $Z_i\subset X_i$, 
the second-level process differs from the first-level only in minor
technical details:
We define sequences of centres
$x_{ij\al}=\phi_{x_{ij}}(q_{ij\al})$ converging to $x_{ij}$
and scales $\la_{ij\al}$ converging to zero, now using the metrics
$\tilde g_{i\al}$
and a coordinate chart $\phi_{x_{ij}}$ on $X_i$ given by
$\phi_{x_{ij}}=\phi_{in}\circ\tau_{q_{ij}}^{-1}$ where 
$\phi_{in}(q_{ij})=x_{ij}$. The blow-up maps are then defined using
coordinate charts on $X_i$ given by
$\phi_{x_{ij\al}}=\phi_{in}\circ\tau_{q_{ij\al}}^{-1}$ and setting
$f_{x_{ij\al}}=\phi_{ijn}\circ c_{\la_{ij\al}}\circ\phi_{x_{ij\al}}^{-1}$.
We then proceed exactly as before and similarly for all higher-level blow-ups.


\section{Bubble tree compactification} 
By analogy with the arguments of [T88, \S5] and [P-W], we define a bubble tree compactification for the moduli space  $M_{X_0,k}(g_0)$ of anti-self-dual connections.
First, we need an appropriate notion of an ``ideal connection'':

\defn\thmlabel{\defTreeIdealConn} A {\it bubble tree ideal 
$g_0$-anti-self-dual connection $A$ of second Chern class $k$} over $X_0$ is determined by the following data.
\item{(a)} An oriented tree $\script{I}$ with a finite set of vertices $\{I\}$, including a base vertex $0$, and a set of edges $\{(I_-,I)\}$. Each vertex $I$ is labelled with an integer $k_I\ge 0$ such that 
\itemitem{(i)} $\sum_{I\in\sI}k_I=k$, 
\itemitem{(ii)} If $I>0$ is a terminal vertex, then $k_I>0$,
\itemitem{(iii)} There are at most $k$ terminal vertices, excluding the base vertex.
\item{(b)} A $(2m-1)$-tuple $(A_I,x_I)_{I\in\script{I}}$, where $m$ is the number of vertices in $\sI$.
\item{(c)} If $I=0$, then
$A_0$ is a $g_0$-anti-self-dual connection on a $G$ bundle $P_0$ over $X_0$ with 
$c_2(P_0)=k_0\ge 0$.
\item{(d)} If $I>0$, then 
\itemitem{(i)} $A_I$ is either the product connection $\Theta$ or
a centred $g_I$-anti-self-dual connection on a $G$ bundle $P_I$ over the sphere $X_I\equiv\SS^4$ with $c_2(P_I)=k_I$, where $g_I$ is the standard round metric, 
\itemitem{(ii)} $x_I$ is a point in $X_0$ if $I_-=0$ and a point in $X_{I_-}\less\{x_{Is}\}$ if $I_->0$.
\item{(e)} If $I>0$ and $A_I=\Theta$, then there are at least $2$ outgoing edges emanating from that vertex.
\medskip

Definition {\defTreeIdealConn} should be compared with the 
construction of approximately anti-self-dual connections in \S
\ref{\secGluCon}.  
The ideal connection $(A_I,x_I)_{I\in\script{I}}$ is often written as
$(A_I)_{I\in\script{I}}$. 
Heuristically, we may view an ideal $g_0$-anti-self-dual connection $A=(A_I)_{I\in\script{I}}$ as a `connection' over the join $\vee_{I\in\script{I}} X_I$, where each sphere $X_I$ is attached to the lower level $X_{I_-}$ by identifying the south pole $x_{Is}$ with the point $x_I\in X_{I_-}$. Let $Z_{I_-}\subset X_{I_-}$ denote the set of `attachment points' $x_I$ in $X_0$, if $I_-=0$, or points $x_I$ in $X_{I_-}\less\{x_{Is}\}$, if $I_->0$. Let $m_I$ be the number of points in $Z_I$, i.e., the number of outgoing edges emanating from vertex $I$.
 
Second, we need an appropriate notion of convergence. Let $X\equiv\#_{I\in\script{I}}X_I$ be the connected sum defined in \S\ref{\secGluCon} by a set of scales $\{\la_{I\al}\}_{I\in\script{I}}$, with $\ol{\la}_\al\to 0$ as $\al\to\8$, and a fixed neck parameter $N$. Similarly, if $\{g_\al\}$ is the corresponding sequence of $C^\8$ metrics on $X$ defined in \S\ref{\secMetBubbleTree}, then $g_\al$ converges to $g_I$ in $C^\8$ on compact subsets of $X_I\less(Z_I\cup\{x_{Is}\})$ for each $I\ge 0$.
Following [D-K, \S7.3.1], we consider the following modes of convergence for sequences of anti-self-dual connections over $X$.

\defn\thmlabel{\defStrongConv} Let $\{A_\al\}_{\al=1}^\8$ be a sequence of $g_\al$-anti-self-dual connections on a fixed bundle $P$ with $c_2(P)=k$ over the connected sum $X=\#_{I\in\script{I}}X_I$. 
\item{(a)} If $Y\in s^k(X)$ is a multiset in 
$\cup_{I\in\script{I}}X_I\less (Z_I\cup \{x_{Is}\})$, the sequence $\{A_\al\}$ {\it converges  weakly} to
$((A_I,x_I)_{I\in\script{I}},Y)$ if 
the gauge equivalence classes $[A_\a]$ converge in $C^\8$ to $([A_I])_{I\in\script{I}}$ over compact subsets of 
$\cup_{I\in\script{I}}X_I\less (Z_I\cup \{x_{Is}\}\cup Y)$
and if the curvature densities converge,
$$
|F(A_\a)|_{g_\al}^2\too \sum_{I\in\script{I}}|F(A_I)|_{g_I}^2+8\pi^2\de_Y,  
$$
over compact subsets of 
$\cup_{I\in\script{I}}X_I\less (Z_I\cup \{x_{Is}\})$.

\item{(b)} The sequence $\{A_\a\}$ {\it converges strongly} to the limit
$(A_I,x_I)_{I\in\script{I}}$
if it converges weakly to $(A_I,x_I)_{I\in\script{I}}$
(with no singular set $Y$) and if $\sum_{I\in\script{I}}c_2(P_I)=c_2(P)$.
Here, the $A_I$ are $g_I$-anti-self-dual connections on $G$ bundles $P_I$ over $X_I$
with $c_2(P_I)=k_I$. 
\medskip

We let $BM_{X_0,k}(g_0)$ denote the set 
bubble tree ideal $g_0$-anti-self-dual connection over $X_0$ of total second Chern
class $k$. Thus, each point of $BM_{X_0,k}(g_0)$ is represented by a  
$(2m-1)$-tuple $(A_I,x_I)_{I\in\sI}$, with $m$ being the total number of
vertices of the tree $\sI$.

\defn\thmlabel{\defTreeCpt} We say that a sequence $\{A_\al\}_{\al=1}^\8$ of $g_0$-anti-self-dual connections on a $G$ bundle $P$ over $X_0$ with $c_2(P)=k$ {\it converges strongly to a bubble tree ideal $g_0$-anti-self-dual connection} $(x_I,A_I)_{I\in\script{I}}$ in $BM_{X_0,k}(g_0)$ if there exist sequences of conformal blow-ups $\{f_{I\al}\}_{I\in\script{I}}$ with the following property. Let $\{g_\al\}$ be the induced sequence of $C^\8$ metrics in the conformal class $[g_0]$ on the connected sum $X=\#_{I\in\script{I}}X_I$. Let 
$\{\check A_\al\}$ denote the induced sequence of $g_\al$-anti-self-dual connections over $X$. Then we require that the sequence of metrics $\{g_\al\}$ converges in $C^\8$ on compact sets of $X_I\less(Z_I\cup\{x_{Is}\})$ to the metric $g_I$, $I\ge 0$, and
that the sequence of connections $\{\check A_\al\}$ converges strongly to the ideal $g_0$-anti-self-dual connection $(A_I,x_I)_{I\in\script{I}}$.
\medskip

This definition of convergence extends to the space of bubble tree ideal connections $BM_{X_0,k}(g_0)$, which is then endowed with a second countable Haussdorf topology. Define the {\it bubble tree compactification} $\ol{M}^\tau_{X_0,k}(g_0)$ to be the closure of $M_{X_0,k}(g_0)$ in $BM_{X_0,k}(g_0)$.

\proclaim{Theorem}{\thmlabel{\thmTreeCpt} The space $\ol{M}^\tau_{X_0,k}(g_0)$ is compact.}

The result follows from the special case below.

\proclaim{Theorem}{Any infinite sequence in $M_{X_0,k}(g_0)$ has a strongly convergent subsequence with limit point in $\ol{M}^\tau_{X_0,k}(g_0)$.}

\pf The argument is similar to the proof of Proposition 5.3 in [T88].
Fix a $G$ bundle $P$ over $X_0$ with $c_2(P)=k>0$ and let $\{A_\al\}_{\al=1}^\8$ be a sequence of $g_0$-anti-self-dual connections on $P$. 
The main point is to repeatedly apply conformal blow-ups $f_{I\al}$ until we obtain a sequence of induced metrics $g_\al$ over a connected sum $X$, with $(X,g_\al)$ conformally equivalent to $(X_0,g_0)$, and a sequence of induced $g_\al$-anti-self-dual connections over $X$, denoted by $\{\check A_\al\}$, which is strongly convergent. We adopt the convention below that subsequences are immediately relabelled. 
 
\step{1.} There is a subsequence $\{A_\al\}$ which converges weakly to an ideal $g_0$-anti-self-dual connection $(A_0,Z_0)$, with $Z_0=\{x_i\}_{i=1}^{m_0}$ corresponding to a point in the symmetric product $s^k(X_0)$. If $Z_0=\emptyset$ then we are done, so assume that $m_0\ge 1$. Let $k_i$ be the multiplicity of $x_i$ and note that $0<k_i\le k$. For each $i$ and large enough $\al$, the connection $A_\al$ determines a set of mass centres $\{x_{i\al}\}_{i=1}^{m_0}$, with $x_{i\al}\to x_i$, and a set of scales $\{\la_{i\al}\}_{i=1}^{m_0}$, with 
$\la_{i\al}\to 0$ as $\al\to\8$. Fix a neck width parameter $N>4$, choose a sequence of frames $v_{i\al}\in FX_0|_{x_{i\al}}$ converging to a frame $v_i\in FX_0|_{x_i}$, and let 
$\{f_{i\al}\}_{i=1}^{m_0}$ be the conformal blow-up maps defined by these centres, frames, scales, and parameter $N$. 
If $X=\#_{i=0}^{m_0} X_{i\al}'$, then $(X,g_\al)$ is the conformal blow-up of $(X_0,g_0)$ determined by the maps $f_{i\al}$. 
Let $P$ now denote the induced $G$ bundle over $X$, let $\check A_\al$ denote the induced $g_\al$-anti-self-dual connection over $X$, and let $A_{i\al}$ be the restriction of $\check A_\al$ to the open subset $X_{i\al}'$.

The sequence $[A_{i\al}]$ has a weakly convergent subsequence, again denoted $[A_{i\al}]$, with weak limit $(A_i,Z_i)$, where $Z_i$ corresponds to a point in $s^{k_i}(X_i)$. Corollary \ref{\CorNeckFEst} implies that no mass is lost over the neck $\Om_i$. Hence, if each $Z_i=\emptyset$, $i>0$, then we have $\sum_{i=0}^{m_i}k_i=k$, the sequence $[A_{i\al}]$ converges strongly to $[A_i]$, and we proceed to the Final Step. Otherwise, $Z_i\ne\emptyset$ for some $i$ and we proceed to Step $2$.

\step{2.} For some $i$, Step 1 produces a non-empty singular set
$Z_i=\{x_{ij}\}_{j=1}^{m_i}$. Let $k_{ij}$ be the multiplicity of the point $x_{ij}$, let $c_2(A_i)=k_{i0}$, and note that $\sum_{j=0}^{m_i} k_{ij}=k_i>0$.
Let $\mu_i$ be the singular measure associated with $(A_i,Z_i)$. We now consider two cases, depending on whether or not $A_i$ is the flat product connection $\Theta$ over $X_i$.

\case{(a)} $A_i=\Theta$. Since $\Scale[\mu_i]=1$,
the diameter of the set $Z_i$ must be positive and so this case can only occur if $m_i>1$. Let $k_{ij}$ be the multiplicity of the point $x_{ij}$ and note that as $m_i>1$ we must have $\max_j k_{ij}\le k-1$.

\case{(b)} $A_i\ne\Theta$. Therefore, $k_{i0}=c_2(A_i)>0$ and so we again must have $\max_j k_{ij}\le k-1$, since $\sum_{j=0}^{m_i} k_{ij}=k_i\le k$. 

For large enough $\al$, the connection $A_{i\al}$ determines a set of mass centres $\{x_{ij\al}\}_{i=1}^{m_i}$, with $x_{ij\al}\to x_{ij}$, and a set of scales $\{\la_{ij\al}\}_{i=1}^{m_i}$, with 
$\la_{ij\al}\to 0$ as $\al\to\8$. Let 
$\{f_{ij\al}\}_{j=1}^{m_i}$ be the conformal blow-up maps defined by these centres, scales, and parameter $N$. Let $P$ denote the induced $G$ bundle over the new connected sum $X=\#_{i=0}^{m_0} X_{i\al}'\#_{j=1}^{m_i} X_{ij\al}'$, let $\check A_\al$ denote the induced $g_\al$-anti-self-dual connection over $X$, and let $\{A_{ij\al}\}$ be the induced sequence of $g_\al$-anti-self-dual connections over the open subsets $X_{ij\al}'$ of the spheres $X_{ij}$.  

The sequence $[A_{ij\al}]$ has a weakly convergent subsequence with weak limit $(A_{ij},Z_{ij})$, with no loss of mass over the necks $\Om_{ij\al}$. If each $Z_{ij}=\emptyset$, $j=1,\dots,m_i$, then we have $\sum_{j=0}^{m_i} k_{ij}=k_i$, the sequence $[A_{ij\al}]$ converges strongly to $[A_{ij}]$, and the blow-up process terminates at the vertices $A_{ij}$. Otherwise, $Z_{ij}\ne\emptyset$ for some $j$ and we proceed to Step $3$.

\step{l: $3\le l\le k$.} For some multi-index $I$ of length $|I|=l-1$, Step $l-1$ produces a non-empty singular set
$Z_I=\{x_{Ij}\}_{j=1}^{m_I}$ contained in the sphere $X_I$. The sequence $[A_{I\al}]$ has a weak limit $(A_I,Z_I)$, where $Z_I$ corresponds to a point in $s^{k_I}(X_I)$. Let $k_{Ij}$ be the multiplicity of the point $x_{Ij}$, let $c_2(A_I)=k_{I0}$, and note that $\sum_{j=0}^{m_I} k_{Ij}=k_I>0$. Let $\mu_I$ be the singular measure associated with $(A_I,Z_I)$. 

\case{(a)} $A_I=\Theta$. Since $\Scale[\mu_I]=1$, 
the diameter of the set $Z_I$ must be positive. Hence, $m_I>1$ and so we have 
$$
\max_j k_{Ij}\le k-l+1,\qquad |Ij|=l,\quad 1\le l\le k.\eqlabel\eqnlabel{\eqMassDecr}
$$

\case{(b)} $A_I\ne\Theta$. Therefore, $k_{I0}=c_2(A_I)>0$ and so  
Eq. \eqref{\eqMassDecr} again holds, since $\sum_{j=0}^{m_I} k_{Ij}=k_I\le k$. 

Eq. \eqref{\eqMassDecr} implies that the conformal blow-up process terminates completely after at most $k$ steps.  

For large enough $\al$, the connection $A_{I\al}$ determines a set of mass centres $\{x_{Ij\al}\}_{j=1}^{m_I}$ in $X_I\less\{x_{Is}\}$, with $x_{Ij\al}\to x_{Ij}$, and a set of scales $\{\la_{Ij\al}\}_{j=1}^{m_I}$, with 
$\la_{Ij\al}\to 0$ as $\al\to\8$. Let 
$\{f_{Ij\al}\}_{j=1}^{m_I}$ be the conformal blow-up maps defined by these centres, scales, and parameter $N$. Let $P$ denote the induced $G$ bundle over the connected sum $X=\#_I X_{I\al}'\#_{j=1}^{m_i} X_{Ij\al}'$, let $\check A_\al$ denote the induced $g_\al$-anti-self-dual connection over $X$, and let
$\{A_{Ij\al}\}$ be the induced sequence of $g_\al$-anti-self-dual connections over the open subsets $X_{Ij\al}'$ of the spheres $X_{Ij}$.  

The sequence $[A_{Ij\al}]$ has a weakly convergent subsequence with weak limit $(A_{Ij},Z_{Ij})$, with no loss of mass over the necks $\Om_{Ij\al}$. If each $Z_{Ij}=\emptyset$, $j=1,\dots,m_I$, then we have $\sum_{j=0}^{m_I} k_{Ij}=k_I$, the sequence $[A_{Ij\al}]$ converges strongly to $[A_{Ij}]$, the blow-up process terminates at the vertices $A_{Ij}$, and we proceed to the Final Step. Otherwise, proceed to Step $l+1$.

\noindent{\bf Final Step.} After performing at most $k$ conformal blow-ups, we obtain a sequence of $g_\al$-anti-self-dual connections $\{\check A_\al\}$ over a connected sum $X=\#_{I\in\script{I}} X_{I\al}'$. The sequence $\{\check A_\al\}$ converges strongly to a bubble tree limit $(A_I,x_I)_{I\in\script{I}}$, since the singular points have all been blown up and there has been no mass loss over the necks $\Om_{I\al}$.
\qed

Plainly, the compactification $\ol{M}^\tau_{X_0,k}(g_0)$ is ``larger'' than
the Uhlenbeck compactification $\ol{M}^u_{X_0,k}(g_0)$. Indeed, there is an
obvious surjective map 
$$
\pi:\ol{M}^\tau_{X_0,k}(g_0)\too \ol{M}^u_{X_0,k}(g_0) \eqlabel
$$
obtained by sending a bubble tree ideal connection $(A_I,x_I)_{I\in\sI}$ to
the corresponding Uhlenbeck ideal connection $(A_0,x_1,\dots,x_{m_0})$. The
multiplicity of $x_i\in X_0$ is the sum of the second Chern classes of the
anti-self-dual connections $A_I$ attached to the subtree lying above the
vertex $i$. 

\proclaim{Corollary}{The map $\pi:\ol{M}^\tau_{X_0,k}(g_0)\to
\ol{M}^u_{X_0,k}(g_0)$ is continuous.} 


\section{$D_q$ convergence and strong convergence}
\seclabel{\secDqStrongConv} 
We will need one further notion of convergence in order to show that every point of the moduli space $M_{X,k}(g)$ lies in the image of the gluing map $\sJ$ constructed in Chapter 5.
Let $P$ be a $G$ bundle over a closed manifold $X$ with metric $g$. 
Following [D-K, \S 7.2.4], fix $4\le q<\8$ and let $D_q$ be the metric on the space $\cB_{X,P}$ given by
$$
D_q([A],[B])=\inf_{u\in\cG}\|A-u(B)\|_{L_q(X,g)}. \eqlabel
$$
We recall the following definition of Donaldson and Kronheimer.

\defn\thmlabel{\defLqConv}{\rm [D-K, p. 308]} Let $\{\l_{I\a}\}_{\a=1}^\8$, for each $I>0$, be sequences of scales satisfying $\ol\l_\a\to 0$, where $\ol\l_\a=\max_I\l_{I\a}$,
and let $\{A_\a\}_{\a=1}^\8$ be a sequence of connections on a fixed $G$ bundle $P\to X$, where $X\equiv\#_{I\in \sI} X_I'$ and $X_I={\Bbb S}^4$ if $I>0$. The connected sum $X$ has a sequence of metrics $\{g_\a\}_{\a=1}^\8$ defined by the sequence of scales $\{\l_{I\a}\}_{\a=1}^\8$, a sequence of points $\{x_{I\a}\}_{\a=1}^\8$, where the 
$x_{I\al}$ converge with respect to the fixed metric $g_I$ to a point $x_I\in X_{I_-}$, and a neck width parameter $N$.
Assume that the connections $A_\a$ are $g_\a$-ASD with respect to the sequence of metrics $\{g_\a\}_{\a=1}^\8$ on $X$. Then   
the sequence $\{A_\a\}_{\a=1}^\8$ is {\it $D_q$-convergent} to $(A_I,x_I)_{I\in\sI}$ if $D_q([A_\a|_{X_I''}], [A_I|_{X_I''}])\to 0$ as $\a\to\8$. 

$D_q$ convergence is called ``$L^q$ convergence'' in [D-K].
The result below explains the relationship between strong convergence and
$D_q$ convergence. 

\proclaim{Theorem}{\thmlabel{\thmStrongLqConv}{\rm [D-K, p. 309]} Let $\{A_\a\}_{\a=1}^\8$ be a sequence of connections on a bundle $P\to X$ which are ASD with respect to the sequence of metrics $\{g_\a\}_{\a=1}^\8$ determined by the sequences of scales $\{\l_{I\a}\}$, where $\ol\l_\a\to 0$. Then the sequence $\{A_\a\}_{\a=1}^\8$ is strongly convergent if and only if it is $D_q$-convergent.\qed} 

\endchapter
\eqnno=0
\secno=0
\thmno=0


\chapter{Differentials of the gluing maps} In this Chapter we obtain $L^2$
estimates for the differentials of the gluing maps $\hat\sJ:\sT/\Ga\to
M^*_{X_0,k}$. These give $C^0$ bounds for the components of the $L^2$ metric
${\bf g}$ on the bubbling ends of $M^*_{X_0,k}(g_0)$ and allow us to
complete the proofs of Theorems \ref{\thmFinVolDiam} and
\ref{\thmMetCompletion}. In particular, {\it for the remainder of the
article, the hypotheses of Theorem \ref{\thmFinVolDiam} are assumed to be
in effect}.   

\section{Construction of the gluing maps}
In this section we construct the gluing maps $\sJ:\sT/\Ga\to M^*_{X,k}(g)$
and $\hat\sJ:\sT/\Ga\to M^*_{X_0,k}(g_0)$, and set up the analytical
framework required for the later sections.
Our first task is to construct a right inverse to the linear 
operator $d_{A'}^{+,g}$ and so we choose suitable Sobolev spaces $L^q$,
$L^p_1$ and for the remainder of this Chapter, fix
$$
2\le p<4\quad\hbox{and}\quad 4\le q<\8\quad\hbox{so that}\quad
1/4+1/q=1/p.\eqlabel 
$$
By hypothesis,
$H^2_{A_I}=0$ for all $I$ and thus the operators $d^{+,g_I}_{A_I}$ have
right inverses $P_I$.
More explicitly, if $\D_{A_I}^{+,g_I}$ is the 
Laplacian $d_{A_I}^{+,g_I}(d_{A_I}^{+,g_I})^*$
and $G_{A_I}^{+,g_I}$ is the corresponding Green's
operator, we may set $P_I=(d_{A_I}^{+,g_I})^*G_{A_I}^{+,g_I}$.
A standard application of the Calderon-Zygmund theory and the Sobolev
inequalities gives the following bounds. 

\proclaim{Lemma}{\thmlabel{\lemPIEst}
Assume $H^2_{A_I}=0$. Then the operators
$P_I:L^p\to L^p_1$ and $P_I:L^p\to L^q$ are bounded and there are constants
$C_i=C_i(A_I,g_I,p)$, $i=1,2$, such that 
$$\leqalignno
{\|P_I\xi\|_{L^q(X_I,A_I,g_I)}
&\le C_1\|P_I\xi\|_{L^p_1(X_I,g_I)}\le C_2\|\xi\|_{L^p(X_I,g_I)},
\qquad \xi\in L^p\Om^2(X_I,\ad P_I). &\cr}
$$}

We next define the $C^\8$ cut-off functions to be used in the construction
of a right parametrix $Q$ for $d_{A'}^{+,g}$ by patching together the
operators $P_I$ over $X$.  

\proclaim{Lemma}{{\rm [Lemma 7.2.10, D-K], [D-S, p. 221]}\thmlabel{\lemBetaEst}
For any $\l>0$ and $N>4$, there exists a $C^\8$ function 
$\b_{\l,N}$ on ${\Bbb R}^4$ and a constant $K$ independent of $\l$, $N$,
such that $\b_{\l,N}(x)=1$ if $|x|\ge \half\l^{1/2}$ and $\b_{\l,N}(x)=0$
if $|x|\le N^{-1}\l^{1/2}$, and $\|d\b_{\l,N}\|_{L^4({\Bbb R}^4,\d)}\le
K(\log N)^{-3/4}$.} 

Define $C^\8$ cut-off functions $\be_I$ on each $X_I$ by setting
$$
\be_I\equiv(\phi_{Is}^{-1})^*\be_{\la_I,N}
\prod_{I_+}(\phi_{I_+}^{-1})^*\be_{\la_{I_+},N}
\quad\hbox{on }X_I, \eqlabel\eqnlabel{\eqBetaI}
$$ 
where the factor $(\phi_{Is}^{-1})^*\be_{\la_I,N}$ is omitted when $I=0$. Here, the cut-off functions comprising $\be_I$ have been extended so that $\be_I=1$ on the complement in $X_I$ of the balls $B_{Is}(\half\la_I^{1/2})$ and $B_{I+}(\half\la_{I_+}^{1/2})$.
Also, $\be_I=0$ on the balls $B_{Is}(N^{-1}\la_I^{1/2})$ and $B_{I+}(N^{-1}\la_{I_+}^{1/2})$ in $X_I$: thus, we may extend $\be_I$ by zero to give $\be_I\in C^\8(X)$. The $L^4$ estimate of Lemma \ref{\lemBetaEst} implies that
$$
\|d\be_I\|_{L^4(X_I,g_I)}\le cK(\log N)^{-3/4}. 
\eqlabel\eqnlabel{\eqDBetaILogNEst}
$$
for some $c=c(g_0,k)$.
For the cut-off functions $\{\ga_I\}$ defined by Eqs. \eqref{\eqGamma} and \eqref{\eqGammaI}, we recall that $\sum_I \c_I=1$ on $X$. Note also that $\b_I=1$ on the support of $\c_I$.

Define operators $Q_I:L^p_1\W^{+,g_I}(X_I,\ad P_I)\to L^p\W^1(X_I,\ad P_I)$
by setting $Q_I=\b_I P_I\c_I$
Define a right parametrix $Q:\sL^p_1\W^{+,g}(X,\ad P)\to L^p\W^1(X,\ad P)$ for the
operator $d_{A'}^{+,g}$ by $Q=\sum_I Q_I$. The 
error operator $R:L^p\W^{+,g}(X,\ad P)\to L^p\W^{+,g}(X,\ad P)$ is then given by
$$
d^{+,g}_{A'}Q=1+R.\eqlabel
$$
Lemmas \ref{\lemMetComparEst} and \ref{\lemPIEst} 
then yield the following estimates for the operators $Q_I$ and $Q$.
 
\proclaim{Lemma}{\thmlabel{\lemQEst}
There are constants 
$C_i=C_i(g_0,p,\sT)$, $i=1,2$, such that for any $t\in\sT$ and any
$\xi\in L^p\Om^{+,g_I}(X_I,\ad P_I)$, for {\rm (a)}, or any $\xi\in L^p\Om^{+,g}(X,\ad P)$, for {\rm (b)}, the following bounds hold.
\item{\rm (a)}
$\|Q_I\xi\|_{L^q(X_I,g_I)}\le C_1\|Q_I\xi\|_{L^p_1(X_I,A_I,g_I)}\le
C_2\|\xi\|_{L^p(X_I,g_I)}$,  
\item{\rm (b)}
$\|Q\xi\|_{L^q(X,g)}\le C_1\|Q\xi\|_{\sL^p_1(X)}\le C_2\|\xi\|_{L^p(X,g)}$.}

Next, there is an analogue of Lemma 7.2.14 [D-K] (see also [D-K, p. 294]),
giving an $L^p$ bound for the operator $R$. The proof follows easily from 
Lemmas \ref{\lemAIprimeEst}, \ref{\lemMetXIEst}, and \ref{\lemPIEst}, and
Eq. \eqref{\eqDBetaILogNEst}. In [D-K] it is assumed that the metrics $g_I$ are flat in small
neighbourhoods of the points $x_I$, but this restriction is easily removed
using Lemma \ref{\lemMetXIEst}.

\proclaim{Lemma}{\thmlabel{\lemREst}
{\rm } There is a constant
$\e=\e(\ol b,N,p)$, with $\e\to 0$ as $N\to\8$ and $\ol b\to 0$
such that for any $t\in \sT$ and $\xi\in L^p\W^{+,g}(X,\ad P)$,
$\|R\xi\|_{L^p(X,g)}\le \e\|\xi\|_{L^p(X,g)}$.} 

Thus, {\it for the remainder of this article}
choose $N_0>4$ large enough and $b_0\le 1$ small enough so that 
$\e(\ol b,N,p)\le 2/3$ for all $\ol b\le b_0$ and $N\ge N_0$,
and fix $N=N_0$ and $b_I=4N\l_I^{1/2}$ for all $I\in\sI$.
We now construct a right inverse $P$ for $d^{+,g}_{A'}$.
Lemma \ref{\lemREst} gives the
$(L^p, L^p)$ operator norm bounds $\|R\|\le 2/3$ and 
$\|(1+R)^{-1}\|\le 3$. 
Since $Q_I=\b_IP_I\c_I$, we have the $(L^p, L^q)$ operator norm
bound $\|Q_I\|\le C_I$, say, giving the $(L^p, L^q)$ operator norm bound
$\|Q\|\le C\equiv\sum_I C_I$. In summary, there is the following version of
Proposition 7.2.35 [D-K]. 

\proclaim{Proposition}{\thmlabel{\propPEst}
There are constants $N_0$ and $b_0$ such that for any $N\ge N_0$, $\ol
b\le b_0$, and $t\in\sT$, the operator 
$P\equiv Q(1+R)^{-1}:\sL^p_1\W^{+,g}(X,\ad P)\to L^p\W^1(X,\ad P)$ is a right
inverse to $d^{+,g}_{A'}$ and there are
constants $C_i=C_i(g_0,p,\sT)$, $i=1,2$ such that for any
$\xi\in L^p\W^{+,g}(X,\ad P)$, 
$$
\|P\xi\|_{L^q(X,g)}\le C_1\|P\xi\|_{\sL^p_1(X)}\le C_2\|\xi\|_{L^p(X,g)}. 
$$}

We next construct families of solutions to the full non-linear
anti-self-dual equation 
over connected sums. For each $t\in\sT$
we seek a solution $A(t)=A'(t)+a(t)$ to 
$F^{+,g}(A'+a)=0$, or equivalently
$$
d^{+,g}_{A'}a+(a\wedge a)^{+,g}=-F^{+,g}(A'),\eqlabel\eqnlabel{\eqASDa}
$$
where $a\in\W^1(X,\ad P)$. 
If $a=P\xi$, with 
$\xi(t)\in\W^{+,g}(X,\ad P)$, then this equation becomes
$$
\xi+(P\xi\wedge P\xi)^{+,g}=-F^{+,g}(A').
\eqlabel\eqnlabel{\eqASDXi}
$$
With the aid of Lemma 7.2.23 [D-K, p. 290] (an application of the Contraction Mapping Theorem to Eq. \eqref{\eqASDXi}) and Proposition \ref{\propPEst}, one easily obtains the version below of Theorem 7.2.24 [D-K].

\proclaim{Theorem}{\thmlabel{\thmDKExist} For sufficiently small $\la_0< 1$, sufficiently large $N_0>4$, and sufficiently small $T_{A_I}$, $I\in\sI$, the following holds. For any $t\in\sT$, there exists an $L^p_1$ $g$-anti-self-dual connection
$A(t)=A'(t)+a(t)$ over $X$, with $a(t)=P\xi(t)$. There are positive constants $C_i=C_i(g_0,p,\sT)$, $i=1,2,3$, such that 
$$\leqalignno
{\|a\|_{L^q(X,g)}&\le C_1\|\xi\|_{L^p(X,g)}\le C_2\|F^{+,g}(A')\|_{L^p(X,g)}\le C_3{\ol b}^{4/p}. &\cr}
$$} 

We pull back the $g$-anti-self-dual connections $A$ on $P\to X$ via the conformal maps $f_I$ to
give $g_0$-anti-self-dual connections $\hat A=\hat A'+\hat a$ on $\hat P\to X_0$, where $\hat A$ is defined by 
$$
\hat A=f_0^*\cdots f_I^*A\quad\hbox{over }f_0^{-1}\cdots f_I^{-1}(X_I'),
\eqlabel
$$
and similarly for $\hat A'$ and $\hat a$. In particular, $\hat A=\hat A'+\hat a$ is a solution to the $g_0$-anti-self-dual equation 
$F^{+,g_0}(\hat A'+\hat a)=0$ over $X_0$, or explicitly
$$
d^{+,g}_{\hat A'}\hat a+(\hat a\wedge \hat a)^{+,g_0}
=-F^{+,g_0}(\hat A'),\eqlabel\eqnlabel{\eqASDhata}
$$
where $\hat a\in\W^1(X_0,\ad \hat P)$. Standard arguments show that
the anti-self-dual connections
$A$ and $\hat A$ are actually $C^\8$ and that 
they are smooth points of the moduli spaces $M_{X,k}(g)$ and $M_{X_0,k}(g_0)$ [D-K]:

\proclaim{Lemma}{\thmlabel{\lemAIrred}
Let $A$ be the $g$-anti-self-dual connection over $X$ produced by\/ {\rm
Theorem \ref{\thmDKExist}} and let $\hat A$ be the corresponding
$g_0$-anti-self-dual connection over $X_0$. Then the following hold:
\item{\rm (a)} The connections $A$ and $\hat A$ are $C^\8$, 
\item{\rm (b)} $H^0_A=0$ and $H^0_{\hat A}=0$, for small enough $b_0$ and large enough $N_0$, 
\item{\rm (c)} $H_A^2=0$ and $H_{\hat A}^2=0$.}

From \S\ref{\secDqStrongConv}, we recall that $D_q$ is the metric $\cB_{X,k}$
given by 
$D_q([A],[B])=\inf_{u\in\cG}\|A-u(B)\|_{L_q(X,g)}$.
In particular, we have the following version of Theorem 7.2.62
[D-K] (compare also Theorem 4.53 [D86]). 

\proclaim{Theorem}{\thmlabel{\thmDKMap}
Let $A_I$ be $g_I$-anti-self-dual connections on $G$ bundles $P_I$ over manifolds $X_I$, $I\in\sI$. If $I=0$, then $X_0$ is a closed, oriented, $C^\8$ four-manifold with generic $C^\8$ metric $g_0$ and negative definite intersection form. If $I>0$, then $X_I=\SS^4$ with standard round metric $g_1$ of radius 1. Let $X=\#_{I\in\sI}X_I$, the connected sum four-manifold with $C^\8$ metric $g$ (conformally equivalent to $g_0$) determined by the choice of points $\{x_I\}$, frames $\{v_I\}$, scales $\{\la_I\}$, and neck width parameter $N$. Let $P$ be the connected sum bundle over $X$, where $c_2(P)=k\ge 1$. Let $\ol\la=\max_{I\in\sI}\la_I$. Let $T_{A_I}$ be open balls centred at $0\in H^1_{A_I}$, $I\in\sI$, let 
$\Ga=\prod_{I\in\script{I}}\Ga_{A_I}$ and
$T=T_{A_0}\times\prod_{I\in\script{I}}(T_{A_I}\times\Gl_{x_I})$,
as in\/ {\rm Eqs. \eqref{\eqT}} and\/ {\rm \eqref{\eqBigStab}}.
Then, for sufficiently small $\la_0< 1$, sufficiently large $N_0>4$, and sufficiently small $T_{A_I}$, $I\in\sI$, the following holds. There is a $C^\8$ homeomorphism onto an open subset:
$$
\sJ:T/\Ga\too U\subset M_{X,P}^*(g), \qquad t\longmapsto [A(t)],
$$
where $A(t)=A'(t)+a(t)$, $a(t)=P\xi(t)$, and $\xi(t)$ are as in 
\/ {\rm Theorem \ref{\thmDKExist}}.
For any $\nu>0$ and $4\le q<\8$, the manifold $T$ and constant $\la_0(\nu)$ can be chosen so that, for all $\ol\la<\la_0(\nu)$,
$U=\{[A]\in M_{X,P}^*(g): D_q([A|_{X_I''}],[A_I])<\nu\}$,
for all $I\in\sI$.}

\pf This is a straightforward generalisation of Theorem 7.2.62 [D-K] to the
case of multiple connected sums (see [D-K, \S7.2.8]) and a restriction to
the case where $G=\SU(2)$ and $b^+(X_0)=0$. The
metric $g_0$ is not required to be flat in small neighourhoods of the
gluing sites $x_I\in X_0$. 
Lemma \ref{\lemAIrred} implies that the image of $\sJ$ lies in the dense open subset $M_{X,P}^*(g)\subset M_{X,P}(g)$. The fact that $\sJ$ is $C^\8$ is a calculation of the type that appears many times in \S\S5.3, 5.4, and 5.5. See also Appendix A [T84b] and Remark 4.24 [D86]. \qed

We refer to $\script{J}$ as a {\it gluing map} over the connected sum and
its image $U\subset M_{X,k}^*(g)$ as a {\it gluing neighbourhood}.
Moreover, $\sJ$  extends to a $C^\8$ gluing map on the larger parameter
spaces $\sT$ and $\sT^0$ 
of Eqs. \eqref{\eqBigT} and \eqref{\eqBigTZero}. Further properties of
these maps are described in the next section. 
Lastly, for the original metric $g_0$ on
the base four-manifold $X_0$, Theorem \ref{\thmDKMap} takes the following
form. 

\proclaim{Corollary}{\thmlabel{\corDKHatMap} Given the hypotheses of\/ {\rm Theorem \ref{\thmDKMap}}, there is a homeomorphism onto an open subset
$$
\hat\sJ:T/\Ga\too V\subset M_{X_0,\hat P}^*(g_0), 
\qquad t\longmapsto [\hat A(t)],
$$
where $V\subset M_{X_0,\hat P}^*(g_0)$ is obtained by pulling back the
subset $U\subset M_{X,P}^*(g)$ of\/ {\rm Theorem \ref{\thmDKMap}}. 
\qed}

Again, $\hat\sJ$ extends to a $C^\8$ map on the larger parameter spaces
$\sT$ and $\sT^0$ and additional properties of $\hat\sJ$ are discussed in the
next section. 


\section{Structure of the compactified moduli spaces}
The bubbling ends of $\ol{M}^u_{X_0,k}(g_0)$ away from the diagonals are
described in [D-K, \S 8.2]. We extend this description to neighbourhoods of
points in the diagonals of the Uhlenbeck compactification. For related
constructions and some further details, we refer to the papers of Taubes
and Donaldson.
 
The proposition below is
the basic result we require in order to parametrise neighbourhoods covering
the ends of $M_{X,k}^*(g)$ away from the reducible connections. See also
[D83, \S III], [D86, \S IV], and [T84b, p. 529] for various special
cases of the statements below. The proof below is similar to arguments in
the proofs of Proposition 5.17 \& 5.22 [D-S] and Theorem 4.53 [D86, p. 316
\& p. 325]. 

\proclaim{Proposition}{\thmlabel{\propGluMapDiff}
Given the hypotheses of\/ {\rm Theorem \ref{\thmDKMap}}, the following hold:
\item{\rm (a)} The approximate gluing map $\script{J}':T/\Ga \to \cB_{X,k}^*$ is a $C^\8$ embedding, 
\item{\rm (b)} The gluing map $\script{J}:T/\Ga\to U\subset M_{X,k}^*(g)$ is a diffeomorphism onto an open subset,
\item{\rm (c)} The extended gluing map $\script{J}:\sT/\Ga\to\sU\subset M_{X,k}^*(g)$ is a $C^\8$ submersion onto an open subset,
\item{\rm (d)} The extended gluing map $\script{J}:\sT^0/\Ga\to\sU^0\subset M_{X,k}^*(g)$ is a diffeomorphism onto an open subset.}

\pf (a) The proof is essentially the same as the argument required for (b) and so is omitted.
(b) From Theorem \ref{\thmDKMap}, $\sJ$ is a $C^\8$ homeomorphism and so it is enough to show that $\sJ$ is also an immersion, since $T/\Ga$ has dimension equal to that of $M_{X,k}^*(g)$.
From the proof of Theorem \ref{\thmDKMap}, there is a $C^\8$ $\Ga$-equivariant gluing map $\tilde\script{J}:T\to\cA_{X,k}^*$, $t\mapsto A(t)$. So, we first show that $\tilde\script{J}$ is an immersion and then conclude that the induced map on quotients is a diffeomorphism. The constant $\la_0$ may be chosen as small as desired and in (a) and (b),the $\la_I$ and $x_I$ may be held fixed.
 
\step{1.} {\it Definition of restriction maps.}
Choose cut-off functions $\psi_I$, as in \S\ref{\secGluCon}, which are zero on the balls 
$B_{Is}(b_I/2)$, $B_{I_+}(b_{I_+}/2)$ and equal to $1$ on the complement in
$X_I$ of the slightly larger balls $B_{Is}(b_I)$, $B_{I_+}(b_{I_+})$.
Define a map
$\pi_{X_I}: L^2\Om^1(X,\ad P)\to L^2\Om^1(X_I,\ad P_I)$
by left multiplication with $\psi_I$, so that
$$
\|\om-\pi_{X_I}\om\|_{L^2(X_I,g_I)}=O(\ol\la),
\qquad\om\in\Om^1(X_I,g_I),\eqlabel\eqnlabel{\eqPiIEst}
$$
since $\psi_I$ is equal to $1$ on the complement of a set in $X_I$ of $g_I$-volume $O(\ol\la^2)$. 
Next, for $I>0$, choose a cut-off function which is zero outside the annulus $\Om_{Is}=\Om(x_{Is},N^{-1}\la_I^{1/2},N\la_I^{1/2})$ in $X_I$ and which is equal to $1$ on the slightly smaller annulus 
$\Om(x_{Is},\half\la_I^{1/2},2\la_I^{1/2})$ containing the supports of the derivatives of the cut-off functions $\ga_{I_-},\ga_I$. Define a map 
$\pi_{\Om_I}: L^2\Om^1(X,\ad P)\to L^2\Om^1(\Om_{Is},\ad P_I)$
by left multiplication with this cut-off function. 
Lastly, let $\Pi=\pi_0\oplus_{I>0}(\pi_{X_I}\oplus\pi_{\Om_I})$ be the
induced map 
$$
L^2\Om^1(X,\ad P)\too L^2\Om^1(X_0,\ad P_0)\bigoplus_{I>0}
\left(L^2\Om^1(X_I,\ad P_I)\oplus L^2\Om^1(\Om_{Is},\ad P_I)\right).
$$

\step{2.} {\it Partial derivatives with respect to lower moduli parameters.}
We have $C^\8$ $\Ga_{A_I}$-equivariant maps $\tilde\vth_I:T_{A_I}\to\cA_{X_I,P_I}^*$, $t_I\mapsto A_I(t_I)$ given by the Kuranishi model.
Let $v$ be a tangent vector to $T_{A_I}$, i.e., suppose $[v]\in H_{A_I}^1$. Then Eq. \eqref{\eqPiIEst} and the estimates of \S 5.4 give the following bounds for the differentials with respect to the lower moduli parameters:
$$\leqalignno
{\|\pi_{X_I}D\script{J}(v)-D\tilde\vth_I(v)\|_{L^2(X_I,g_I)}
&=O(\ol\la^{1/2}). &\numeqn\eqnlabel{\eqDpiJEst}\cr}
$$
The map $\tilde\vth_I$ is an immersion and so the range of $D\tilde\vth_I$ has dimension equal to $\dim H^1_{A_I}$. For small enough $\ol\la$, Eq. \eqref{\eqDpiJEst} implies that the range of $\pi_{X_I}D\script{J}$ also has dimension equal to $\dim H^1_{A_I}$.

\step{3.} {\it Partial derivatives with respect to gluing parameters.}
Let $v$ be a tangent vector to $\Gl_I$. The estimates of \S 5.5 give the following bounds for the differentials with respect to the gluing parameters:
$$\leqalignno
{\|\pi_{\Om_I}D\script{J}(v)-D\script{J}'(v)\|_{L^4(X_I,g_I)}
&=O(\ol\la^2),  &\numeqn\eqnlabel{\eqvpiDJEst}\cr}
$$
recalling that $D\script{J}'(v)$ is supported on 
$\Om_I(\half\la_I^{1/2},2\la_I^{1/2})$.
But from Proposition \ref{\propDAprimeDvEst} we have 
$$
\|D\script{J}'(v)\|_{L^4(X_I,g_I)}\ge c|v|,
\eqlabel\eqnlabel{\eqvpiDJlowBnd}
$$ 
for some constant $c>0$ independent of $\ol\la$. In particular, the range of $\pi_{\Om_I}D\script{J}'$ has dimension equal to $\dim\Gl_I$. So, for sufficiently small $\ol\la$, Eqs. \eqref{\eqvpiDJEst} and \eqref{\eqvpiDJlowBnd} imply that the range of $\pi_{\Om_I}D\script{J}$ also has dimension equal to $\dim \Gl_I$. 

\step{4.} {\it The quotient map.} Combining these observations, we find that the range of $\Pi D\script{J}$ has dimension equal to 
$\dim H_{A_0}+\sum_{I>0}(\dim H_{A_I}^1+\dim\Gl_I)=\dim T$,
so that $\ker \Pi D\script{J}=0$ and $\script{J}$ is an immersion. 
From Theorem \ref{\thmDKMap}, the open subset $\tilde U\equiv\tilde\script{J}(T)$ in $\cA_{X,k}^*$ projects to an open subset $U\equiv\script{J}(T)$ in $M_{X,k}^*(g)$ and composing $\tilde\script{J}$ with the projection $\cA_{X,k}^*\to\cA_{X,k}^*/\cG$, we obtain a submersion $\script{J}:T\to M_{X,k}^*(g)$. The group $\Ga$ acts freely on $T$, $\tilde\script{J}$ is $\Ga$-equivariant, and $\dim T/\Ga=\dim M_{X,k}^*(g)$, and the gluing map descends to a diffeomorphism $\script{J}:T/\Ga\to M_{X,k}^*(g)$, as required. 
(c) This follows from (b). For the derivatives with respect to $\la_I$ or
$x_I$, the cut-off functions required to define $\Pi$
should be replaced by cut-offs with similar supports and which are {\it
fixed} with respect to small variations in the scales and centres.  

\noindent (d) This is similar to the proof of (c) and uses Proposition \ref{\propCentModDiff}.\qed

In order to parametrise neighbourhoods of boundary points in $\ol{M}_{X_0,k}^u(g_0)$, we use the following corollary to Proposition \ref{\propGluMapDiff}.

\proclaim{Corollary}{\thmlabel{\corGluMapHatDiff}
Given the hypotheses of\/ {\rm Theorem \ref{\thmDKMap}}, 
the following hold:
\item{\rm (a)} The approximate gluing map $\hat\script{J}':T/\Ga \to \cB_{X_0,k}^*$ is a $C^\8$ embedding,
\item{\rm (b)} The gluing map $\hat\script{J}:T/\Ga\to V\subset M_{X_0,k}^*(g_0)$ is a diffeomorphism onto an open subset,
\item{\rm (c)} The extended gluing map $\hat\script{J}:\sT/\Ga\to\sV\subset M_{X_0,k}^*(g_0)$ is a $C^\8$ submersion onto an open subset,
\item{\rm (d)} The extended gluing map $\hat\script{J}:\sT^0/\Ga\to\sV^0\subset M_{X_0,k}^*(g_0)$ is a diffeomorphism onto an open subset.}

Taken together,
Theorems 7.3.2 and 7.2.62 in [D-K] imply that if $A$ is any $g$-anti-self-dual
connection on a fixed $G$ bundle $P$ over the connected sum $X$ and the
necks $\Om$ are all sufficiently pinched (so that $\ol\la$ is small), then
$[A]$ lies in the image of the gluing map. The corresponding statement in
our application is given below. 

\proclaim{Theorem}{\thmlabel{\thmGluMapOnto} Given the hypotheses of\/ {\rm Theorem \ref{\thmDKMap}}, then the following holds. Let $\{A_\al\}_{\a=1}^\8$ be a sequence of connections on a $G$ bundle $P$ over the connected sum $X=\#_{I\in\sI}X_I$ which are anti-self-dual with respect to the sequence of metrics $\{g_\al\}_{\a=1}^\8$ determined by the sequences of scales $\{\la_{I\al}\}$ with $\ol\l_\a\to 0$, a fixed neck width parameter $N$, sequences of points $\{x_{I\al}\}$ converging to $\{x_I\}$, and frames in $FX_0|_{x_{I\al}}$ converging to frames in $FX_0|_{x_I}$. Suppose the sequence $\{A_\al\}_{\a=1}^\8$ is strongly convergent to $(A_I)_{I\in\sI}$, where $A_I$ is a $g_I$-anti-self-dual connection over each summand $X_I$.
For $\al_0$ sufficiently large, there exists a gluing neighbourhood $\sU$ such that $[A_\al]\in\sU$, for all $\al\ge\al_0$.}

\pf See [D-K, \S7.3.1]. Theorem \ref{\thmStrongLqConv} implies that the
sequence $\{A_\al\}$ is $D_q$ convergent (for any $4\le q<\8$) to
$(A_I)_{I\in\sI}$. So, Theorem \ref{\thmDKMap} implies that the points
$[A_\al]$ are contained in a gluing neighbourhood $\sU$, for all
$\al\ge\al_0$ if $\al_0$ is sufficiently large. \qed 

Recall that $\Gl_{x_I}=\SU(2)\simeq \SS^3$, a copy of the standard three-sphere, and let $\ol{\Gl}_{x_I}$ be the closure of $\Gl_{x_I}\times (0,\la_0)$ in 
the cone $(\Gl_{x_I}\times [0,\la_0))/\sim$, where $(\rho,0)\sim (\rho',0)$ if $\rho,\rho'\in\Gl_{x_I}$. 
Then, by analogy with [D-K, \S 8.2] and [D86, \S V], we set
$$\leqalignno
{\ol\script{T}
&\equiv T_{A_0}\times\prod_{I\in\script{I}}
\left(T_{A_I}\times B(x_I,r_0)\times\ol{\Gl}_{x_I})\right), 
&\numeqn\eqnlabel{\eqClosedBigT}\cr}
$$
and likewise, define $\ol{\sT}^0$. It is also convenient to define
$$
\rd\script{T}
\equiv \{t_\8=(t_I,y_I,\rho_I,\la_I)_{I\in\sI}\in \ol\script{T}: \la_I=0\hbox{ for some }I\},
\eqlabel\eqnlabel{\eqBdryBigT}
$$
where the 4-tuple $(t_I,y_I,\rho_I,\la_I)$ above is replaced by $t^0$, if $I=0$.
The space $\rd\script{T}^0$ is defined similarly.
Moreover, the gluing map $\sJ$ has a natural definition on the boundary $\rd\script{T}$. 
Suppose $t_\8\in\rd\sT$ and let $(\la_1,\dots,\la_c)$ denote the corresponding scales in Eq. \eqref{\eqBdryBigT} which have been set equal to zero.
By cutting the edges with $\la_i=0$, we may view the tree $\sI$ as a union of subtrees $\cup_{i=1}^c\sI^i$. If $t_\8\in\rd\sT$, we write $t_\8=(t^1,\dots,t^c)$, with $t^i\in\sT^i$,   
and set 
$$
\sJ(t_\8)\equiv (\sJ(t^1),\dots,\sJ(t^c)), \qquad t_\8\in\sT,\eqlabel
$$ 
where each $\sJ(t^i)$ is an
anti-self-dual connection over a connected sum $Y_i=\#_{I\in\sI^i}X_I$, say, and $X=\#_{I\in\sI}X_I=\#_{i=1}^c Y_i$. The relationship between the gluing maps $\sJ$ and $\sJ^i$ is explained by the continuity result below, which we just state in the special case $X=X_0\#X_1\#X_2$, for the sake of clarity. The argument required for this case carries over with no significant change to the more general cases just described.

\proclaim{Proposition}{\thmlabel{\propGluMapMultConn} Let $X=X_0\#X_1\#X_2$, let $Y=X_0\#X_1$, and let $Y''=Y\less B(x_1,\half\la_1^{1/2})$.
Assume the hypotheses of\/ {\rm Theorem \ref{\thmDKMap}} and let $\sJ_X$, $\sJ_Y$ be the gluing maps over the connected sums $X$ and $Y$, respectively. Then there is an $\eps=\eps(q)>0$ and a constant $C=C(g_0,q,\sT)$ such that
$\|\sJ_X(t)|_{Y''}-\sJ_Y\|_{L^q(Y,g)} \le C\la_1^{\eps_0}$.}

The proof is similar to that of Proposition 7.2.64 [D-K] and the arguments
in \S5.3, and so is omitted. 
It now follows that $\sJ$ extends continuously to $\ol\script{T}$.

\proclaim{Proposition}{\thmlabel{\propGluMapExtend}
Assume the hypotheses of\/ {\rm Theorem \ref{\thmDKMap}}.  Let 
$\{t_\al\}_{\al=1}^\8$ be a sequence in $\script{T}$ which converges to $t_\8\in\rd\sT$. Then the sequence $\{\sJ(t_\al)\}_{\al=1}^\8$ converges strongly to $\sJ(t_\8)$.}

\pf Let $\{\la_i\}_{i=1}^c$ denote the scales, determined by $t_\8$, which have been set equal to zero in Eq. \eqref{\eqBdryBigT}. The points $t_\al\in\sT$ are then naturally written as $t_\al=(t_\al^1,\dots,t_\al^c)$, with the sequences $t_\al^i$ converging to $t^i\in\sT^i$, say.
According to Proposition \ref{\propGluMapMultConn}, the sequence $\sJ(t_\al)$ is then $D_q$ convergent to $(\sJ(t^1),\dots,\sJ(t^c))$ and hence, strongly convergent by Theorem \ref{\thmStrongLqConv}. \qed 

It remains to show
that $M_{X_0,k}(g_0)$ has a finite cover consisting of gluing
neighbourhoods. Of course, away from the bubbling ends, the moduli space is
covered by the standard Kuranishi charts. In addition, the geometry of
these charts around the reducible connections has already been analysed in
[G-P89], so our focus here is on the bubbling ends.   
Given any Uhlenbeck boundary point 
$(A_0,x_1,\dots,x_l)\in {\ol M}^u_{X_0,k}(g_0)$, where $c_2(A_0)=k-l$ and each $x_i$ has multiplicity 1, Theorem 8.2.3 [D-K] provides an open neighbourhood $\ol\sV$ of $(A_0,x_1,\dots,x_l)$ in ${\ol M}^u_{X_0,k}(g_0)$, a parameter space $\sT^0/\Ga$, and a gluing map $\hat\sJ$ giving a homeomorphism of $\sT^0/\Ga$ with $\sV=\ol\sV\cap M_{X_0,k}^*(g_0)$. Theorem 8.2.4 in [D-K] states that this gluing map extends to a homeomorphism $\hat\sJ:\ol\sT^0/\Ga\to\ol\sV$. Thus, away from the diagonals, the ends of ${\ol M}^u_{X_0,k}(g_0)$ are covered by gluing neighbourhoods. 
The generalisations below provide a covering of the ends of ${\ol M}^u_{X_0,k}(g_0)$, including the diagonal boundary points.

\proclaim{Theorem}{\thmlabel{\thmGluMapUhlNbd}
Let $(A_0,x_1,\dots,x_{m_0})$ be a boundary point in ${\ol M}^u_{X_0,k}(g_0)$. Under the hypotheses of\/ {\rm Theorem \ref{\thmDKMap}}, there exist neighbourhoods $\ol\sV\subset{\ol M}_{X_0,k}^u(g_0)$ of $(A_0,x_1,\dots,x_{m_0})$ and a parameter space $\sT^0$ such that the following holds. If $\sV=\ol\sV\cap M_{X_0,k}^*(g_0)$, then
the gluing map $\hat\script{J}:\sT^0/\Ga\to \sV$ is a diffeomorphism.}

\pf Suppose $\{[A_\al]\}_{\al=1}^\8$ is a sequence in $M_{X_0,k}(g_0)$, converging weakly to the Uhlenbeck limit $(A_0,x_1,\dots,x_{m_0})$. Let $\{[\check A_\al]\}_{\al=1}^\8$ be the corresponding strongly convergent sequence in $M_{X,k}(g_\al)$ with the bubble tree limit $(A_I,x_I)_{I\in\sI}$. Then Theorem \ref{\thmGluMapOnto} produces a gluing neighbourhood $\sJ(\sT^0/\Ga)=\sU\subset M_{X,k}(g_\al)$ and an $\al_0$ such that $[\check A_\al]\in\sU$ for all $\al\ge \al_0$. Let $\sV$ be the corresponding neighbourhood in $M_{X_0,k}(g_0)$. Then the conclusions follow from Corollary \ref{\corGluMapHatDiff}.\qed

\proclaim{Theorem}{\thmlabel{\thmGluMapUhlNbdExt} Given the hypotheses of\/ 
{\rm Theorem \ref{\thmGluMapUhlNbd}}, the gluing map $\hat\sJ$ extends to a homeomorphism of $\ol{\sT}^0/\Ga$ with a neighbourhood $\ol\sV$ of 
$(A_0,x_1,\dots,x_{m_0})$ in ${\ol M}^u_{X_0,k}(g_0)$.}

\pf This follows from Proposition \ref{\propGluMapExtend} and Theorem \ref{\thmGluMapUhlNbdExt}. \qed

\rmk So, every boundary point in ${\ol M}^u_{X_0,k}(g_0)$ has a neighbourhood constructible by gluing. Plainly, the same statement holds for boundary points in ${\ol M}^\tau_{X_0,k}(g_0)$. 


\section{Derivatives with respect to scales and centres}
\seclabel{\secDHatADLaq} 
The main purpose of this section is to obtain $L^2$ estimates for the partial derivatives of the family of $g_0$-anti-self-dual connections $\hat A$ with respect to the scales $\l_I$ and centres $x_I$. 

Unless noted otherwise, throughout this section and for the remainder of this article, we assume that $p$ and $q$ are Sobolev exponents satisfying the strict inequalities $2<p<4$ and $4<q<\8$, where $q$ is determined by $1/p=1/4+1/q$. The constant $\la_0>0$ is assumed small and may be decreased as needed. We use $C=C(g_0,p,\sT)$ to denote constants which are independent of the points $t=(t_I,\rho_I,x_I,\la_I)\in\sT$. As usual, we abbreviate the derivative with respect to the centre parameters, $p_I^\mu\rd/\rd q_I^\mu$ (where $|p_I|\le 1$) by $\rd/\rd p_I$.

Denoting $\eta\equiv -F^{+g}(A')$ in Eq. \eqref{\eqASDXi}, we have the following preliminary estimate for the derivatives of $\hat a$ with respect to the $\l_I$ and $x_I$ parameters.

\proclaim{Lemma}{\thmlabel{\lemDaDLaqPreEst}
Let $\xi$ and $a=P\xi$ be as in \/ {\rm Theorem \ref{\thmDKMap}}, and assume that the conditions of that theorem hold. Then, for small enough $\la_0>0$, 
there is a constant $C=C(g_0,\sT)$ such that for any $t\in\sT$,
\item{\rm (a)}
$\|\rd{\hat a}/\rd{\l_I}\|_{L^2(X_0,g_0)}
\le C\left(\|\pd{P}{\l_I}\xi\|_{L^2(X,g)}
+\|\rd{\xi}/\rd{\l_I}\|_{L^2(X,g)}+\ol\la\l_I^{-1/2}\right)$,
\item{\rm (b)}
$\|\rd{\hat a}/\rd{p_I}\|_{L^2(X_0,g_0)}
\le C\left(\|\pd{P}{p_I}\xi\|_{L^2(X,g)}
+\|\rd{\xi}/\rd{p_I}\|_{L^2(X,g)}+\ol\la\right)$.}

\pf From Proposition \ref{\propDhatOmDLaq}, we have
$$
\left\|\pd{\hat a}{\l_I}\right\|_{L^2(X_0,g_0)}\le 
C\left\|\pd{a}{\l_I}\right\|_{L^2(X,g)}
+C\la_I^{-1/2}\|a\|_{\sL^2_1(X)},
$$
where $a=P\xi$ and 
$\pd{P\xi}{\l_I}=\pd{P}{\l_I}\xi+P\pd{\xi}{\l_I}$.
The estimates of Proposition \ref{\propPEst} and
Theorem \ref{\thmDKExist} then gives (a). The proof of (b) is similar.
\qed

We now differentiate the $g$-anti-self-dual equation and obtain a priori estimates for the partial derivatives of $\xi$ with respect to $\l_I$ and $x_I$.

\proclaim{Lemma}{\thmlabel{\lemDXiDLaqPreEst} 
Let $\xi$ be as in \/ {\rm Theorem \ref{\thmDKMap}} and assume that the conditions of that theorem hold. Then, for small enough $\la_0>0$,
there is a constant $C=C(g_0,\sT)$ such that for any $t\in\sT$,
\item{\rm (a)}
$\|\rd{\xi}/\rd{\l_I}\|_{L^2(X,g)}\le 
C\left(1+\ol{\la}^2\la_I^{-1/2}+\ol\la\|\pd{P}{\l_I}\xi\|_{L^4(X,g)}\right)$,
\item{\rm (b)}
$\|\rd{\xi}/\rd{p_I}\|_{L^2(X,g)}\le 
C\left(1+\ol\la\|\pd{P}{p_I}\xi\|_{L^4(X,g)}\right)$.}

\pf  Differentiating Eq. \eqref{\eqASDXi} with respect to
$\l_I$ gives
$$\leqalignno
{\pd{\xi}{\l_I}&=\pd{\eta}{\l_I}-\pd{*_g}{\l_I}(P\xi\wedge P\xi) -\left(\pd{P\xi}{\l_I}\wedge P\xi\right)^{+,g}
-\left(P\xi\wedge \pd{P\xi}{\l_I}\right)^{+,g}. &\cr}
$$
The estimates of Lemma \ref{\lemDHodgeStarDLaq} 
and Proposition \ref{\propPEst} imply that 
$$\leqalignno
{\left\|\pd{\xi}{\l_I}\right\|_{L^2} 
&\le \left\|\pd{\eta}{\l_I}\right\|_{L^2}+C\|\xi\|_{L^2}^2\la_I^{-1/2}
+C\|\xi\|_{L^2}\left(\left\|\pd{\xi}{\l_I}\right\|_{L^2}
+\left\|\pd{P}{\l_I}\xi\right\|_{L^4}\right). &\cr}
$$
Proposition \ref{\propDFAprimeDLaq} gives $\|\rd\eta/\rd\l_I\|_{L^2}\le C$, while Theorem \ref{\thmDKExist} gives
$\|\xi\|_{L^2}\le C\ol\la$. Thus, for $\la_0$ small enough, we may assume  $C\|\xi\|_{L^2}\le 1/2$. 
Part (a) then follows by combining the above estimates and rearranging,
and the proof of (b) is similar.
\qed

To complete our task, we need an estimate for the derivatives of $P$ with respect to $\la_I$ and $x_I$. Before proceeding, we first record some bounds for the derivatives of the cut-off functions $\be_I$ and $\ga_I$.
Suppose $1\le p<\8$. From the definition of $\be_I$
there is a constant $C=C(g_I,N,p)$ such that  
$$
|d\b_I|_{g_I}\le C\l_I^{-1/2}\quad\hbox{on }\Om_I, \Om_{Is}\quad\hbox{and}\quad
\|d\b_I\|_{L^p(X_I,g_I)}\le C\l^{2/p-1/2}.
\eqlabel\eqnlabel{\eqMoreBetaEst}
$$
Second, for the derivatives of $\b_J$ with respect to $\la_I$, one has
$$\leqalignno
{|\rd\b_J/\rd\l_I|_{L^\8(X_J,g_J)}\le C\l_I^{-1}&\quad\hbox{and}\quad
|\rd d\b_J/\rd\l_I|_{L^\8(X_J,g_J)}\le C\l_I^{-3/2},  
&\numeqn\eqnlabel{\eqDBetaDLaEst}\cr
\|\rd\b_J/\rd\l_I\|_{L^p(X_J,g_J)}\le \l_I^{2/p-1}&\quad\hbox{and}\quad
\|\rd d\b_J/\rd\l_I\|_{L^p(X_J,g_J)}\le C\l_I^{2/p-3/2}, &\cr}
$$
for $J=I_-$ or $I$,
these derivatives being zero otherwise.
Third, for the derivatives of $\b_J$ with respect to $x_I$, one has
$$\leqalignno
{\|\rd\b_J/\rd p_I\|_{L^\8(X_J,g_J)}\le C\l_I^{-1/2}&\quad\hbox{and}\quad
\|\rd d\b_J/\rd p_I\|_{L^\8(X_J,g_J)}\le C\l_I^{-1},  
&\numeqn\eqnlabel{\eqDBetaDqEst}\cr
\|\rd\b_J/\rd p_I\|_{L^p(X_J,g_J)}\le \l_I^{2/p-1/2}&\quad\hbox{and}\quad
\|\rd d\b_J/\rd p_I\|_{L^p(X_J,g_J)}\le C\l_I^{2/p-1}, &\cr}
$$
for $J=I_-$ or $I$, these derivatives being zero otherwise. The cut-off functions $\ga_J$ also satisfy the bounds of Eqs. \eqref{\eqMoreBetaEst}, \eqref{\eqDBetaDLaEst} and \eqref{\eqDBetaDqEst}.

\proclaim{Proposition}{\thmlabel{\propDPDLaqEst}
For any $0<\d<\half$ and $2<p<4$ defined by $p=4/(1+2\d)$, and small enough $\la_0$, there is a constant $C=C(\de,g_0,\sT)$ such that for any $t\in\sT$ and $\xi\in L^p\Om^{+,g}(X,\ad P)$,
\item{\rm (a)}
$\|\pd{P}{\l_I}\xi\|_{L^4(X,g)}
\le C\l_I^{-1/2-\d}\|\xi\|_{L^p(X,g)}$,
\item{\rm (b)}
$\|\pd{P}{p_I}\xi\|_{L^4(X,g)}
\le C\l_I^{-\d}\|\xi\|_{L^p(X,g)}$.}

\pf (a) As $P=Q(1+R)^{-1}$, we first obtain
operator bounds for $\rd Q/\rd\l_I$, $\rd R/\rd\l_I$, and then deduce an
operator bound for $\rd P/\rd\l_I$. 

\step{1.} {\it Estimate for $\rd Q/\rd\l_I$}. Recall that $Q\xi=\sum_J Q_J\xi$, where $Q_J=\b_JP_J\c_J$ is independent of $\l_I$ for $J\ne I_-,I$, and so
$$\leqalignno
{\pd{Q}{\l_I}&=\pd{Q_{I_-}}{\l_I}+\pd{Q_I}{\l_I},\quad\hbox{where}\quad
\pd{Q_I}{\l_I}=\pd{\b_I}{\l_I}P_I\c_I+\b_IP_I\pd{c_I}{\l_I}, &\cr}
$$
with the analogous expression for $\rd Q_{I_-}/\rd \l_I$. 
Choose $4<q,q_1<\8$ and $2<p,p_1<4$ by setting
$$\leqalignno
{p&=4/(1+2\d)\quad\hbox{and}\quad q=4/(1-2\d) &\numeqn\cr
1/p&=1/4+1/q_1\quad\hbox{and}\quad 1/2=1/p_1+1/q_1, &\cr}
$$
and observe that
$1/4=1/q+1/q_1$ and $1/2=1/p+1/q$, while
$2/p=1/2+\d$ and $2/q=1/2-\d$. 
Applying H\"older's inequality, the operator bounds for $P_I$ of Lemma \ref{\lemPIEst}, 
and the fact that $\|\rd\b_I/\rd\l_I\|_{L^q}$ and
$\|\rd\c_I/\rd\l_I\|_{L^q}$ are bounded by $C\l_I^{2/q-1}$ from Eq. \eqref{\eqDBetaDLaEst}, we find 
$$\leqalignno
{\left\|\pd{Q_I}{\l_I}\xi\right\|_{L^4}
&\le C\left\|\pd{\b_I}{\l_I}\right\|_{L^q}\|\xi\|_{L^p}
+C\left\|\pd{\c_I}{\l_I}\right\|_{L^q}\|\xi\|_{L^p}
\le C\l_I^{2/q-1}\|\xi\|_{L^p}. &\cr}
$$
Combining the above estimate with the analogous bound for the $\rd Q_{I_-}/\rd \l_I$ term, we see that 
$$
\left\|\pd{Q}{\l_I}\xi\right\|_{L^4}\le 
C\l_I^{2/q-1}\|\xi\|_{L^p}, \eqlabel
$$
completing Step 1.

\step{2.}  {\it Estimate for $\rd R/\rd\l_I$}. We have $R=d_{A'}^{+,g}Q-1$ on $X$ and so differentiating with respect to $\l_I$ gives
$$\leqalignno
{\pd{R}{\l_I}\xi
&=\pd{*_g}{\l_I}d_{A'}Q\xi +\left[\pd{A'}{\l_I},Q\xi\right]^{+,g} + d_{A'}^{+,g}\pd{Q}{\l_I}\xi. &\cr}
$$
Using our $L^\8$ bound for $\rd *_g/\rd\l_I$ of Lemma \ref{\lemDHodgeStarDLaq}, the $L^4$ bound for $\rd A'/\rd\l_I$ of Proposition \ref{\propDAprimeDLaqEst}, and the operator norm bounds for $Q$ of Lemma \ref{\lemQEst}, we see that
$$\leqalignno
{\left\|\pd{R}{\l_I}\xi\right\|_{L^2} &\le C\l_I^{-1/2}\|d_{A'}Q\xi\|_{L^2}
+ C\|\xi\|_{L^2} 
+ \left\|d_{A'}^{+,g}\pd{Q}{\l_I}\xi\right\|_{L^2}. &\numeqn\cr}
$$
For the $d_{A'}Q$ term above, note that
$d_{A'}Q\xi=\sum_J d_{A_J'}Q_I\xi$ and writing $A_J'=A_J+a_J$ over $X_J'$, we have
$$
d_{A_J'}Q_J\xi
=d\b_J\wedge P_J\c_J\xi+\b_Jd_{A_J}P_J\c_J\xi+\b_J[a_J,P_J\c_J\xi]. 
$$
Using the bounds $\|d\b_J\|_{L^4}\le C$ of Eq. \eqref{\eqMoreBetaEst}, 
$\|a_J\|_{L^4}\le C\ol\la$ of Lemma \ref{\lemAIprimeEst}, H\"older's inequality, and the operator bounds for $P_J$ of Lemma \ref{\lemPIEst}, we find that 
$$
\|d_{A'}Q\xi\|_{L^2}\le C\|\xi\|_{L^2}. \eqlabel
$$
For the $d_{A'}^{+,g}\rd Q/\rd\l_I$ term, note that
$$
d_{A'}^{+,g}\pd{Q}{\l_I}=d_{A_{I_-}'}^{+,g}\pd{Q_{I_-}}{\l_I}
+d_{A_I'}^{+,g}\pd{Q_I}{\l_I},
$$
where, using $d_{A_I}^{+,g_I}P_I=1$ and $\b_I=1$ on $\supp\c_I$, we have
$$\leqalignno
{d_{A_I'}^{+,g}\pd{Q_I}{\l_I}\xi
&=\left(d\pd{\b_I}{\l_I}\wedge P_I\c_I\xi\right)^{+,g} + \pd{\b_I}{\l_I}[a_I,P_I\c_I\xi]^{+,g}  &\cr
&\quad +\left(d\b_I\wedge P_I\pd{\c_I}{\l_I}\xi\right)^{+,g} + \pd{\c_I}{\l_I}\xi+\b_I\left[a_I,P_I\pd{\c_I}{\l_I}\xi\right]^{+,g} &\cr
&\quad + {1\over 2}\pd{\b_I}{\l_I}(*_g-*_{g_I})d_{A_I}P_I\c_I\xi
+ {1\over 2}\b_I(*_g-*_{g_I})d_{A_I}P_I\pd{\c_I}{\l_I}\xi, &\cr}
$$
with the analogous expression for $d_{A_{I_-}'}^{+,g}\rd Q_{I_-}/\rd\l_I$. 
From Lemma \ref{\lemDHodgeStarDLaq} and Lemma \ref{\lemPIEst}, we see that
$$\leqalignno
{\left\|d_{A_I'}^{+,g}\pd{Q_I}{\l_I}\xi\right\|_{L^2}
&\le C\left\|d\pd{\b_I}{\l_I}\right\|_{L^{p_1}}\|\xi\|_{L^p}
+C\left\|\pd{\b_I}{\l_I}\right\|_{L^q}\|a_I\|_{L^4}\|\xi\|_{L^p} &\cr
&+C\|d\b_I\|_{L^4}\left\|\pd{\c_I}{\l_I}\right\|_{L^q}\|\xi\|_{L^p} 
+\left\|\pd{\c_I}{\l_I}\right\|_{L^q}\|\xi\|_{L^p} 
+C\|a_I\|_{L^4}\left\|\pd{\c_I}{\l_I}\right\|_{L^q}\|\xi\|_{L^p} &\cr
&+C\ol\la\left\|\pd{\b_I}{\l_I}\right\|_{L^q}\|\xi\|_{L^p} 
+C\ol\la\left\|\pd{\c_I}{\l_I}\right\|_{L^q}\|\xi\|_{L^p}, &\cr}
$$
Now $\|a_I\|_{L^4}\le C\ol\la$ by Lemma \ref{\lemAIprimeEst}, and from Eq. \eqref{\eqDBetaDLaEst}, we have that $\|\rd\b_I/\rd\l_I\|_{L^q}$ and
$\|\rd\c_I/\rd\l_I\|_{L^q}$ are bounded by $C\l_I^{2/q-1}$. Hence,
$$\leqalignno
{\left\|d_{A_I'}^{+,g}\pd{Q_I}{\l_I}\xi\right\|_{L^2}
&\le C\l_I^{2/q-1}\|\xi\|_{L^p}, &\cr}
$$
with the analogous bound for the $d_{A_{I_-}'}^{+,g}\rd Q_{I_-}/\rd\l_I$
term. Therefore,  
$$
\left\|d_{A'}^{+,g}\pd{Q}{\l_I}\xi\right\|_{L^2}
\le C\l_I^{2/q-1}\|\xi\|_{L^p}. \eqlabel
$$
Combining the above inequalities and noting that $\|\xi\|_{L^2(X,g)}\le C\|\xi\|_{L^p(X,g)}$, we have
$$
\left\|\pd{R}{\l_I}\xi\right\|_{L^4(X,g)}\le C\l_I^{2/q-1}\|\xi\|_{L^p},
\eqlabel
$$
which completes Step 2.

\step{3.}  {\it Estimate for $\rd P/\rd\l_I$}. Differentiating $P=Q(1+R)^{-1}$ with respect to $\l_I$ gives
$$
\pd{P}{\l_I}=\pd{Q}{\l_I}(1+R)^{-1}-Q(1+R)^{-1}\pd{R}{\l_I}(1+R)^{-1},
$$
and thus applying the bounds from Steps 1 and 2, we have
$$\leqalignno
{\left\|\pd{P}{\l_I}\xi\right\|_{L^4(X,g)}
\le C\l_I^{2/q-1}\|\xi\|_{L^p(X,g)}, &\cr}
$$
which yields (a) since $2/q-1=-1/2-\d$.
For (b), the strategy of (a)
shows that $\|\pd{Q}{p_I}\xi\|_{L^4}$ and 
$\|\pd{R}{p_I}\xi\|_{L^4(X,g)}$ are bounded by 
$C\l_I^{2/q-1/2}\|\xi\|_{L^p}$, giving
$$
\left\|\pd{P}{p_I}\xi\right\|_{L^4(X,g)}\le C\l_I^{2/q-1/2}\|\xi\|_{L^p},
\eqlabel
$$
and so (b) follows.\qed

As is readily verified, Lemma \ref{\lemDXiDLaqPreEst} and Proposition
\ref{\propDPDLaqEst} then provide the 
the following estimates for the derivatives of
$\xi$ and $a$ with respect to $\l_I$ and $x_I$:

\proclaim{Corollary}{\thmlabel{\corDXIDLaqEst} Let $\xi$ and $a=P\xi$ be as
in \/ {\rm Theorem \ref{\thmDKMap}} and assume that the conditions of that
theorem hold. Then, for small enough $\la_0>0$, 
there is a constant $C=C(\d,g_0,\sT)$ such that for any $t\in\sT$,
\item{\rm (a)} $\|\rd{\xi}/\rd{\l_I}\|_{L^2(X,g)}
\le C(1+\ol\la^{3/2+\d}\l_I^{-1/2-\d})$,
\item{\rm (b)} $\|\rd{\xi}/\rd{p_I}\|_{L^2(X,g)}
\le C(1+\ol\la^{3/2+\d}\l_I^{-\d})$,
\item{\rm (c)} $\|\rd{a}/\rd{\l_I}\|_{L^2(X,g)}
\le C(1+\ol\la^{1/2+\d}\l_I^{-1/2-\d})$,
\item{\rm (d)} $\|\rd{a}/\rd{p_I}\|_{L^2(X,g)}
\le C(1+\ol\la^{1/2+\d}\l_I^{-\d})$.}

With bounds for the derivatives of $\xi$ and $P$ with respect to $\la_I$
and $x_I$ at hand, we obtain our final estimates for the derivatives of the
anti-self-dual connections $A$ and $\hat A$. Since $A=A'+a$ and
combining Proposition \ref{\propDAprimeDLaqEst} and Corollary
\ref{\corDXIDLaqEst}, we have:

\proclaim{Corollary}{\thmlabel{\corDADLaq} Assume that the conditions of \/ {\rm Theorem \ref{\thmDKMap}} hold.
Then, for any $0<\d<1/2$ and small enough $\la_0>0$, 
there is a constant $C=C(\d,g_0,\sT)$ such that for any $t\in\sT$, the following bounds hold:
\item{\rm (a)} $\|\rd A/\rd{\l_I}\|_{L^2(X,g)} 
\le C(1+\ol\la^{1/2+\d}\l_I^{-1/2-\d})$, 
\item{\rm (b)} $\|\rd A/\rd{p_I}\|_{L^2(X,g)} 
\le C(1+\ol\la^{1/2+\d}\l_I^{-\d})$.}

\proclaim{Theorem}{\thmlabel{\thmDHatADLaq} Assume that the conditions of \/ {\rm Theorem \ref{\thmDKMap}} hold.
Then, for any $0<\d<1/2$ and small enough $\la_0>0$, 
there is a constant $C=C(\d,g_0,\sT)$ such that for any $t\in\sT$, the following bounds hold:
\item{\rm (a)} $\|\rd{\hat a}/\rd{\l_I}\|_{L^2(X_0,g_0)} 
\le C(1+\ol\la^{1/2+\d}\l_I^{-1/2-\d})$, 
\item{\rm (b)} $\|\rd{\hat A}/\rd{\l_I}\|_{L^2(X_0,g_0)} 
\le C(1+\ol\la^{1/2+\d}\l_I^{-1/2-\d})$,
\item{\rm (c)} $\|\rd{\hat a}/\rd{p_I}\|_{L^2(X_0,g_0)} 
\le C(1+\ol\la^{1/2+\d}\l_I^{-\d})$,
\item{\rm (d)} $\|\rd{\hat A}/\rd{p_I}\|_{L^2(X_0,g_0)} 
\le C(1+\ol\la^{1/2+\d}\l_I^{-\d})$.}

\pf Using the bound $\|\xi\|_{L^p}\le C\ol\la^{2/p}$ of Theorem
\ref{\thmDKExist}, the equality $2/p=1/2+\d$, the $L^2$ estimate for
$\rd\hat a/\rd\l_I$ in Lemma \ref{\lemDaDLaqPreEst}, the $L^2$ estimate for
$\rd\xi/\rd\l_I$ in Corollary \ref{\corDXIDLaqEst}, and the 
operator estimate for $\rd P/\rd\l_I$ in Proposition \ref{\propDPDLaqEst},
we obtain 
$$
\left\|\pd{\hat a}{\l_I}\right\|_{L^2(X_0,g_0)} \le
C\left(\ol\la^{1/2+\d}\l_I^{-1/2-\d}
+1+\ol\la^{3/2+\d}\l_I^{-1/2-\d}+\ol\la^{1/2+\d}\l_I^{-1/2}\right),
$$
which gives (a). Then (b) follows from (a) and the estimate $\|\rd\hat A'/\rd\l_I\|_{L^2(X_0,g_0)}\le C$ of Proposition \ref{\propDHatAprimeDLaqEst}.
The proofs of (c) and (d) are similar.
\qed


\section{Derivatives with respect to lower moduli}
In this section we obtain estimates for the derivatives of the family of
$g_0$-anti-self-dual connections $\hat A$ with respect to the lower moduli
parameters $t_I\in T_{A_I}$.   
Just as in \S\ref{\secDHatADLaq}, the strategy is to use the $g$-anti-self-dual equation of Eq. \eqref{\eqASDXi}, together with its derivatives with respect to the $t_I$ parameters, to first obtain estimates for the derivatives of   
$a$ and $\xi$, and then the required derivatives of $\hat a$ and $\hat A'$. 
The Sobolev exponents $p,q$ are fixed so that $2\le p<4$ and $4\le q<\8$, where $q$ is determined by $1/p=1/4+1/q$.
We have the following preliminary estimates for the derivatives of $\xi$ and $a$.

\proclaim{Lemma}{\thmlabel{\lemDaDtDxiDt} 
Let $\xi$ and $a=P\xi$ be as in \/ {\rm Theorem \ref{\thmDKMap}}, and assume that the conditions of that theorem hold. Then, for small enough $\la_0>0$, 
there is a constant $C=C(g_0,p,\sT)$ such that for any $t\in\sT$,
\item{\rm (a)}
$\|\rd{a}/\rd{t_I}\|_{L^p(X,g)}\le
C\|\rd{\xi}/\rd{t_I}\|_{L^p(X,g)}+\|\pd{P}{t_I}\xi\|_{L^p(X,g)}$, 
\item{\rm (b)}
$\|\rd{\xi}/\rd{t_I}\|_{L^p(X,g)}
\le C\left(\ol\la^{2/p-1/2}+\ol\la^{2/p}
\|\pd{P}{t_I}\xi\|_{L^4(X,g)}\right)$.}

\pf Differentiating Eq. \eqref{\eqASDXi} with respect to $t_I$ gives
$$\leqalignno
{\pd{\xi}{t_I}&=\pd{\eta}{t_I}-\left(\pd{P\xi}{t_I}\wedge P\xi\right)^{+,g}
-\left(P\xi\wedge \pd{P\xi}{t_I}\right)^{+,g}, &\cr
\pd{P\xi}{t_I} &=\pd{P}{t_I}\xi + P\pd{\xi}{t_I}. &\cr}
$$
The proofs of (a) and (b) are then similar to those of Lemmas \ref{\lemDaDLaqPreEst} and \ref{\lemDXiDLaqPreEst}. \qed

Thus, an operator estimate for $\rd P/\rd t_I$ is required. Since $P=Q(1+R)^{-1}$, we have
$$
\pd{P}{t_I}=\pd{Q}{t_I}(1+R)^{-1}-Q(1+R)^{-1}\pd{R}{t_I}(1+R)^{-1}. 
\eqlabel\eqnlabel{\eqDPDt}
$$
We recall that $P_I=d_{A_I}^{*,g_I}G_{A_I}^{+,g_I}$. Differentiating with
respect to $t_I$, we obtain 
$$
\pd{P_I}{t_I}=\pd{d_{A_I}^{*,g_I}}{t_I}G_{A_I}^{+,g_I}
-d_{A_I}^{*,g_I}G_{A_I}^{+,g_I}\pd{\De_{A_I}^{+,g_I}}{t_I}G_{A_I}^{+,g_I}.
$$
The derivatives of $d_{A_I}^{+,g_I}$ and $d_{A_I}^{*,g_I}$ with respect to $t_I$ are given by
$$\leqalignno
{\pd{d_{A_I}^{+,g_I}}{t_I}\w &=\left[\pd{A_I}{t_I},\w\right]^{+,g_I} 
=\left[\pd{A_I}{t_I},\cdot\right]^{+,g_I}\w, &\cr
\pd{d_{A_I}^{*,g_I}}{t_I}\xi &=-*\left[\pd{A_I}{t_I},*\xi\right] 
=\left[\pd{A_I}{t_I},\cdot\right]^*\xi, &\cr}
$$
for any $\w\in \W^1(X_I,\ad P_I)$ and $\xi\in \W^{+,g_I}(X_I,\ad P_I)$.
Therefore,
$$\leqalignno
{\pd{\De_{A_I}^{+,g_I}}{t_I} 
&=\left[\pd{A_I}{t_I},\cdot\right]^{+,g_I}d_{A_I}^{*,g_I}
+d_{A_I}^{+,g_I}\left[\pd{A_I}{t_I},\cdot\right]^*, &\cr}
$$
and so we find that
$$ 
\pd{P_I}{t_I}=(1-P_Id_{A_I}^{+,g_I})
\left[\pd{A_I}{t_I},\cdot\right]^*G_{A_I}^{+,g_I}
-P_I\left[\pd{A_I}{t_I},\cdot\right]^{+,g_I}P_I. 
\eqlabel\eqnlabel{\eqDPIDt}
$$
Note that $1-P_Id_{A_I}^{+,g_I}$ is a bounded $(L^q,L^q)$ operator on $\Om^{+,g_I}(X_I,\ad P_I)$ by the Calderon-Zygmund theory. 

\proclaim{Lemma}{\thmlabel{\lemDPIDtDQDt}
There is a constant $C=C(g_0,p,\sT)$ such that for any $t\in\sT$,
\item{\rm (a)}
$\|\pd{P_I}{t_I}\xi\|_{L^q(X_I,g_I)}
\le C\|\xi\|_{L^p(X_I,g_I)}$, 
for $\xi\in L^p\W^{+,g_I}(X_I,g_I)$,
\item{\rm (b)}
$\|\pd{Q}{t_I}\xi\|_{L^q(X,g)}\le C\|\xi\|_{L^p(X,g)}$, 
for $\xi\in L^p\Om^{+,g}(X,g)$.}

\pf Since $1-P_Id_{A_I}^{+,g_I}$ is bounded on $L^q(X_I,g_I)$, then Eq. \eqref{\eqDPIDt} and the H\"older inequalities show that
$$\leqalignno
{\left\|\pd{P_I}{t_I}\xi\right\|_{L^q}
&\le \left\|\pd{A_I}{t_I}\right\|_{L^\8}\|G_{A_I}^{+,g_I}\xi\|_{L^q}
+C\left\|\pd{A_I}{t_I}\right\|_{L^4}\|P_I\xi\|_{L^q}. &\cr}
$$
But $G_{A_I}^{+,g_I}$ and $P_I$ are bounded $(L^p, L^q)$ operators and noting that the family $A_I(t_I)$ is smoothly parametrised by  
$t_I\in T_{A_I}$, we obtain (a). Now $Q_I=\b_IP_I\c_I$ and
$Q=\sum_I Q_I$, so inequality (b) follows. \qed

It remains to estimate the derivative of $R$ with respect to $t_I$. 

\proclaim{Lemma}{\thmlabel{\lemDRDt}
There exists a constant $C=C(g_0,p,\sT)$ such that for any $t\in\sT$ and
$\xi\in\Om^{+,g}(X,g)$, we have
$\|\pd{R}{t_I}\xi\|_{L^p(X,g)}\le C\|\xi\|_{L^p(X,g)}$.}

\pf We recall that $R=d_{A'}^{+,g}Q-1$ over $X$ and 
$R=d_{A_I'}^{+,g}Q_I-1$ over $X_I$.
Writing $A_I'=A_I+a_I$, we find that
$$\leqalignno
{R&=d\b_I\wedge P_I\c_I + \b_Id_{A_I}^{+,g_I}P_I\c_I
+ \b_I[a_I,\cdot\ ]^{+,g}P_I\c_I 
+ {1\over 2}(*_g-*_{g_I})\b_Id_{A_I}P_I - 1. &\cr}
$$
Noting that $d_{A_I}^{+,g_I}P_I=1$ and
differentiating with respect to $t_I$, we have
$$\leqalignno
{\pd{R}{t_I}\xi
&=d\b_I\wedge\pd{P_I}{t_I}\c_I\xi
+ \b_I\left[\pd{a_I}{t_I},P_I\c_I\xi\right]^{+,g}
+\b_I\left[a_I,\pd{P_I}{t_I}\c_I\xi\right]^{+,g} 
&\numeqn\eqnlabel{\eqDRDt}\cr
&\qquad + {1\over 2}(*_g-*_{g_I})\b_I\left[\pd{A_I}{t_I},P_I\c_I\xi\right] 
+ {1\over 2}(*_g-*_{g_I})\b_Id_{A_I}\pd{P_I}{t_I}\c_I\xi, 
&\cr}
$$
and so
$$\leqalignno
{\left\|\pd{R}{t_I}\xi\right\|_{L^p}
&\le C\|d\b_I\|_{L^4}\left\|\pd{P_I}{t_I}\c_I\xi\right\|_{L^q}
+ C\left\|\pd{a_I}{t_I}\right\|_{L^4}\|P_I\c_I\xi\|_{L^q} 
&\numeqn\eqnlabel{\eqLpDRDt}\cr
&\qquad + C\|a_I\|_{L^4}\left\|\pd{P_I}{t_I}\c_I\xi\right\|_{L^q}
+ C\ol\l\left\|\pd{A_I}{t_I}\right\|_{L^4}\|P_I\c_I\xi\|_{L^q} &\cr
&\qquad + C\ol\l\left\|d_{A_I}\pd{P_I}{t_I}\c_I\xi\right\|_{L^p}, &\cr}
$$
where $a_I=(\psi_I-1)\s_I^*A_I$ and 
$\rd{a_I}/\rd{t_I}=(\psi_I-1)\s_I^*\rd{A_I}/\rd{t_I}$. Aside from
the self-dual projection and factor $\be_I$, 
the last term on the right-hand side of Eq. \eqref{\eqDRDt} is given by
$$\leqalignno
{d_{A_I}\pd{P_I}{t_I}\c_I\xi
&=-d_{A_I}*\left[\pd{A_I}{t_I},G_{A_I}^{+,g_I}\c_I\xi\right]
+d_{A_I}P_Id_{A_I}^{+,g_I}*
\left[\pd{A_I}{t_I},G_{A_I}^{+,g_I}\c_I\xi\right] &\cr
&\qquad -d_{A_I}P_I\left[\pd{A_I}{t_I},P_I\c_I\xi\right]^{+,g}. &\cr}
$$
Since $P_I$ is a bounded operator from $L^p$ to $L^p_1$  
and using the bounded inclusion $L^p_1\to L^q$, we see that
$$\leqalignno
{\left\|d_{A_I}\pd{P_I}{t_I}\c_I\xi\right\|_{L^p}
&\le C\left\|d_{A_I}*\left[\pd{A_I}{t_I},
G_{A_I}^{+,g_I}\c_I\xi\right]\right\|_{L^p}
+\left\|\left[\pd{A_I}{t_I},P_I\c_I\xi\right]\right\|_{L^p} &\cr
&\le C\left\|\pd{A_I}{t_I}\right\|_{L^\8}
\|G_{A_I}^{+,g_I}\c_I\xi\|_{L^p}
+\left\|\]^{A_I}\pd{A_I}{t_I}\right\|_{L^\8}
\|G_{A_I}^{+,g_I}\c_I\xi\|_{L^p} &\cr
&\qquad + C\left\|\pd{A_I}{t_I}\right\|_{L^\8}
\|\]^{A_I}G_{A_I}^{+,g_I}\c_I\xi\|_{L^p}
+ \left\|\pd{A_I}{t_I}\right\|_{L^4}\|P_I\c_I\xi\|_{L^q}. &\cr}
$$
Since the family $A_I(t)$ is smoothly parametrised by $t_I\in T_{A_I}$ and as $G_{A_I}^{+,g_I}$ is a bounded operator from $L^p$ to $L^p_2$, we have
$$
\left\|d_{A_I}\pd{P_I}{t_I}\c_I\xi\right\|_{L^p} \le C\|\xi\|_{L^p}. 
\eqlabel\eqnlabel{\eqdDPIDtgaxi}
$$
Eqs. \eqref{\eqLpDRDt}, \eqref{\eqdDPIDtgaxi} and Lemma \ref{\lemDPIDtDQDt} then yield the required bound for $\rd{R}/\rd{t_I}$. \qed

Thus, Eq. \eqref{\eqDPDt}, together with
Lemmas \ref{\lemDPIDtDQDt} and \ref{\lemDRDt}, provides an estimate for 
the derivative of $P$ with respect to $t_I$:

\proclaim{Proposition}{\thmlabel{\propDPDt}
There is a constant $C=C(g_0,p,\sT)$ such that for any $t\in\sT$ and
$\xi\in L^p\Om^{+,g}(X,\ad P)$, we have
$\|\pd{P}{t_I}\xi\|_{L^q(X,g)}\le C\|\xi\|_{L^p(X,g)}$.}

This leads to our final estimates for the derivatives of $\xi$ and $a$
with respect to $t_I$. 

\proclaim{Corollary}{\thmlabel{\corDaDtDxiDt}
Let $\xi$ and $a=P\xi$ be as in \/ {\rm Theorem \ref{\thmDKMap}}, and assume
that the conditions of that theorem hold. Then, for small enough $\la_0>0$, 
there is a constant $C=C(g_0,p,\sT)$ such that for any $t\in\sT$,
\item{\rm (a)} $\|\rd{\xi}/\rd{t_I}\|_{L^p(X,g)}\le C\ol\l^{2/p-1/2}$,
\item{\rm (b)} $\|\rd{a}/\rd{t_I}\|_{L^p(X,g)}\le C\ol\l^{2/p-1/2}$.}

\pf Inequality (a) follows from Lemma \ref{\lemDaDtDxiDt} and Proposition \ref{\propDPDt}, since $\|\xi\|_{L^p}\le C\ol\l^{2/p}$ by Theorem \ref{\thmDKExist}. Inequality (b) then follows from (a) and 
Lemma \ref{\lemDaDtDxiDt}.\qed

By combining Proposition \ref{\propDAprimeDt} and Corollary \ref{\corDaDtDxiDt}
we obtain an estimate for the derivatives of the connections $A=A'+a$ over
$X$: 

\proclaim{Corollary}{\thmlabel{\corDADtI} Assume that the conditions of \/ {\rm Theorem \ref{\thmDKMap}} hold. Then,
for any $2\le p<4$ and sufficiently small $\l_0>0$, there is a constant $C=C(g_0,p,\sT)$ such that for any $t\in\sT$, 
\item{\rm (a)} 
$\|\rd{A}/\rd{t_I}-\rd{A_I}/\rd{t_I}\|_{L^p(X_I'',g_I)}
\le C\ol\la^{2/p-1/2}$,
\item{\rm (b)} $\|\rd{A}/\rd{t_I}\|_{L^p(X,g)}\le C$.}

We now come to the main result of this section.

\proclaim{Theorem}{\thmlabel{\thmDHatADtI} Assume that the conditions of \/ {\rm Theorem \ref{\thmDKMap}} hold. Then,
for any $2\le p<4$ and sufficiently small $\l_0>0$, there is a constant $C=C(g_0,p,\sT)$ such that for any $t\in\sT$, 
\item{\rm (a)} 
$\|\rd{\hat a}/\rd{t_I}\|_{L^p(X_0,g_0)}\le C\ol\l^{2/p-1/2}$
\item{\rm (b)} $\|\rd{\hat A}/\rd{t_I}\|_{L^p(X_0,g_0)}\le C$.}

\pf Let $U\equiv f_0^{-1}\cdots f_I^{-1}(X_I')\subset X_0$ and note that
$\pd{\hat a}{t_I}= \sum_I f_0^*\cdots f_I^*\pd{a}{t_I}$ on $U$.
Now Lemma \ref{\lemXferEst} gives
$$\leqalignno
{\left\|\pd{\hat a}{t_I}\right\|_{L^p(X_0,g_0)}
&\le C\sum_I\left\|f_0^*\cdots f_I^*\pd{a}{t_I}\right\|_{L^p(X_0,g_0)} 
\le C\left\|\pd{a}{t_I}\right\|_{L^p(X,g)}, &\cr}
$$
and so Part (a) follows from Corollary \ref{\corDaDtDxiDt}. Then Part (b) follows from (a), the estimate
$\|\rd\hat A'/\rd t_I\|_{L^p(X_0,g_0)}\le C$ of Proposition \ref{\propDhatAprimeDt}. \qed


\section{Derivatives with respect to bundle gluing parameters}
We obtain estimates for the partial derivatives of the family of $g_0$-anti-self-dual connections $\hat A(t)$ with respect to the bundle gluing parameters $\rho_I\in\Gl_I$. The Sobolev exponents $p,q$ are fixed so that $2\le p<4$, with $4\le q<\8$ determined by $1/4+1/q=1/p$. 
We first recall the estimate of Donaldson and Kronheimer for the derivative of $a=P\xi$ with respect to the gluing parameters $\rho_I$. As described in \S\ref{\secDAprimeDGlu}, we work with an equivalent family of $g$-anti-self-dual connections $A=A'+a$ on a fixed bundle $P$. Thus, considering only the gluing parameters, we have a diffeomorphism $B_\frak{g}\ni v\to A(\bar\rho_I,v)\in\cA_{X,P}^*$ (where $B_\frak{g}$ is the unit ball in $\frak{g}$), giving a family of $C^\8$ connections on a fixed bundle $P=P(\bar\rho_I)$, as in Eq. \eqref{\eqAprimeGluChart}. Here,
$B_\frak{g}\ni v\to\rho_I(v)=\bar\rho_I\exp(v)\in\Gl_I$ 
is a coordinate chart centred at $\bar\rho_I\in\Gl_I$, as in Eq. \eqref{\eqGluChart}. This understood, one has the following bounds.

\proclaim{Proposition}{{\rm [D-K, p. 303]}\thmlabel{\propDaDvEst} 
Let $a$ be as in \/{\rm Theorem \ref{\thmDKMap}}, and assume that the conditions of that theorem hold. Then, for small enough $\la_0>0$, 
there is a constant $C=C(g_0,p,\sT)$ such that for any $t\in\sT$,
$\|\rd{a}/\rd{v}\|_{L^q(X,g)}\le C\ol\la^{2/p+1}$.}

\pf The proof in [D-K] deals only with single connected sums $X=X_0\# X_1$, but the argument adapts without significant change to the general case of multiple connected sums $\#_{I\in\sI}X_I$. Likewise, the assumptions in [D-K] that $\Ga_I=1$ and $H_{A_I}^0=H_{A_I}^1=0$, for all $I$, do not affect the relevant estimates.\qed

\proclaim{Corollary}{\thmlabel{\corDADvEst} 
Let $A$ be as in \/{\rm Theorem \ref{\thmDKMap}}, and assume that the conditions of that theorem hold. Then, for small enough $\la_0>0$, 
there is a constant $C=C(g_0,p,\sT)$ such that for any $t\in\sT$,
$\|\rd{A}/\rd{v}\|_{L^p(X,g)}\le C\ol\la^{2/p-1/2}$.}

\pf Combine Propositions \ref{\propDAprimeDvEst} and \ref{\propDaDvEst}. \qed

Moreover, we have the following estimates for the derivatives of the $g_0$-anti-self-dual connections $\hat A=\hat A'+\hat a$ on 
the fixed bundle $\hat P$ over $X_0$.

\proclaim{Theorem}{\thmlabel{\thmDHatADv}
Assume that the conditions of \/{\rm Theorem \ref{\thmDKMap}} hold. Then, for small enough $\la_0>0$, 
there is a constant $C=C(g_0,p,\sT)$ such that for any $t\in\sT$,
\item{\rm (a)}
$\|\rd{\hat a}/\rd{v}\|_{L^p(X,g)}\le C\ol\la^{2/p+1}$,
\item{\rm (b)}
$\|\rd{\hat A}/\rd{v}\|_{L^p(X,g)}\le C\ol\la^{2/p-1/2}$.}

\pf Since $\hat a=f_0^*\cdots f_I^*a$ on 
$U\equiv f_0^{-1}\circ\cdots\circ f_I^{-1}(X_I')$,
Lemma \ref{\lemXferEst} gives
$$
\left\|f_0^*\cdots f_I^*\pd{a}{v}\right\|_{L^p(U,g_0)} 
\le C\left\|\pd{a}{v}\right\|_{L^p(X_I',g_I)},
$$
and Proposition \ref{\propDaDvEst} gives (a). Similarly, (b) follows from (a) and Proposition \ref{\propDhatAprimeDvEst}.\qed


\section{Differentials of the gluing maps and final arguments}
We summarise the results of the preceding sections and record our bounds for the differentials of the approximate gluing maps $\script{J}$ and $\hat\script{J}$. The estimates for $D\hat\sJ$ then give bounds for the diagonal (and so all) components of the $L^2$ metric ${\bf g}$ and completes the proof of Theorem \ref{\thmFinVolDiam}. Combining these metric bounds with results of Donaldson in [D89] then completes the proof of  Theorem \ref{\thmMetCompletion}.
The following two theorems summarise the estimates obtained in \S\S 5.3 to
5.5, the first following from Corollaries \ref{\corDADLaq},
\ref{\corDADtI}, and \ref{\corDADvEst} and the second from
Theorems \ref{\thmDHatADLaq}, \ref{\thmDHatADtI}, and \ref{\thmDHatADv}. 

\proclaim{Theorem}{\thmlabel{\thmDJEst} Let $\script{J}:\script{T}/\Ga\to M_{X,k}^*(g)$ be a gluing map and assume that the conditions of\/ {\rm Theorem \ref{\thmDKMap}} hold.
Then for sufficiently small $\l_0>0$ and any $t\in \script{T}$, there exists a constant $C=C(g_0,\script{T})$ such that the following bounds hold:
\item{\rm (a)} $\|D\sJ(\rd/\rd t_I^\al)\|_{L^2(X,g)}\le C$,
\item{\rm (b)} $\|D\sJ(\rd/\rd\rho_I^\be)\|_{L^2(X,g)} 
\le C\ol\la^{1/2}$, 
\item{\rm (c)} $\|D\sJ(\rd/\rd x_I^\mu)\|_{L^2(X,g)} 
\le C(1+\ol\la^{1/2+\de}\la_I^{-\de})$, 
\item{\rm (d)} $\|D\sJ(\rd/\rd\la_I)\|_{L^2(X,g)} 
\le C(1+\ol\la^{1/2+\de}\la_I^{-1/2-\de})$.}

\proclaim{Theorem}{\thmlabel{\thmDhatJEst} Let $\hat\sJ:\script{T}/\Ga\to M_{X_0,k}^*(g_0)$ be a gluing map and assume that the conditions of\/ {\rm Theorem \ref{\thmDKMap}} hold. 
Then for any $0<\de<1/2$, sufficiently small $\l_0>0$ and any $t\in \script{T}$, there exists a constant $C=C(\de,g_0,\script{T})$ such that the following bounds hold:
\item{\rm (a)} $\|D{\hat\sJ}(\rd/\rd t_I^\al)\|_{L^2(X_0,g_0)} \le C$, 
\item{\rm (b)} $\|D{\hat\sJ}(\rd/\rd \rho_I^\be)\|_{L^2(X_0,g_0)} 
\le C\ol\la^{1/2}$, 
\item{\rm (c)} $\|D{\hat\sJ}(\rd/\rd x_I^\mu)\|_{L^2(X_0,g_0)} 
\le C(1+\ol\la^{1/2+\de}\la_I^{-\de})$, 
\item{\rm (d)} $\|D{\hat\sJ}(\rd/\rd \la_I)\|_{L^2(X_0,g_0)} 
\le C(1+\ol\la^{1/2+\de}\la_I^{-1/2-\de})$.}

It remains to reinterpret the bounds of Theorem \ref{\thmDhatJEst} in terms of the corresponding bounds for the diagonal components of the $L^2$ metric ${\bf g}$. 

\proclaim{Corollary}{\thmlabel{\corL2MetBounds} Under the hypotheses of\/ {\rm Theorem \ref{\thmDhatJEst}}, the following bounds hold:
\item{\rm (a)} ${\bf g}(\rd/\rd t_I^\al,\rd/\rd t_I^\al) \le C$,
\item{\rm (b)} ${\bf g}(\rd/\rd \rho_I^\be,\rd/\rd \rho_I^\be) \le C\ol\la$,
\item{\rm (c)} ${\bf g}(\rd/\rd x_I^\mu,\rd/\rd x_I^\mu) \le C(1+\ol\la^{1+2\de}\la_I^{-2\de})$, 
\item{\rm (d)} ${\bf g}(\rd/\rd \la_I,\rd/\rd \la_I) \le
C(1+\ol\la^{1+2\de}\la_I^{-1-2\de})$.} 

Recall that the ${\bf g}$-length of a path $(s_0,s_1)\ni s\to A(s)\in M_{X_0,k}^*(g_0)$ is computed by
$$
\int_{s_0}^{s_1}\sqrt{\bf g}\left(\pd{A}{s},\pd{A}{s}\right)\,ds
\le \int_{s_0}^{s_1}\left\|\pd{A}{s}\right\|_{L^2(X_0,g_0)}\,ds.
$$
The proofs of our main results are now essentially complete. 
\medskip

\noindent{\it Proof of Theorem \ref{\thmFinVolDiam}}: Since $0<\de<1/2$, the bounds of Theorem \ref{\corL2MetBounds} imply that the gluing neighbourhoods $\sV=\hat\sJ(\sT^0/\Ga)$ have finite ${\bf g}$-volume and ${\bf g}$-diameter. Therefore, the bubbling ends of $M_{X_0,k}^*(g_0)$ have finite ${\bf g}$-volume and ${\bf g}$-diameter since the entire moduli space is covered by finitely many such neighbourhoods. Away from the Uhlenbeck boundary, gluing neighbourhoods consist simply of $C^\8$ Kuranishi charts. The conical ends corresponding to Kuranishi charts around the reducible connections have finite ${\bf g}$-volume and ${\bf g}$-diameter by
Theorem 1 [G-P89]. \qed

Next we consider the relationship between the metric completion and the
Uhlenbeck compactification of the anti-self-dual moduli space.
Let $d_2$ be the distance function on $M_{X_0,k}^*(g_0)$ defined by the
$L^2$ metric ${\bf g}$. Thus, if $[A]$, $[B]$ are two points in
$M_{X_0,k}^*(g_0)$, then $d_2([A],[B])$ is the infimum over all ${\bf
g}$-lengths of paths in $M_{X_0,k}^*(g_0)$ joining $[A]$, $[B]$. If the two
points lie in different path components of the moduli space, then set
$d_2([A],[B])=\8$. Since $b^+(X_0)=0$, the moduli space has at most
finitely many path components; we say that $M_{X_0,k}^*(g_0)$ has {\it
finite {\bf g}-diameter} if the sum of the {\bf g}-diameters of the
connected components is finite. 
In [D89], Donaldson constructs two other distance functions, $D_2$ and $D_2^\eps$, for any fixed $\eps>0$. First, given points $[A]$, $[B]$ in $\cB_{X_0,k}^*$, set
$$
D_2([A],[B])\equiv \inf_{u\in\cG}\|A-u^*B\|_{L^2(X_0,g_0)}. 
$$
Lemma 2 [D89] (or Lemma 4.2.4 [D-K]) show that $D_2$ is a well-defined distance function on $\cB_{X_0,k}^*$. Moreover, Lemma 1 [D89] shows that $D_2([A],[B])$ is equal to the distance function defined in the usual way by the $L^2$ metric on $\cB_{X_0,k}^*$ as the infimum over ${\bf g}$-lengths of paths in $\cB_{X_0,k}^*$ joining $[A]$ and $[B]$. One then obtains a second distance function on $M_{X_0,k}^*(g_0)$ by restriction.
Define an $\eps$-neighbourhood of $M_{X_0,k}^*(g_0)$ in $\cB_{X_0,k}^*$ by 
$$
\cB_{X_0,k}^{*,\eps}\equiv\{[A]\in\cB_{X_0,k}^* :\|F_A^{+,g_0}\|_{L^2(X_0,g_0)}<\eps\}.
$$ 
Then $D_2^\eps([A],[B])$ is defined as infimum of the ${\bf g}$-lengths of paths in $\cB_{X_0,k}^{*,\eps}$ joining two points $[A]$ and $[B]$ in $\cB_{X_0,k}^{*,\eps}$. One now obtains a third distance function on $M_{X_0,k}^*(g_0)$ by restriction.
The three distance functions $d_2$, $D_2$, and $D_2^\eps$ on $M_{X_0,k}^*(g_0)$ are related by
$$
D_2([A],[B])\le D_2^\eps([A],[B])\le d_2([A],[B]),\qquad [A], [B]\in M_{X_0,k}^*(g_0). \eqlabel\eqnlabel{\eqDistFnComp}
$$
To show that the $d_2$-completion of $M_{X_0,k}^*(g_0)$ is homeomorphic to the Uhlenbeck compactification $\ol{M}_{X_0,k}^u(g_0)$, it is enough to show that a sequence $[A^\al]$ in $M_{X_0,k}^*(g_0)$ is $d_2$-Cauchy if and only if it is convergent in the Uhlenbeck topology. For the metric $D_2^\eps$, one has

\proclaim{Theorem}{\thmlabel{\thmD2EpsCompletion}{\rm [D89, Theorem 4]} For any $\eps>0$, the $D_2^\eps$-completion of $M_{X_0,k}^*(g_0)$ is homeomorphic to $\ol{M}_{X_0,k}^u(g_0)$.\qed}

Thus Donaldson's result gives part of the proof of
Theorem \ref{\thmMetCompletion}: Suppose a sequence
$[A^\al]$ in $M_{X_0,k}^*(g_0)$ is $d_2$-Cauchy. According to Eq.
\eqref{\eqDistFnComp}, it must also be $D_2^\eps$-Cauchy and so is
convergent in the Uhlenbeck topology by Theorem \ref{\thmD2EpsCompletion}
or simply by Proposition 6 [D89]. The proof of the reverse direction,
namely that a sequence $[A^\al]$ which is convergent in the Uhlenbeck
toplogy is also $d_2$-Cauchy, is included in [F94].

\endchapter


\centerline{\bf References}
\bigskip
\item{[A-H-S]} M. F. Atiyah, N. J. Hitchin, and I. M. Singer, {\it
Self-duality in four-dimensional Riemannian geometry}, Proc. Royal
Soc. London {\bf A 362} (1978) 425-461.

\item{[D-M-M]} H. Doi, Y. Matsumoto, and T. Matumoto, {\it An
explicit formula of the metric on the moduli space of
BPST-instantons over ${{\Bbb S}^4}$},  A F\^ete of Topology,
(Y. Matsumoto et al., eds.), Academic, New York, 1988, pp. 543-556. 

\item{[D86]} S. K. Donaldson, {\it Connections, cohomology and the
intersection forms of four manifolds}, J. Differential Geometry
{\bf 24} (1986) 275-341.

\item{[D87]} S. K. Donaldson, {\it The orientation of Yang-Mills moduli
spaces and 4-manifold topology}, J. Differential Geometry
{\bf 26} (1987) 397--428.

\item{[D89]} S. K. Donaldson, {\it Compactification and completion
of Yang-Mills moduli spaces}, Differential Geometry, Lecture
Notes in Mathematics {\bf 1410} (F.~J. Carreras et al., eds.)
Springer-Verlag, New York, 1989, pp. 145-160.

\item{[D90]} S. K. Donaldson, {\it Instantons in Yang-Mills
theory}, Proceedings of the IMA Conference on Geometry and
Particle Physics, Oxford 1988, (F. Tsou, ed.), Oxford University Press,
Oxford, 1990, pp. 59-75.

\item{[D-K]} S. K. Donaldson and P. B. Kronheimer, {\it The Geometry
of Four-Manifolds}, Oxford University Press, Oxford, 1990.

\item{[F92]} P. M. N. Feehan, {\it Geometry of the moduli space of self-dual connections over the four-sphere}, Ph.D. Thesis, Columbia University, 1992.

\item{[F94]} P. M. N. Feehan, {\it On the geometry of the moduli
space of anti-self-dual connections}, preprint.

\item{[F-U]} D. S. Freed and K. K. Uhlenbeck, {\it Instantons and
Four-Manifolds}, second edition,\hfill\break Springer-Verlag, New
York, 1991.

\item{[G90]} D. Groisser, {\it The geometry of the moduli space of $\cp^2$
instantons}, Invent. Math. {\bf 99} (1990) 393-409. 

\item{[G92]} D. Groisser, {\it Curvature of Yang-Mills moduli spaces near the
boundary, I}, Comm. Anal. Geom. {\bf 1} (1993), 139-215.

\item{[G-P87]} D. Groisser and T. H. Parker, {\it The Riemannian
geometry of the Yang-Mills moduli space}, Commun. Math. Phys. {\bf
112} (1987) 663-689.

\item{[G-P89]} D. Groisser and T. H. Parker, {\it The geometry of
the Yang-Mills moduli space for definite manifolds}, J. Differential Geometry
{\bf 29} (1989) 499-544.

\item {[Ha]} L. Habermann, {\it On the geometry of the space of
$Sp(1)$-instantons with Pontrjagin index $1$ on the $4$-sphere}, Ann.
Global Anal. Geom. {\bf 6} (1988) 3-29.

\item{[Mas]} H. Masur, {\it The extension of the
Weil-Petersson metric to the boundary of Teichm\"uller space}, Duke
Math. J. {\bf 43} (1976) 623-635.

\item{[P-W]} T. H. Parker and J. W. Wolfson, {\it Pseudo-holomorphic maps and bubble trees}, J. Geometric Analysis {\bf 3} (1993).

\item{[Pe]} Peng, X-W., {\it Asymptotic behavior of the $L^2$ metric on
moduli spaces of Yang-Mills connections}, Math. Z. {\bf 220} (1995) 127-158.

\item{[S-U]} J. Sacks and K. K. Uhlenbeck, {\it The existence of minimal 2-spheres}, Ann. Math. {\bf 113} (1981) 1-24.

\item{[T82]} C. H. Taubes, {\it Self-dual Yang-Mills
connections on non-self-dual 4-manifolds}, J. Differential Geometry
{\bf 17} (1982) 139-170.

\item{[T84a]} C. H. Taubes, {\it Self-dual connections on
4-manifolds with indefinite intersection matrix}, 
J. Differential Geometry {\bf 19} (1984) 517--560.

\item{[T84b]} C. H. Taubes, {\it Path-connected Yang-Mills moduli spaces},
J. Differential Geometry {\bf 19}, 337--392.

\item{[T88]} C. H. Taubes, {\it A framework for Morse theory for
the Yang-Mills functional}, Invent. Math. {\bf 94} (1988) 327-402. 

\item{[T92]} C. H. Taubes, {\it The existence of anti-self-dual
conformal structures}, J. Differential Geometry {\bf 36} (1992) 163-253.

\item{[U82]} K. K. Uhlenbeck, {\it Removable singularities in Yang-Mills fields}, Commun. Math. Phys. {\bf 83} (1982) 11-29.
\bigskip

\leftline{Department of Mathematics, Harvard University, Cambridge, MA 02138.} 
\leftline{{\it E-mail}: feehan@math.harvard.edu.}

\endchapter
\end